\documentclass[12pt]{article}

\usepackage[T1]{fontenc}
\usepackage[utf8]{inputenc}
\usepackage{amsmath, amssymb, amsthm, amsfonts}
\usepackage{graphicx}
\usepackage{epstopdf}
\usepackage{psfrag}

\usepackage{geometry}
\usepackage{mathrsfs}
\usepackage{enumitem}

\usepackage{titlesec}
\usepackage{hyperref}
\usepackage[dvipsnames]{xcolor} 

\usepackage{makeidx}
\usepackage{eucal}
\usepackage{pdfsync}
\usepackage{tikz}
\usepackage{comment}

\hypersetup{
    colorlinks=true,
    linkcolor=MidnightBlue,
    citecolor=OliveGreen,
    filecolor=magenta,
    urlcolor=MidnightBlue,
}

\newtheorem{theorem}{Theorem}
\newtheorem{hypothesis}[theorem]{Hypothesis}
\newtheorem{lemma}[theorem]{Lemma}

\newtheorem{proposition}[theorem]{Proposition}
\newtheorem{corollary}[theorem]{Corollary}
\newtheorem{remark}[theorem]{Remark}
\newtheorem{defn}[theorem]{Definition}

\def\beq{\begin{equation}}
\def\eeq{\end{equation}}

\def\bH{\mathbf{H}}
\def\bI{\mathbf{I}}

\def\bs{\mathbf{s}}

\def\cA{\mathcal{A}}
\def\cB{\mathcal{B}}

\def\cF{\mathcal{F}}
\def\cG{\mathcal{G}}
\def\cH{\mathcal{H}}
\def\cI{\mathcal{I}}
\def\cJ{\mathcal{J}}
\def\cK{\mathcal{K}}

\def\cM{\mathcal{M}}

\def\cO{\mathcal{O}}
\def\cP{\mathcal{P}}

\def\cS{\mathcal{S}}

\def\cV{\mathcal{V}}
\def\cW{\mathcal{W}}

\def\cZ{\mathcal{Z}}

\newcommand{\eps}{\varepsilon}

\def\({\left(}
\def\rt){\right)}
\def\[{\left[}
\def\]{\right]}

\def\v{{\bf v}}

\def\R{\mathbb{R}}
\def\Z{\mathbb{Z}}
\def\N{\mathbb{N}}

\def\eps{\epsilon}

\def\sign{\text{sign}}
\def\ff{\infty}

\def\({\left(}
\def\){\right)}

\def\beq{\begin{equation}}
\def\ee{\end{equation}}
\def\bea{\begin{eqnarray}}
\def\eea{\end{eqnarray}}
\def\bes{\begin{eqnarray*}}
\def\ees{\end{eqnarray*}}

\def\aaa{\mathfrak{a}}

\def\asy{\mathrm{asy}}
\def\vel{\boldsymbol{v}}
\def\nnu{\widehat{\nu}}

\numberwithin{equation}{section}

\begin{document}

\title{Limit Theorems for Dispersing Billiards \\ with Infinite Cusps}

\author{Marco Lenci$^{1}$,\quad Fran\c{c}oise P\`ene$^{2}$,\quad  Hong-Kun Zhang$^{3}$}

\date{Preliminary version \\[3pt] 
June 4, 2026}

\maketitle

\footnotetext[1]{Dipartimento di Fisica e Astronomia, Universit\`a di Bologna, Italy; Email: marco.lenci@unibo.it}

\footnotetext[2]{Univ Brest, Universit\'e de Brest, LMBA,  UMR CNRS 6205, France; Email: francoise.pene@univ-brest.fr}

\footnotetext[3]{ Department of Mathematics \& Statistics, University of Massachusetts Amherst, USA; Email: hongkun@math.umass.edu}

\begin{abstract}
Extending Lenci's work in \cite{ML02,ML03}, we present new insights into semi-dispersing billiards with infinite cusps.
We investigate the Convergence to Stable Processes (CSP) and Functional Central Limit Theorems (FCLT) both for the collision map and for the billiard flow. Whereas the billiard map always preserves an infinite measure equivalent to Lebesgue, the billiard flow 
preserves the Lebesgue measure, which is either finite or infinite depending on the finiteness of the area of the billiard table. 
When the reference measure is infinite, the usual LLN fails, but we are able to prove the CSP and FCLT. We also prove classical FCLT for the flow when the area of the billiard table is finite.  
These results are obtained through a precise study of the properties of the billiard map, combined with a general approach based on induction and on coupling involving point processes. 
\end{abstract}

\tableofcontents

\part{Introduction}

\section{Background on billiards}

A billiard is a dynamical system defined by the motion of a point particle inside a domain $Q$, known as \emph{table}, where the particle experiences inertial motion in the interior of $Q$ and totally elastic collisions at its boundaries: this means that the normal component of the velocity changes sign, while the tangential component is unchanged. For planar billiards, i.e., when $Q \subset \R^2$, this results in the celebrated Fresnel-Descartes Law: the angle of incidence equals the angle of reflection. 

A semidispersing (planar) billiard is one where $\partial Q$ is given by a finite or countable number of smooth components, which are flat (segments, half-lines or lines) or convex, as seen from the interior of the table. Semidispersing billiards were the first billiards whose chaotic behavior was studied within the context of rigorous ergodic theory. Building on the foundational work of Sinai \cite{Sin70}, it was found that, in great generality, semidispersing billiards in bounded domains are hyperbolic (in the sense of piecewise smooth systems) and Bernoulli (thus K-mixing, mixing and ergodic), cf.~\cite{GO, SC87, CH96} and references therein. Furthermore, for ample collections of such systems, the rate of decay of correlations is known for large classes of observables and is exponential when the hyperbolicity is uniform. 

Among billiard systems exhibiting nonstandard decay of correlations are those with \emph{cusps}, i.e., zero-angle corners. Chernov and Markarian~\cite{CM05} gave the first rigorous treatment of the decay of correlations in dispersing billiards with regular cusps, and Chernov and Zhang~\cite{CZ3} later established the optimal bound $O(1/n)$, where $n$ denotes time for the collision map. B\'alint, Chernov, and Dolgopyat~\cite{Balint2011,BCD2} proved further statistical properties, including a central limit theorem with the nonstandard normalization $\sqrt{n\log n}$. Zhang~\cite{Z2016b} introduced billiards with cusps at zero-curvature points, 
determining the corresponding decay rate $O(1/n^\gamma)$, $0 < \gamma < 1$. Jung and Zhang~\cite{JZ18} showed convergence to a stable law, and subsequent work~\cite{JPZ,MV,JMPVZ} established functional limit theorems yielding convergence to a Lévy process, w.r.t.\ the appropriate Skorokhod topology.

\begin{figure}[ht]
\centering
\includegraphics[width=0.5\textwidth]{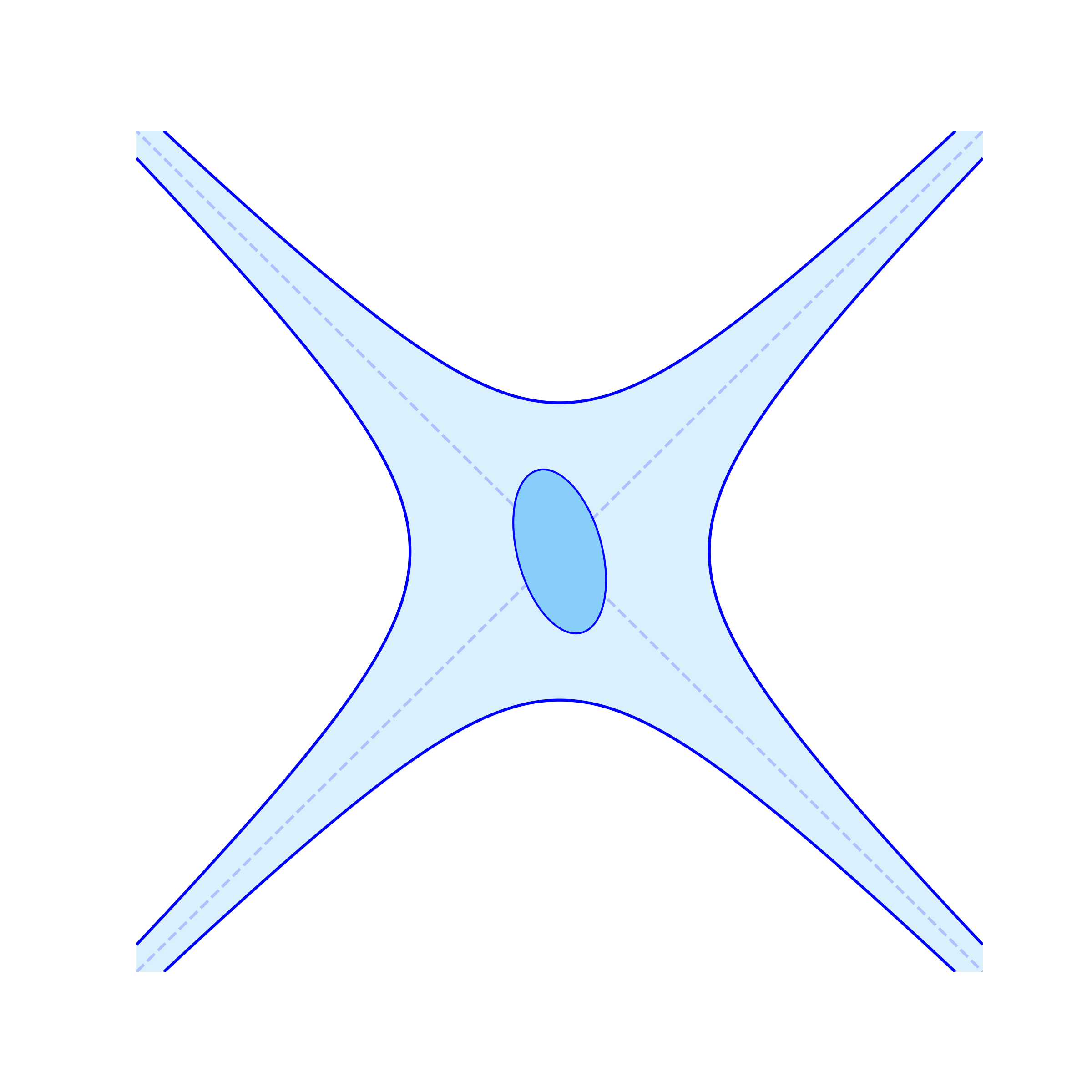}
\caption{Example of a billiard table with four infinite cusps.}
\label{fig1}
\end{figure}

The state of the art is significantly worse for semidispersing billiards with infinite cusps, where we define an infinite cusp as the part of a billiard table delimited by two infinite dispersing boundary components which lie at either side of a half-line and are asymptotic to it; see Fig.~\ref{fig1}. Lenci proved hyperbolicity and ergodicity for a number of semidispersing billiards with an infinite cusp \cite{ML02, ML03}, but little is known about higher ergodic properties, the only exception being a paper by Arbieto, Markarian, Pacifico and Soares~\cite{AMPS} showing that the billiards studied in \cite{ML03} are mixing in the sense of Krengel and Sucheston \cite{KS1985}, i.e., of \emph{zero type} \cite{DS2023}, a.k.a.\ \emph{local-local mixing} \cite{ML13}.

This unsatisfactory state of affairs is due to the fact that billiards in non-compact tables fully belong to the realm of infinite ergodic theory \cite{Aaronson1997, KMS2016, DS2023}, a much more intricate and less developed subject than ordinary (finite) ergodic theory. Even basic definitions and elementary results of finite ergodic theory are problematic and/or false \cite{Krengel1985, Petersen1983, Aaronson1997, ML10}, let alone the (few) standard techniques for the stochastic properties of dynamical systems, which, naturally, rely heavily on the relevant measure being a probability. Essentially, at the moment, the only method to study the stochastic properties of hyperbolic or expanding dynamical systems in infinite measure (outside the particular case of group extensions) is to consider the induced dynamics on a suitable finite-measure subset. This, however, requires a tight control on the various quantities involved in the statistical description of the dynamics (e.g., return times, distortion bounds, decorrelation rates, etc.). The present paper corroborates this ``empirical principle'' to an especially strong degree. 

That notwithstanding, a quite significant amount of literature was produced in the last decades on a natural and physically relevant class of unbounded dispersing billiards, the Lorentz gas. Here we refer to both periodic and aperiodic systems. The former are of course easier to study, as they can be thought of as extensions of Lorentz gases on the torus, a.k.a.\ \emph{Sinai billiards}; a nice survey of classical and more recent results on Sinai billiards and periodic (infinite) Lorentz gases is found in \cite{FP25}. We mention in particular two recent papers by Pène and Thomine \cite{PT1, PT2}, in which central limit theorems (in the strong distributional sense) are proved for the first time for certain classes of non-periodic observables, thus showing specific nonstandard statistical properties of infinite-measure Lorentz gases. See also the very recent \cite{FP23}.

As for aperiodic Lorentz gases, results are understandably weaker but still very much nontrivial. An admittedly incomplete list of references includes \cite{ML06, CLS2010, DSV2009, LenciTroubetzkoy, DolgopyatNandori2018, AiminoLiverani, GalliLenci, TZ23}.

In this paper we present a thorough study of the statistical properties of a family of dispersing billiards with a finite number of infinite cusps, such as the one illustrated in Fig.~\ref{fig1}. In more detail, we assume that the boundary of the billiard table $Q$ consists of finitely many disjoint $C^3$-smooth curves, a.k.a.\ \emph{boundary components}, whose curvature is never zero. Some of the boundary components might correspond to internal scatterers, given by bounded convex $C^3$ sets. These assumptions imply in particular that $\partial Q$ has no corners or finite cusps. As for the infinite cusps, we assume the following:

\begin{itemize}
\item[\textbf{(h1)}] 
There are $\cJ \in \Z^+$ infinite cusps and there exist positive constants $\beta, \widehat{t}, \mathfrak{c}_i$ ($i = 1, \ldots, \cJ$) such that the $i$-th cusp is delimited by the two curves $\Gamma^{(i)}_{\pm}$. In a suitable Cartesian coordinate system $(t, z)$ (depending on $i$), $\Gamma^{(i)}_{\pm}$ are given by the equations $z = \pm f_{\mathfrak{c}_i}(t)$, $t \ge \widehat{t}$, where
\beq\label{z1s}
f_{\mathfrak{c}}(t) = \frac{\mathfrak{c}}{\beta t^{\beta}} \, .
\eeq  
We call any right half-line of the $t$-axis an \emph{asymptote} for the cusp.%
\footnote{Here and in the rest of the paper, when we refer to a cusp we always mean the asymptotic part of a cusp, that is, in $(t, z)$ coordinates, the region
$-f_{\mathfrak{c}_i}(t) \le z \le f_{\mathfrak{c}_i}(t)$, for $t$ large enough depending on the context.}

\item[\textbf{(h2)}] For any $i=1,...,\cJ$, for the $i$-th cusp, there exists a boundary point $q^{(i)}_\asy$ and a unit velocity $\vel^{(i)}_\asy$ such that the forward trajectory with initial conditions $x^{(i)}_\asy = (q^{(i)}_\asy, \vel^{(i)}_\asy) \in \partial Q \times \mathbb{S}^1$ is an asymptote for the cusp. We assume that $\vel_\asy^{(i)}$ (equivalently, the $t$-axis) is not tangent to $\partial Q$ at $q_\asy^{(i)}$.
\end{itemize}

Thus, all cusps are \emph{symmetric} and \emph{similar}, the latter term meaning that they coincide up to a dilation in the direction of the asymptote. We point out that our main results, cf.~Section \ref{1_Sec:2}, remain true \emph{verbatim} even if $Q$ possesses cusps defined by an exponent $\beta' < \beta$.%
\footnote{Specifically, the individual study of any such cusp does not change, but its contribution will vanish in the tail probability estimate given in Proposition~\ref{Mm2} and thus in Condition~(A) of Hypothesis~\ref{hypinductionbis}.}
As for the symmetry assumption, we fully expect that a cusp defined, in the suitable Cartesian system,  by $z = f_{\mathfrak{c}^+}(t)$ and $z = -f_{\mathfrak{c}^-}(t)$, $\mathfrak{c}^+ \ne \mathfrak{c}^-$, contributes to our estimates like a symmetric cusp, up to inessential  constants. These assumptions were made for the pure sake of technical simplicity, as was the case for the absence of corner points. 

What is really relevant for the statistical properties of our billiards is the parameter $\beta$, which we also refer to as the \emph{order} of the infinite cusps defined above.

\section{Limit theorems for the billiard map and flow} \label{1_Sec:2}

There are traditionally two ways to study the dynamics of a billiard system: in terms of physical time, thus using the so-called \emph{billiard flow} $(\Phi^t)_{t \in \R}$, or by looking at the system only at collision times (as the trajectory between a collision and the next is trivial for billiards): this amounts to using the cross section for $\Phi^t$ defined by all post-collisional position-velocity pairs. The resulting map is called the \emph{billiard map} or the \emph{collision map}. As is customary, the speed of the point particle is assumed to be 1. The conservation of energy ensures that this assumption does not result in a loss of generality.

The billiard flow $\Phi^t$ is thus defined on $Q \times \mathbb{S}^1$, interpreted as the space of all admissible position-velocity pairs $(q,\vel)$. An implicit way to define it is to say that it describes inertial motion on $Q \times \mathbb{S}^1$, subject to the identifications $(q,\vel_-) \cong (q,\vel_+)$, where $q\in \partial Q$ and $\vel_-, \vel_+$ is a pair of pre- and post-collisional unit vectors, relative to the boundary collisions occuring in $q$. It is a well-known fact, amounting to Liouville's Theorem (since billiards are non-smooth Hamiltonian systems, after all) that $\Phi^t$ preserves $\nu$, the Lebesgue measure on $Q \times \mathbb{S}^1$.%
\footnote{More precisely, $\nu$ is the product of the Lebesgue measure on $Q$ and the Haar measure on $\mathbb{S}^1$. Observe that $Q$ may have finite or infinite area, thus $\nu$ can be finite or infinite.}

We denote $\cM \subset Q \times \mathbb{S}^1$ the set of all post-collisional pairs $(q,\vel)$, meaning that $q \in \partial Q$ and $\vel$ points towards the interior of $Q$. Alternatively, if the identification of pre- and post-collisional velocities is implemented, $\cM$ may be defined as $\partial Q \times \mathbb{S}^1$, modulo such identification. A convenient way to represent an element $x  = (q,\vel) \in \cM$ is by means of the pair $(r, \varphi)$, where $r$ is a clockwise arc-length coordinate identifying $q \in \partial Q$, and $\varphi \in [-\pi/2, \pi/2]$ is the angle that $\vel$ makes with the inward normal to the boundary at $q$. With an understandable abuse of notation, we shall write $x = (r, \varphi)$. $\cM$ is called the \emph{collision space}. It is clearly a global cross section for the billiard flow, namely, all forward trajectories of the flow, save for the ones running on the asymptotes of the cusps (a null-measure set) intersect $\cM$. The billiard map $\cF: \cM \to \cM$ is then defined as the Poincar\'e map of the billiard flow on $\cM$. A consequence of the invariance of $\nu$ for the flow is that $\cF$ preserves the measure $\mu$, which, in $(r, \varphi)$ coordinates, is defined by \cite{CM}
$$
  d\mu = \cos\varphi \, dr \, d\varphi\, .
$$
Observe that, for our billiards, $\mu$ is always an infinite, $\sigma$-finite measure. We now present our main results, first for the collision map and then for the flow.

\medskip

Given an observable $f \in L^1(\cM, \mu)$, we write $(S_n f)_{n\in\N}$ for the sequence of Birkhoff sums with respect to the billiard map, i.e.,
\begin{equation} \label{defSnf}
  S_n f(x) := \sum_{k=0}^{n-1} f \circ \cF^k(x) \ \ (n\ge 1) \,, \qquad S_0 f \equiv 0 \,.
\end{equation}
One is interested in the asymptotic behavior of $(S_n f)_n$, w.r.t.\ the infinite measure $\mu$. It is a general fact that, for chaotic dynamical systems preserving an infinite measure, the classical Law of Large Numbers (LLN) fails. Instead, integrable observables satisfy a \emph{degenerate LLN}:

\begin{theorem}[Billiard map: Degenerate Law of Large Numbers] \label{LLNfails}
The infinite-measure-preserving dynamical system $(\cM, \cF, \mu)$ is (Poincar\'e) recurrent and ergodic. In particular, for all $f \in L^1(\cM, \mu)$,
$$
\lim_{n \to \infty} \frac{S_n f} n = 0 \,, \qquad \mbox{$\mu$-almost everywhere}.
$$
\end{theorem}

This result is proved in Section~\ref{sec44}. 
The degeneracy of the above limit naturally prompts one to look for a non-degenerate limit theorem for $S_n f$, with a different normalization. We do prove such a theorem, showing \emph{strong distributional convergence} to an explicit limit (Theorem \ref{thm:LLN}). Strong distributional convergence, or \emph{strong convergence in distribution}, is the natural notion of distributional convergence relative to an infinite measure, although it can be defined both for finite and infinite measures, as we recall here.

\begin{defn} \label{def-sid}
If $(Z_n)_{n\in\N}$ is a sequence of measurable functions on a (finite or infinite) measure space, we say that $Z_n$ converges to the random variable $Z$ {\bfseries strongly in distribution}, as $n \to \infty$, if $Z_n$ converges in distribution to $Z$ relative to any probability measure absolutely continuous w.r.t.\ the reference measure. We shall then write $Z_n \Longrightarrow Z$.
\end{defn}

Let us introduce a stochastic process that will appear in our upcoming limit theorem. Set 
\beq \label{alpha}
  \alpha := \frac{\beta}{\beta+1}
\eeq 
and, for $b>0$,
\beq \label{formulaI}
  I_{b} := \int_0^{\frac{\pi}{2}} (\sin \varphi)^{b} \, d\varphi 
\eeq
(this is an important quantity that will be used throughout the paper). Also, recalling that, for each cusp, $\mathfrak{c}_i$ is the constant given in \textbf{(h1)}, denote 
\begin{equation}\label{Calpha}
  \mathfrak{C} := \sum_{i=1}^{\cJ} \mathfrak{c}_i^{\frac 1{\beta+1}} \,.
\end{equation}
We denote by $(\Delta_1(s))_{s \ge 0}$ a positive $\alpha$-stable L\'evy process (with independent stationary increments) such that $\Delta_1(1)$ has characteristic function $\Psi_{\alpha,c_1}$, where
\beq \label{stable1}
\Psi_{\alpha,c}(u) := \exp\!\left( -|cu|^\alpha \, \Gamma(1-\alpha)e^{-i\frac{\pi\alpha}{2} 
\sign(u)} \right)
\eeq 
and
\begin{equation}\label{:c1}
c_1 := \left( \frac{4\mathfrak{C}} \beta \right)^{\frac{1}{\alpha}} I_{1/\alpha} \,,
\end{equation}
cf.~\eqref{formulaI}-\eqref{Calpha}. Recall that $\Delta_1(s)$ has the same distribution as $s^{\frac{1}{\alpha}} \Delta_1(1)$. We shall especially need the generalized inverse process $(\Delta_1^{-1}(t))_{t\ge 0}$, defined by
$$
\Delta_1^{-1}(t) := \inf \{s > 0 \, : \, \Delta_1(s) \ge t\} \, .
$$

\begin{theorem}[Billiard Map: Convergence to Stable Process (CSP)] \label{thm:LLN} 
For all $f \in L^1(\cM, \mu)$ and $T>0$, the following limit holds, as $n \to \infty$:
$$
\left( \frac{S_{\lfloor nt \rfloor} f} {n^\alpha} \right)_{t \in [0,T]} \Longrightarrow \left( \left( \int_{\cM} f \, d\mu \right) \Delta_1^{-1}(t)\right)_{t\in[0,T]}\, ,
$$
denoting strong distributional convergence (see Definition \ref{def-sid}) of the above processes, in the uniform topology, w.r.t.\ $\mu$ (equivalently, the Lebesgue measure) on $\cM$.
\end{theorem}

An easy corollary is the strong distributional convergence of the corresponding random variables:
$$
\frac{S_n f}{n^\alpha} \Longrightarrow \Delta_1^{-1}(1) \int_{\cM} f \, d\mu \, ,
$$
as $n \to \infty$. Theorem \ref{thm:LLN} is proved at the end of Section~\ref{proofLLN}. The above limits vanish when $\int f \, d\mu = 0$, that is, when $f$ is \emph{$\mu$-centered}. One naturally wonders what the growth rate is for Birkhoff sums of centered observables, and if a limit exists upon the corresponding rescaling. Such a result may be rightfully be called a Central Limit Theorem for our infinite-measure-preserving billiard map. We prove it in its functional form for the natural case of compactly supported H\"older observables.

\begin{theorem}[Billiard Map: Functional Central Limit Theorem (FCLT)] \label{thm:TCL}
Let $f:\cM \rightarrow \mathbb{R}$ be a $\mu$-centered H\"older function with compact support. There exists $\sigma_f \ge  0$ such that, for any $T > 0$,
$$
\left(\frac{S_{\lfloor nt \rfloor} f }{n^{\frac{\alpha}{2}}}\right)_{t \in [0,T]} \Longrightarrow \left(\sigma_f B_{\Delta_1^{-1}(t)}\right)_{t \in [0,T]} \, ,
$$
as $n \to \infty$, in the uniform topology, w.r.t.~$\mu$ (or Lebesgue), where $B$ is a standard Brownian motion, independent of $\Delta_1$. 
\end{theorem}

We present the proof of Theorem \ref{thm:TCL} at the end of Section~\ref{proofTCLmap}.
The asymptotic variance $\sigma^2_f$ appearing in the limit above is given by the following Green-Kubo formula
\begin{equation} \label{diffusioncoef}
\sigma^2_f := \sum_{k \in \mathbb{Z}} \int_{M} f \cdot f \circ F^k \, d\mu \, ,
\end{equation}
where $M$ is some naturally chosen bounded subset of $\cM$, containing the support of $f$, and $F := \cF_M$ is the first-return map of $\cF$ to $M$ (see Section~\ref{Sec:3} for the precise assumptions on $M$).

\medskip\noindent\textbf{Remark (Non-degeneracy of the diffusion coefficient).}
Because the reference set $M$ contains the support of $f$, the induced observable $\widetilde f:=\sum_{k=0}^{R-1}f\circ\cF^{k}$ coincides with $f$ on $M$; hence \eqref{diffusioncoef} is, up to the factor $\mu(M)$, the Green--Kubo variance of $f$ for the \emph{probability}-preserving, exponentially mixing induced map $(M,F,\mu_M)$, cf.~\eqref{correT}. The series therefore converges absolutely and $\sigma_f\ge 0$. By the classical coboundary dichotomy for invertible systems with summable correlations,
$$
\sigma_f=0\quad\Longleftrightarrow\quad \text{$f$ is an $L^2(\mu_M)$ coboundary for $F$,}
$$
i.e.\ $f=\psi\circ F-\psi$ holds $\mu_M$-a.e.\ for some $\psi\in L^2(M,\mu_M)$ (equivalently, the Birkhoff sums of $f$ under $F$ stay bounded in $L^2(\mu_M)$); see e.g.~\cite[Ch.~7]{CM}. Consequently the limit in Theorem~\ref{thm:TCL} is non-degenerate for every $f$ outside this exceptional cohomology class. A convenient sufficient condition for $\sigma_f>0$ is the existence of two $F$-periodic points $p,q$ of equal period $m$ with $\sum_{j=0}^{m-1}f(F^{j}p)\neq\sum_{j=0}^{m-1}f(F^{j}q)$, provided the transfer function admits a continuous version (as Liv\v{s}ic-type regularity suggests for the uniformly hyperbolic $F$).
\par\medskip

Let us turn to the billiard flow $(\Phi^t)_{t \in \R}$. Observe via \eqref{z1s} that the billiard table $Q$ has finite area if and only if $\beta>1$. Thus, $(Q \times \mathbb{S}^1, (\Phi^t)_t, \nu)$ is a finite- or infinite-measure-preserving dynamical system depending on $\beta$ being larger than 1 or otherwise. Naturally, the limit theorems will vary accordingly. When $\beta>1$, let us denote by $\nnu$ the normalized version of $\nu$:
$$
\nnu := \frac{\nu} {\nu(Q \times \mathbb{S}^1)} = \frac{\nu} {2\pi \, \mathrm{Area}(Q)} \, . 
$$

\begin{theorem}[Billiard flow in finite area: Standard FCLT, a.k.a.\ Weak Invariance Principle] \label{thm:TCLflow}
Assume $\beta>1$. The billiard flow $\Phi^t$ with the measure $\nnu$ (equivalently, $\nu$) is ergodic, that is, for all $f \in L^1 (Q \times \mathbb{S}^1, \nnu)$,   
$$
\lim_{T \to +\infty} \frac 1 T\int_0^T f \circ \Phi^s\, ds = \int_{Q \times \mathbb{S}^1} f\, d\nnu,\quad \nnu\mbox{-almost everywhere} \,.
$$
Moreover, let $g: Q \times \mathbb{S}^1 \rightarrow \mathbb{R}$ be H\"older continuous, compactly supported and $\nnu$-centered (i.e., $\int g\, d\nnu=0$). There exists $\boldsymbol{\sigma}_g \ge  0$ such that, for all $T > 0$,
$$
\left( \frac 1 {\sqrt{n}} \int_0^{nt} g \circ \Phi^s \, ds \right)_{t \in [0,T]} \Longrightarrow \left(\boldsymbol{\sigma}_g B_t \right)_{t \in [0,T]} ,
$$
as $n\to \infty$, in the uniform topology, w.r.t.~$\nnu$, where $B$ is a standard Brownian motion. 
\end{theorem}

This result is proved in Section~\ref{:proofThm5}. Without giving technical details at this point, we anticipate that the diffusion coefficient $\boldsymbol{\sigma}_g$ can be expressed via some Green-Kubo formula, as was to be expected, cf.~Remark~\ref{remsigmagflow}. We now turn to the infinite-area case. In order to avoid an unpleasant number of technical complications, we exclude the case $\beta=1$.

\begin{theorem}[Billiard flow in infinite area: CSP and FCLT] \label{thm:TCLflow2} 
Assume $0<\beta< 1$. For all $f \in L^1 (Q \times \mathbb{S}^1, \nu)$ and $T>0$, 
$$
\left( \frac 1 {n^\beta} \int_0^{nt} f \circ \Phi^s\, ds \right)_{t\in[0,T]} \Longrightarrow \left( \left( \int_{Q\times\mathbb S^1} f\, d\nu \right) \Delta_2^{-1}(t)\right)_{t\in[0,T]}
$$
as $n\to \infty$, in the uniform topology, w.r.t.~$\nu$. Here $\Delta_2$ is an increasing $\beta$-stable L\'evy process (with independent stationary increments) such that the characteristic function of $\Delta_2(1)$ is $\Psi_{\beta,c_2}$, as defined in \eqref{stable1}, with
$$
c_2 := \left(2\sum_{i=1}^\cJ \mathfrak c_i \right)^{\frac{1}{\beta}} \left( \frac 2 \beta \right)^{\frac{\beta+1}{\beta}} I_{\frac{1}{\beta}-1} \, ,
$$
where $I_{\frac{1}{\beta}-1}$ is given by \eqref{formulaI} (recall the meaning of $\mathfrak c_i$ and $\cJ$ in \textbf{(h1)}).

Also, let $g: Q \times \mathbb{S}^1 \rightarrow \mathbb{R}$ be a $\nu$-centered, compactly supported, H\"older continuous function. There exists $\boldsymbol{\sigma}_g \ge 0$ such that, for all $T > 0$,
$$
\left( \frac1 {n^{\beta/2}} \int_0^{nt} g \circ \Phi^s \, ds\right)_{t \in [0,T]} \Longrightarrow \left( \boldsymbol{\sigma}_g B_{\Delta_2^{-1}(t)} \right)_{t \in [0,T]} \, ,
$$
as $n\to \infty$, in the uniform topology, w.r.t.~$\nu$, where $B$ is a standard Brownian motion independent of $\Delta_2$.
\end{theorem}

Theorem \ref{thm:TCLflow2} is proved in Section~\ref{:proofThm6}. Provided the support of $g$ is made of points (i.e., position-velocity pairs) with previous collision in $M$ ($M$ was introduced after the statement of Theorem \ref{thm:TCL}), the quantity $\boldsymbol{\sigma}_g$ of Theorem~\ref{thm:TCLflow2} corresponds to $\sigma_{U_g}$, defined in \eqref{diffusioncoef}, with
$$
U_g(x) := \int_0^{\tau(x)}g(\Phi^s(x))\, ds\, ,
$$
where $\tau(x)$ is the free flight
$$
\tau(q,\vel) := \inf \{s>0 \, :\, q+s\vel \in \partial Q\}\, .
$$
This identification $\boldsymbol{\sigma}_g=\sigma_{U_g}$ is justified in Section~\ref{proofTCLnonstandardflow}: the support condition guarantees that, on each excursion, $g$ is integrated only along the single free flight issued from $M$, so that the suspension/time-change argument reduces the flow Green--Kubo sum to that of $U_g$ for the induced map $F$.

These four limit theorems are proved in a general framework in Sections~\ref{proofformap} and~\ref{proofTCLflow}. The fact that the billiards we are considering satisfy this general framework rely on estimates stated in Sections~\ref{cornerseries},~\ref{Sec:3} and~\ref{secF}, completed by Appendices~\ref{proofprop3} and~\ref{appendKBbd}. More precisely, general versions of Theorems~\ref{thm:LLN} and~\ref{thm:TCL} (limit theorems for the map) are stated and proved respectively in Section~\ref{proofLLN} and~\ref{proofTCLmap}, whereas general versions of Theorems~\ref{thm:TCLflow} and~\ref{thm:TCLflow2} (limit theorems for the flow) are proved in Section~\ref{proofstandardTCLflow} and~\ref{proofTCLnonstandardflow}. The proofs rely on results stated in Section~\ref{secF} for a general context, and use induction and a convergence result for the pair made by a point process and the (continuous-argument) process associated to a Birkhoff sum. 
\medskip

\part{Properties of billiards with infinite cusps}

\section{The Corner Series}\label{cornerseries}

Building upon the methods developed by Young \cite{Y98}, a comprehensive framework for analyzing the statistical properties of nonuniformly hyperbolic systems was formulated by Markarian, Chernov, and Zhang \cite{M04, CZ, CZ2, CZ3, CZ09}. This framework was used to establish slow rates of mixing for hyperbolic systems with singularities and was applied to various models.

We fix $\widehat{t} \in \R^+$ and consider the subdomain $\widetilde{Q} = Q\setminus \bigcup_{i=1}^{\cJ}\bar Q_i$, where $\bar Q_i$ is the $i$-th cusp area corresponding to 
$$ \{(t,z)\, :\, |z| \le  f_{{\mathfrak{c}}_i}(t), \ t > \widehat{t}\}\, ,$$
in the $( t,z)$ coordinates of assumption {\bf (h1)}. Assuming $\widehat{t}$ is sufficiently large, all trajectories exiting a cusp area will hit $\widetilde{Q}$ before returning to a cusp area. We define the induced domain as the set
$$M := \partial \widetilde{Q} \times [-\pi/2, \pi/2].$$
Note that $M$ is a bounded set, thus $\mu(M) < \infty$. For any $x \in M$, we define the first return time function as
$$R(x) = \min\{n \geq 1: \cF^n(x) \in M\},$$
and thus the first return map $F \colon M \to M$ is defined by
\beq \label{Fdef}
F(x) = \cF^{R(x)}(x), \quad \forall x \in M.
\eeq
The induced map $F$ preserves a probability measure $\mu_M$ obtained by conditioning $\mu$ on $M$:
$$d\mu_{M} = c_M \cos \varphi \, dr \, d\varphi,$$
where $c_M := (2|\partial \widetilde{Q}|)^{-1}$ is a normalizing constant.

In this section, following \cite{CZ3, Z2016b}, we investigate the geometry of the corner series in the cusps. We fix some $i=1,...,\cJ$. 
Let $M^{(i)}$ be the subset of $M$ containing points first entering the $i$-th cusp.  We define $M^{(i)}_N := M^{(i)}\cap \{R=N\}$, representing billiard trajectories entering the $i$-th cusp and experiencing $(N-1)$ reflections before exiting.
We indicate that all constants chosen in this section, unless otherwise stated, depend neither on $N$ nor on the choice of $x \in M_N^{(i)}$, nor on $n\in\{1,...,N-1\}$.

For any $x = (r, \varphi) \in M^{(i)}_N$, we denote $x_N = F (x)$. Note that $x_N = \cF^N (x)$, and denote $x_n := \cF^n (x) = (r_n, \varphi_n)$ for $n = 1, \ldots, N-1$, representing the orbit in the cusp. We call $\{x_n, n = 1, \ldots, N-1\}$ a corner series of length $N-1$. We denote $t_n$ as the $t$-coordinate of the base point of $x_n = (r_n, \varphi_n)$, with the notation of~{\bf (h1)}.
Observe that
\beq
\label{rntn}
r_n = \cO(1) + \int_{\widehat t}^{t_n} \sqrt{1 + |f_{\mathfrak{c}}'(t)|^2} \, dt = \cO(1) +\int_{\widehat t}^{t_n} (1 + \cO(t^{-2\beta-2})) \, dt = t_n + \cO(1)\, .
\eeq

We also set $\tau_n=\tau(x_n)$ as the distance between the respective positions of $x_n$ and $x_{n+1}$, and $\mathcal{K}_n$ as the curvature of $\partial Q$ at the position of $x_n$. 
We also use more convenient coordinates:
$$\gamma_n := \pi/2 - |\varphi_n|, \quad \alpha_n := \arctan(|f'_{{\mathfrak{c}}_i}(t_n)|) \, ,$$
where $\alpha_n$ is the angle of the tangent line to $\Gamma^{(i)}_{\pm}$ at $(t_n, \pm f_{{\mathfrak{c}}_i}(t_n))$ made with the horizontal axis, and $\gamma_n$ is the angle between this tangent line and the post-collisional vector at this point.
We also set $v_n$ for the angle made between the velocity vector (after the $n$-th reflection)
with the horizontal line in $(t,z)$ coordinates:
\begin{equation}\label{defvn}
v_n:=\gamma_n+\alpha_n\, .
\end{equation}
Note that $\alpha_n$ and $\gamma_n$ are positive for $2 \leq n \leq N-1$; $\alpha_n$ are all small if $N$ is large enough. While $\gamma_n$ are initially small, they slowly grow to about $\pi/2$ for $n \sim N/2$, and then decrease again, becoming small once more. As in~\cite{CM05}, we denote by $N_2$ the turning time given by
$$
t_{N_2+1} := \max\{t_n : 1 \leq n \leq N\}.
$$
This means that
\begin{equation}\label{gammaN2}
v_{N_2}
=\gamma_{N_2} + \alpha_{N_2} \le \frac{\pi}{2} < v_{N_2}+2\alpha_{1+N_2}.
\end{equation}
Using the symmetry property of the billiard table boundary in the cusp area, we can show that the number of collisions in the corner series occurring before $N_2$ and after $N_2$ differ by no more than one. That is,
\begin{equation}\label{N2}
|N_2 - N/2| \leq 2,
\end{equation}
with a proof identical to the one given in \cite[Lemma 3.1]{CM05} (we will not repeat this proof here). Thus, we further subdivide the corner series into three segments. 
We fix a very small number $\bar{\gamma}$ and define
\beq\label{N1N3}
N_1 = \max\{n \leq N_2 : \gamma_n < \bar{\gamma}\}, \quad N_3 = \max\{n \geq N_2 : \gamma_n > \bar{\gamma}\}.
\eeq
We call the segment on $[1, N_1]$ the "entering period" in the corner series, the segment $[N_1 + 1, N_3 - 1]$ the "turning period", and the segment $[N_3, N]$ the "exiting period". Clearly, $|N_3 + N_1 - N| \leq 2$.

Recall that the billiard map is time-reversible. Indeed $\mathcal{F}^{-1} = \sigma_0 \circ \mathcal{F} \circ \sigma_0$ where $\sigma_0$ is the involution given by $\sigma_0(r, \varphi) = (r, -\varphi)$ in coordinates $(r,\varphi)$. Thus it is sufficient to consider the first half of the series $x_n$, $1 \leq n \leq N_2$. 
The next result is proved in Appendix~\ref{proofprop3}. 
\begin{proposition}\label{prop3} 
Let $i=1,...,\cJ$. 
For every $N$
large enough and any 
$x\in M^{(i)}_N$, 
the following holds true uniformly:\footnote{
the notation $a_{n,N}(x)\approx b_{n,N}(x)$ means that there exist positive constants $c'_1,c'_2$ 
such that, for every $N$ large enough, for every $x\in M_N$, for every $n$,
$c'_1 a_{n,N}\le b_{n,N}\le c_2a_{n,N}$.}
\begin{itemize}
\item[(1)]  $N_1\approx N_2-N_1\approx N_3-N_2\approx N-N_3\approx N.$ i.e. all  three segments in the corner series have length of order $N$;
\item[(2)] $\alpha_n
\approx (n^{\beta+1}N^{ \beta})^{-\frac 1{2\beta+1}}$, for $n\in [1,N_2]$, and so in particular  $\alpha_1\approx N^{-\frac{\beta}{2\beta+1}}$ 
and $\alpha_n\approx n^{-1}\approx N^{-1}$ for $n\in[N_1,N_2]$ ; furthermore $\alpha_{n}/\alpha_{n+1}-1=\cO(n^{-1})$;
\item[(3)] $v_1\approx N^{-\frac\beta{2\beta+1}}$; 
$\gamma_1= \cO\left(N^{-\frac\beta{2\beta+1}}\right)$,
$\gamma_n 
\approx (nN^{-1})^{\frac\beta{2\beta+1}}\approx n\alpha_n
 $ for $n\in [2,N_2]$, and so in particular $\gamma_2\approx N^{-\frac {\beta}{2\beta+1}}$;
\item[(4)] $t_n\approx  (nN^{\frac \beta{\beta+1}})^{\frac 1{2\beta+1}}$ 
for $n\in [1,N_1]$, and so in particular $t_1\approx N^{\frac{\beta}{(2\beta+1)(\beta+1)}}$ and $t_{n}\approx N^{\frac{1}{\beta+1}}$ for $n\in [N_1, N_2]$; furthermore 
$\tau_n\approx t_n/n$ for $n\in [1,N_1]$;
\item[(5)] $\frac{2\cK_n\tau_n}{\sin\gamma_n}
=\frac{(\beta+1)\beta}{(2\beta+1)^2 n^2}\left(1+\cO\left(
\frac{
1}n
+v_n^2
\right)\right)$ for $n\in [2,N_1]$;
\item[(6)]$\frac{\tau_{n+1}}{\tau_{n}}=1-\frac{2\beta}{(2\beta+1)n}\left(1+ \cO\left(\frac{
1}n+v_n^2\right)\right)
$, for $n\in [1,N_1]$;
\item[(7)]
  The  quantities  $H_n=t_n^{-\beta}\sin v_n$ (for $n=1,\cdots, N-1$) are almost equal and invariant along a corner series of length $N$:
\begin{equation}\label{controlHn}
\forall n=1,...,N_2,\quad H_n=H_{N_2}(1+\cO(n^{-1}))=H_{N_2}+\cO(t_n^{-2\beta-1})\, ,
\end{equation}
and
\begin{equation}\label{HN2}
H_{N_2}=
(I_{1/\alpha}/(\mathfrak c_iN))^{\frac\beta{\beta+1}}
\left(1+\cO(N^{-1})\right)\, ,
\end{equation}
with $I_{1/\alpha}:=\int_{0}^{\frac \pi 2
}(\sin v)^{1+\frac 1\beta}\, dv$ (see~\eqref{alpha}--\eqref{formulaI}).
This is similar to the Clairaut relation for geodesic on revolution surface.
 \end{itemize}
\end{proposition}

\section{Properties of the induced dynamical system }\label{Sec:3}

\begin{remark}\label{rktower}
From our detailed study of trajectories in cusps outlined in Proposition~\ref{prop3}, it follows that any long sequence of iterations based on a cusp will eventually exit the cusp and thus visit $M$. Therefore, the infinite measure system $(\mathcal M, \mathcal F, \mu)$ can be modeled by a tower system over $(M, F, \mu(M)\mu_M)$, by identifying $(x, \ell) \in M \times \mathbb{N}_0$ with $\mathcal{F}^\ell(x) \in \mathcal{M}$. Since $\mu_M(M)=1$, due to the Poincaré theorem applied to $(M,F,\mu_M)$, this already ensures the recurrence of $(\mathcal M, \mathcal F, \mu)$. Furthermore, the ergodicity of $(\mathcal M, \mathcal F, \mu)$ will follow from that of $(M, F, \mu_M)$.
\end{remark}

Our study of long cusp trajectories also implies the following estimate on the return time $\mathcal{R}$ to $M$. This proposition will be crucial for our proof of limit theorems based on convergence of point processes (see Part~\ref{part3}).
\begin{proposition}\label{Mm2}
For any $i =1,..., \cJ$ and any $y > 1$, we have
\begin{equation}\label{mes1}
\mu_M(M^{(i)} \cap \{R \ge y\}) \sim \frac{4 I_{1/\alpha}^\alpha\mathfrak{c}_i^{\frac 1{\beta+1}}}{\beta \mu(M)}y^{-\alpha},
\end{equation}
where $\mathfrak{c}_i$ is defined in \textbf{(h1)} for the $i$-th cusp and where $I_{1/\alpha}$ is given by~\eqref{formulaI}. In particular,
\begin{equation}
\mu_M(M \cap \{R \ge y\}) \sim \frac{4 \mathfrak{C}}{\beta \mu(M)} \frac{(I_{1/\alpha})^\alpha}{y^\alpha},
\end{equation}
where $\mathfrak{C}$ was defined in \eqref{Calpha}.
\end{proposition}

By the $\mathcal{F}$-invariance of $\mu$, the left hand side
of~\eqref{mes1} is equal to $\mu\left( \bigcup_{N \ge y} \mathcal{F}^{N_2}(M_N^{(i)})\right)/\mu(M)$, where $N_2$ is the turning time as defined before Proposition~\ref{prop3}.
The measure of this set will be approximated with the use of next lemma.

\begin{lemma}\label{L0}
Let $i=1,...,\cJ$. 
The set $\bigcup_{N \ge y} \mathcal{F}^{N_2}(M_N^{(i)})$ corresponds to the set of configurations parametrized by $(r, \varphi)$ based on $\Gamma_{\mathfrak{s}}^{(i)}$ (with $\mathfrak{s} = \pm$) such that
\begin{equation}\label{bornphi}
{\mathfrak{c}}_i r^{-\beta-1} (1 + \mathcal{O}(r^{-1})) \le \mathfrak{s} \varphi < 3 {\mathfrak{c}}_i r^{-\beta-1} (1 + \mathcal{O}(r^{-1})),
\end{equation}
and
\begin{equation}\label{Cond2}
r\ge \left(\frac{{\mathfrak{c}}_i y}{I_{1/\alpha}}\right)^{\frac{1}{\beta+1}} \left(1 + \mathcal{O}\left(y^{-\frac{1}{\beta+1}}\right)\right).
\end{equation}
\end{lemma}
\begin{proof}

The fact that $(r, \varphi)$ is a turning point $(r_{N_2}, \varphi_{N_2})$ (with the notations of Proposition~\ref{prop3}) means that 
$$
\alpha_{N_2} + \frac{\pi}{2} - \mathfrak{s} \varphi = v_{N_2} \le \frac{\pi}{2} < 
\alpha_{N_2} + \frac{\pi}{2} - \mathfrak{s} \varphi + 2 \alpha_{1+N_2},
$$
(see~\eqref{gammaN2}), i.e., that
$$
\alpha_{N_2} \le \mathfrak{s} \varphi < \alpha_{N_2} + 2 \alpha_{1+N_2}.
$$
It follows from~\eqref{rntn} that $t_{n} = r_{n} (1 + \mathcal{O}(r_{n}^{-1}))$ and thus 
\begin{align*}
\alpha_n &= \arctan|f'_{{\mathfrak{c}}_i}(t_{n})| = \arctan({\mathfrak{c}}_i t_{n}^{-\beta-1}) = \arctan\left({\mathfrak{c}}_i |r_{n}|^{-\beta-1} (1 + \mathcal{O}(|r_{n}|^{-1}))\right)\\
&= {\mathfrak{c}}_i |r_{n}|^{-\beta-1} (1 + \mathcal{O}(|r_{n}|^{-1}))\, .
\end{align*}
Furthermore, from
Item~(2) of Proposition~\ref{prop3}, we have 
$$u_{N_2} := \frac{\alpha_{N_2}}{\alpha_{1+N_2}} = 1 + \mathcal{O}(N_2^{-1}) = 1 + \mathcal{O}(\alpha_{N_2})$$ which leads to~\eqref{bornphi}. 

Considering again a turning point $(r_{N_2}, \varphi_{N_2})$ of a trajectory (with the previous notations), it follows from~\eqref{rntn} and Item~(4) of Proposition~\ref{prop3} that
$$
r_{N_2}^{-\beta} = t_{N_2}^{-\beta} (1 + \mathcal{O}(r_{N_2}^{-1})) = t_{N_2}^{-\beta} (1 + \mathcal{O}(N_2^{-\frac{1}{\beta+1}})),
$$
and from $v_{N_2} \le \frac{\pi}{2} < v_{N_2} + 2 \alpha_{1+N_2}$, combined with Item~(2) of Proposition~\ref{prop3}, that
$$
\sin v_{N_2} = 1 + \mathcal{O}(\alpha_{1+N_2}^2) = 1 + \mathcal{O}(N_2^{-2}).
$$
Combining these estimates with Item~(7) of Proposition~\ref{prop3}, we obtain
$$r_{N_2}^{-\beta} = t_{N_2}^{-\beta} \sin v_{N_2} (1 + \mathcal{O}(N_2^{-\frac{1}{\beta+1}})) = I_{1/\alpha}^{\frac{\beta}{\beta+1}} ({\mathfrak{c}}_i N)^{-\frac{\beta}{\beta+1}} \left(1 + \mathcal{O}(N^{-\frac{1}{\beta+1}})\right).
$$
In particular, $N \ge y$ means that
$$
r_{N_2}^{-\beta} \le \left(\frac{I_{1/\alpha}}{{\mathfrak{c}}_i y}\right)^{\frac{\beta}{\beta+1}} \left(1 + \mathcal{O}\left(y^{-\frac{1}{\beta+1}}\right)\right)\, ,
$$
which is equivalent to~\eqref{Cond2}.
\end{proof}
\begin{proof}[Proof of Proposition~\ref{Mm2}]
Since $\mathcal{F}$ preserves $\mu$ and $\mathcal{F}^{N_2}$ defines a bijection from $\bigcup_{N \ge y} M_N$ onto its image, it follows that
\begin{align}
\mu\left( \bigcup_{N \ge y} M_N^{(i)}\right) \nonumber
&= \mu\left(\bigcup_{N \ge y} \mathcal{F}^{N_2}(M_N^{(i)})\right) \\
&= 2 \int_{r \ge \left(\frac{{\mathfrak{c}}_i y}{I_{1/\alpha}}\right)^{\frac{1}{\beta+1}} (1 + \mathcal{O}(y^{-\frac{1}{\beta+1}}))} \left( \int_{{\mathfrak{c}}_i r^{-\beta-1} (1 + \mathcal{O}(r^{-1}))}^{3 {\mathfrak{c}}_i r^{-\beta-1} (1 + \mathcal{O}(r^{-1}))} \cos \varphi \, d\varphi \right) \, dr \nonumber \\
&= 2 \int_{r \ge \left(\frac{{\mathfrak{c}}_i y}{I_{1/\alpha}}\right)^{\frac{1}{\beta+1}} (1 + \mathcal{O}(y^{-\frac{1}{\beta+1}}))} 2 {\mathfrak{c}}_i r^{-\beta-1} (1 + \mathcal{O}(r^{-1})) \, dr \nonumber \\
&= 4 \mathfrak c_i \beta^{-1} \left(\frac{{\mathfrak{c}}_i y}{I_{1/\alpha}}\right)^{-\frac{\beta}{\beta+1}} (1 + \mathcal{O}(y^{-\frac{1}{\beta+1}})).
\label{muRnB}
\end{align}
This concludes the proof of Proposition~\ref{Mm2}, since 
$\mathfrak c_i^{1-\frac\beta{\beta+1}}=\mathfrak c_i^{\frac 1{\beta+1}}$ and since $\alpha=\frac\beta{\beta+1}$. 
\end{proof}

\subsection{Singular curves of the induced map}

In this section, we investigate the singular curves of the induced system $(F, M)$. 
We will use standard facts from the theory of dispersing billiards \cite{BSC90, BSC91, C99, CM, CM05}. 
For any $x \in \cM$, let $\cK(x)$ represent the curvature of the boundary at the base point of $x$. Consider any $x = (r, \varphi) \in \cM$, with a well-defined image $x_1:=\cF(x) = (r_1, \varphi_1) \in \cM$. Recall that the differential of the billiard map $\mathcal{F}$ at $x$ satisfies:
\begin{align}\label{DTdiff}
D_x \cF = \frac{-1}{\cos \varphi_1} \left(\begin{array}{cc} \tau(x) \cK(x) + \cos \varphi & \tau(x) \\ \tau(x) \cK(x) \cK(x_1) + \cK(x) \cos \varphi_1 + \cK(x_1) \cos \varphi & \tau(x) \cK(x_1) + \cos \varphi_1 \end{array}\right)\, ,
\end{align}
where $\tau(x)$ is the length of the free flight between $x$ and $\mathcal{F}(x)$. 
We define the singular set of $F^{\pm 1}$ as
\begin{equation}\label{:sing}
\mathcal S_{\mp 1} = \partial M\cup
\{x_\asy^{(1)},...,x_\asy^{(\cJ)}\} \cup (M\cap \cF^{\mp 1}(\partial \cM)) \cup \bigcup_{n=1}^{\infty} \partial M_n\, .
\end{equation}
Due to the time-reversibility of the billiard dynamics, $\mathcal S_1$ can be obtained by reflecting $\mathcal S_{-1}$ across the line $\varphi = 0$. Moreover, $F M_n$ is obtained by reflecting $M_n$ across the line $\varphi = 0$. More precisely, a point $(r, \varphi) \in M$ is a singularity point for the map $F$ (i.e., $F(r, \varphi)$ or its differential is not well-defined) if and only if $(r, -\varphi) \in M$ is a singularity point for its inverse $F^{-1}$ (i.e., $F^{-1}(r, -\varphi)$ or its differential is not well-defined).

Consider an infinite cusp with boundary parametrized by $z=\pm f_{\mathfrak c}(t)$ and the configuration $x_\asy = (q_\asy, \vel_\asy) \in \partial Q \times \mathbb{S}^1$ pointing to the common asymptote of 
the cusp.
Note that the point $x_\asy$ belongs to a singular curve of $F$ called $\mathbf{s}_0$, which consists of all points near $x_\asy$ whose forward images by $\mathcal{F}$ are grazing collisions on $\partial Q$. One can check using \eqref{DTdiff} that the slope of the tangent vector at $x = (r, \varphi) \in \mathbf{s}_0$ satisfies
$$
d\varphi/dr = -(\mathcal{K}(x) + \cos \varphi / \tau(x)).
$$
Thus, in the vicinity of $x_\asy$, the curve $\mathbf{s}_0$ can be approximated by a line with slope $-\mathcal{K}(x_\asy)$.
\begin{figure}[h]
\center 
\includegraphics[height=3in]{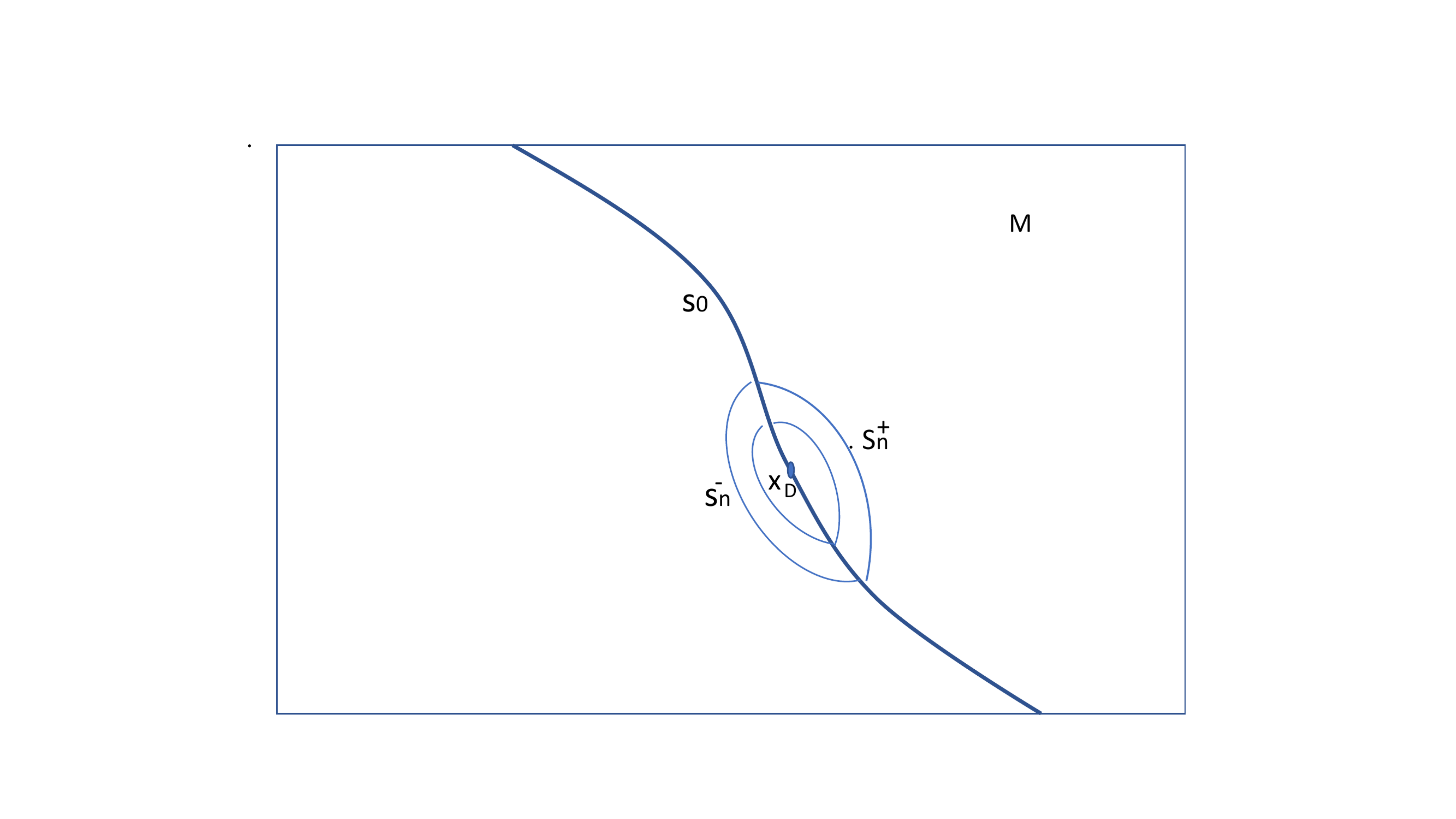}
\renewcommand{\figurename}{Figure.}
\caption{Singularity curves $\bs_0$ and $\bs_n^{\pm}$ of $F$ in the  vicinity of $x_\asy$ in $M$.}\label{FigSing1}
\end{figure}
The singular curves of $F$ near $x_\asy$ consist of two 
sequences of singularity curves $\{\mathbf{s}^+_n\}$ and $\{\mathbf{s}^-_n\}$, with $\mathbf{s}^{\pm}_n \subset \mathcal S_{-1}$, approaching $x_\asy$ from both sides of $\mathbf{s}_0$. More precisely, $\mathbf{s}^{\pm}_n \subset \partial M_n$ consists of points whose trajectories enter the cusp by hitting $\Gamma_{\pm 1}$ first, experiencing $n$ collisions inside the infinite cusp area such that the last collision in the corner series is grazing. Denote by $M^\pm_n$ the strip bounded between $\mathbf{s}^{\pm}_n$ and $\mathbf{s}^{\pm}_{n+1}$. Then $M_n = M^+_n \cup M^-_n$. For symmetry reasons, it is enough to concentrate on $M^+_n$. Figure \ref{FigSing1} shows the structure of the singular curves near $x_\asy$.
\begin{figure}[h]
\center 
\includegraphics[width=6.5in,height=2in]{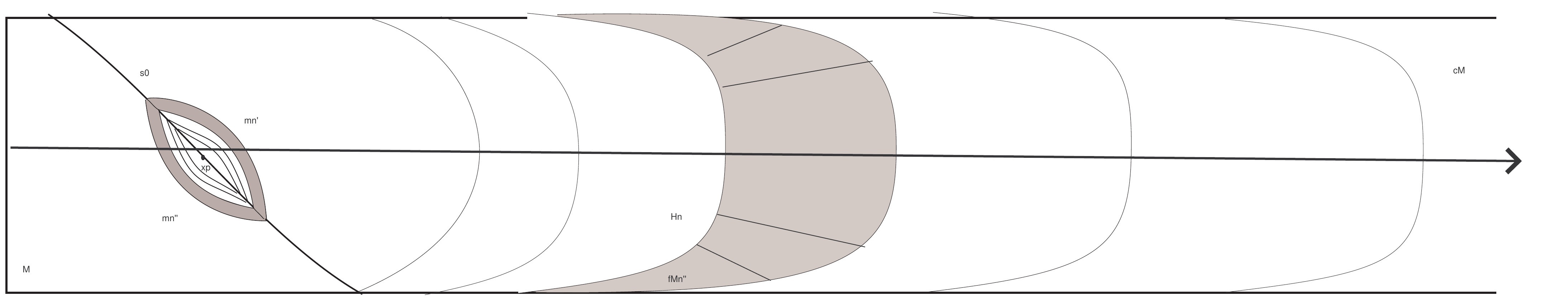}
\renewcommand{\figurename}{Figure.}
\caption{Singularity curves of $F$ in the  vicinity of $x_\asy$; and their forward images under $\cF$ near the cusp in $\cM$ (that are bounded by $H_n$).}\label{FigSing1bis}
\end{figure}
Figure \ref{FigSing1bis} shows the forward iterations of the singular curves near $x_\asy$.
Furthermore, it follows from Proposition \ref{prop3} that, for any $N\geq K_0$, the forward images
$\bigcup_{n=1}^{N-1}\cF^n M_N$
fill entirely the regions squeezed between two curves $H_N$ and $H_{N+1}$ with
equations of the form:
 \beq 
 \sin v= r^{\beta} (I_{1/\alpha}/N)^{\alpha}(1+
\mathcal O(r+ N^{-\frac 1{\beta+1}}))\eeq

\subsection{Properties of stable/unstable manifolds}
\label{pptstable}

Next, we introduce the concept of the wave front. 
Recall that the differential of $\mathcal F$ is given in~\eqref{DTdiff}. 
Due to our assumptions on the table and the definition of $M$, there exist $\tau_{\min}  > 0$ and $0 < \cK_{\min} < \cK_{\max} < \infty$ (depending on our choice of $\widehat t$)  such that:
\beq\label{takck}
\forall x\in M,\quad\tau(x) \ge \tau_{\min}, \quad \cK_{\min} \leq \cK(x) \leq \cK_{\max}\, .
\eeq
Let $V \in \mathcal{T}_x \cM$ be a tangent vector in the tangent space of $\cM$. For $\epsilon > 0$ small, consider an infinitesimal curve $\gamma = \gamma(s) \subset \cM$, where $s \in (-\epsilon, \epsilon)$ is a parameter, such that $\gamma(0) = x$ and $\frac{d}{ds} \gamma(0) = V$. The forward (backward) trajectories of the points $y \in \gamma$,  form a bundle of directed lines in $Q$, called the forward (backward) "wave front." Let $$\mathcal{B} = \mathcal{B}(x),\,\, \mathcal{B}^- = \mathcal{B}^-(x)$$ be the curvature of the orthogonal cross-section of the forward (backward) wave front at the point $x$ with respect to the vector $V$. We have
\begin{equation}\label{cBt}
\mathcal{B}^-(x_1) = \frac{\mathcal{B}(x)}{1 + \tau(x) \mathcal{B}(x)} = \frac{1}{\tau(x) + \frac{1}{\mathcal{B}(x)}}
\,\,\,\text{ and }\,\,\,\,\, \mathcal{B}(x) = \mathcal{B}^-(x) + \frac{2 \mathcal{K}(x)}{\cos \varphi},
\end{equation}
where $x_1 = \mathcal{F} x$.
Denote by $\cV = d\varphi/dr$ the slope of the tangent line of $\gamma$ at $x$. Then $\cV$ satisfies
\beq\label{cVx}
\cV = \cB^{-} \cos \varphi + \cK(r) = \cB \cos \varphi - \cK(r).
\eeq
We use the cone method developed by Wojtkowski \cite{Wo85} to establish hyperbolicity for phase points $x = (r, \varphi) \in M$. In particular, we study stable and unstable wave fronts. The relations in \eqref{cBt} imply that a dispersing wave front remains bounded away from zero.

\begin{defn}\label{defn:1}
The unstable cone $\mathcal{C}_{x}^u$ contains all tangent vectors based at $x$ whose images generate dispersing wave fronts:
$$
\mathcal{C}^u_x = \{(dr, d\varphi) \in T_{x} M \mid \mathcal{K}(r) \leq d\varphi/dr \leq \mathcal{K}(r) + \tau_{\min}^{-1}\}.
$$
Similarly, the stable cones are defined as
$$
\mathcal{C}^s_x = \{(dr, d\varphi) \in {T}_{x} M \mid -\mathcal{K}(r) \geq d\varphi/dr \geq -\mathcal{K}(r) - \tau_{\min}^{-1}\}.
$$
We say that a smooth curve $W \subset M$ is an unstable (resp. stable) curve for $\cF$  if at every point $x \in W$ the tangent line ${T}_x W$ belongs to the unstable (resp. stable) cone $\mathcal{C}^u_x$ (resp. $\mathcal{C}^s_x$).
\end{defn}
We consider a short unstable curve $W \subset M$ with the equation $\varphi = \phi(r)$, such that $FW$ is smooth. Let $x = (r, \varphi) \in W$ and $V = (dr, d\varphi)$ be a tangent vector at $x$ of $W$. Combining with \eqref{cBt} and \eqref{cVx}, we know that for any $x \in M$, the slope $\mathcal V(\mathcal Fx)$ of the tangent vector of $\mathcal{F} W$ at $\mathcal{F} x$ satisfies:
\beq \label{Noooo}
\mathcal{K}(\mathcal{F} x) < \mathcal{V}(\mathcal{F} x) = \mathcal{K}(\mathcal{F} x) + \frac{\cos \varphi(\mathcal{F} x)}{\tau(x) + \frac{\cos \varphi(x)}{\mathcal{K}(x) + \mathcal{V}(x)}} < \mathcal{K}(\mathcal{F} x) + \frac{1}{\tau_{\min}}.
\eeq
This implies that for any $x \in M$, the unstable cone $\mathcal{C}^u_x$ is already strictly invariant under $\mathcal{F}$. Thus, we get $DF \mathcal{C}^u_{x} \subset \mathcal{C}^u_{Fx}$. Similarly, one can check that $DF(\mathcal{C}^s_x) \supset \mathcal{C}^s_{F x}$.

Next, we will show that any unstable vector in $\mathcal{C}^u$ gets uniformly expanded under $DF$. We introduce the $p$-metrics on the tangent space $\mathcal{T} M$. Let $dx = (dr, d\varphi)$, and we define
$$
|dx|_p = \cos \varphi \, |dr|.
$$
Let $\vartheta>\beta+1$. 
 As Sinai did in \cite{Sin70} (the classical homogeneity strips correspond to $\vartheta=2$), we divide $\mathcal{M}$ into horizontal {\bf homogeneity strips} defined by
 $$
 \bH_k := \{(r, \varphi) \in \mathcal{M} \colon \pi/2 - k^{-\vartheta} \leq \varphi \leq \pi/2 - (k+1)^{-\vartheta}\},
 $$
 and
 $$
 \bH_{-k} := \{(r, \varphi) \in \mathcal{M} \colon -\pi/2 + (k+1)^{-\vartheta} \leq \varphi \leq -\pi/2 + k^{-\vartheta}\},
 $$
 for all $k > k_0$, and 
 \begin{equation}
 \label{bbH0}
 \bH_0= \bH_{k_0}:= \{(r, \varphi) \in \mathcal{M} \colon -\pi/2 + (k_0+1)^{-\vartheta} \leq \varphi \leq \pi/2 - (k_0+1)^{-\vartheta}\}.
 \end{equation}
 A technical issue in the study of chaotic billiards is the H\"older continuity of $\varphi\mapsto\log \cos\varphi$ in homogenity strips. This H\"older continuity is ensured in homogeneity strips thanks to the next remark.
 \begin{remark}
 Let $x=(q,\varphi)$, $x'=(q',\varphi')$ be contained in a same homogeneity strip $ \mathbf H_k$, then $\cos\varphi\approx \cos\varphi'\approx k^{-\vartheta}$ and $|\varphi-\varphi'|=\mathcal O\left( k^{-\vartheta-1}\right)$ (uniformly in $k,x,y$ as before), thus $|\varphi-\varphi'|=\mathcal O\left((\cos\varphi)^{\frac{\vartheta+1}{\vartheta}}\right)$ and so
 \begin{equation}\label{controlcos}
 \frac{|\cos \varphi-\cos\varphi'|}{\cos\varphi}
     \le \frac{|\varphi-\varphi'|}{\cos\varphi}=\mathcal O \left(|\varphi-\varphi'|^{1-\frac{\vartheta}{\vartheta+1}}\right)=\mathcal O \left(|\varphi-\varphi'|^{\frac{1}{\vartheta+1}}\right)\, .
 \end{equation}
 For the classical homogeneity strips, $\vartheta=2$ and the exponent is
 ${\frac{1}{\vartheta+1}}=\frac 13$.
\end{remark}
 Here, $k_0 \geq 1$ is a fixed (and usually large) constant.
 For any $N \geq 1$, an unstable curve $W \subset F M_N$ is called a {\bf homogeneous unstable curve} for the induced map $F$ if, for any $k\in\{0,1,N-1,N\}$, the curve
$\mathcal{F}^{-k}W$ is contained in a single homogeneous strip. 
We will state a uniform estimate for the expansion factor $\Lambda_N$ of unstable curves contained in $M_N$, ensuring that its inverse $\Lambda_N^{-1}$ is summable, which will allow us to establish the one-step expansion property. 
To this end, we state the following result, the proof of which is postponed to Appendix~\ref{appendKBbd} and uses the precise study of long corner series done in Proposition~\ref{prop3}.
\begin{proposition}[see Lemma~\ref{propdetails} for details]\label{KBbd}
There exist $N_0\geq 1$ and $c'_0>0$ for all $N \geq N_0$, for all unstable curve $V\subset M_N$ with the properties that: $$\min\cos\varphi(\mathcal FV)\approx\max \cos\varphi(\mathcal FV),\,\,\text{ and }\,\,\,\min\cos\varphi(\mathcal F^{N-1}V)\approx\max \cos\varphi(\mathcal F^{N-1}V)$$ Then the following holds:
\beq\label{P1}
\frac{|FV|_p}{|V|_p}\approx\frac {|FV|}{|V|}
\approx \frac{N^{1+\frac{\beta}{(2\beta+1)(\beta+1)}-\frac{2\beta}{2\beta+1} }}{\cos\varphi(\mathcal FV)\cos\varphi(\mathcal F^{N-1}V)}
\ge c'_0 N^{1+\frac{\beta}{(2\beta+1)(\beta+1)}}\, .
\eeq
Furthermore
\beq\label{P3}
\frac{|\mathcal FV|_p}{|V|_p}\approx N^{\frac\beta{(\beta+1)(2\beta+1)}}
\quad\mbox{and}\quad
|FV|=|\mathcal F^NV|
\approx N^{-\frac{\beta}{2\beta+1}}|\mathcal F^{N-1}V|\, .
\eeq
\end{proposition}
Of course, using time-reversibility, corresponding properties hold true for stable curves under the action of $F^{-1}$.
\begin{lemma}\label{dimMN}
$M_N$ is made of two strips of width ${\mathfrak w}_N=\mathcal O \left(N^{-1 - \frac{\beta}{2\beta+1}}\right)$ (in the unstable direction) and length $|\mathbf{s}_N^\pm|\approx N^{-\frac{\beta^2}{(2\beta+1)(\beta+1)}}$ (in the stable direction). Furthermore
$$\mu(M_N) \approx N^{-\frac{2\beta+1}{\beta+1}}.$$
\end{lemma}
\begin{proof}[Proof of Lemma~\ref{dimMN}]
We use Proposition~\ref{prop3} throughout this proof. 
Let us consider $\mathbf{s}_N^\pm$ the two curves corresponding to $(\partial M_{N}\cap\partial M_{N+1})\setminus \mathbf s_0$. 
Observe that the lengths of these curves
have same order as their $q$-length ${\mathfrak l}_N$, and so has same order as the
$q$-length  of the curve $\mathbf s_0\cap\partial \{R\ge N\}$ of $(q,v)$ such that $\mathcal F(q,\vel)$ is grazing at some position $(t,\pm f_{\mathfrak c}(t))$ with $t\ge t':=c N^{\frac \beta{(2\beta+1)(\beta+1)}}$ (for some $c>0$), which corresponds to a variation of order $N^{-\frac {\beta^2}{(2\beta+1)(\beta+1)}}$
in the $z$-coordinate. Thus
\begin{equation}\label{length}
|\mathbf{s}_N^\pm|\approx {\mathfrak l}_N\approx 
N^{-\frac{\beta^2}{(2\beta+1)(\beta+1)}}\, ,
\end{equation}
This ends the proof of the fact that the length of $M_N$ in the stable direction         
 $|\mathbf{s}_N^\pm|\approx N^{-\frac{\beta^2}{(2\beta+1)(\beta+1)}}$.

Due to the time-reversibility of the billiard dynamics, the singular curves in $F M$ have a similar structure. Furthermore, the short sides of $M_N$ stretch completely under $F M_N$, and are transformed into long sides of $F M_N$. 

Let $W \subset M_N$ be an unstable curve stretches completely in $M_N$ (from $\partial M_{N-1}$ to $\partial M_{N+1}$).
Then  $|FW|\approx {\mathfrak l}_N^{\pm}$. Let us denoted by $c'_1$ and $c'_2$ two positive constants 
such that
\begin{equation}
c'_1{\mathfrak l}_N^{\pm}\le |FW|<c'_2{\mathfrak l}_N^{\pm},  
\end{equation}
for any such $W$. 
It follows from Proposition~\ref{KBbd} that
$$
C N^{1 + \frac{\beta}{(2\beta+1)(\beta+1)}}\le \frac{|FW|}{|W|}=\mathcal O\left(\frac{{\mathfrak l}_N^\pm}{|W|}\right)\, ,
$$
and so that 
\begin{equation}\label{FW/W}
|W|=\mathcal O\left({\mathfrak l}_N^\pm
/N^{1 + \frac{\beta}{(2\beta+1)(\beta+1)}}\right)=\mathcal O\left(N^{-1-\frac\beta{2\beta+1}}\right).
\end{equation}
This provides the announced estimate about 
the width of $M_N$ in the unstable direction. We infer from the above bounds on the width and on the length that
\begin{equation}\label{dommuMN}
\mu(M_N)=\mathcal O\left( N^{-1-\frac\beta{2\beta+1}}N^{-\frac{\beta^2}{(2\beta+1)(\beta+1)}}\right)=\mathcal O\left( N^{-\frac{2\beta+1}{\beta+1}}\right)\, .
\end{equation}
To prove that $\mu(M_N)\approx N^{-\frac{2\beta+1}{\beta+1}}$, we will prove that
$|FW|/|W|\approx \frac{{\mathfrak l}_N^\pm}{|W|}$ for most of the unstable curves $W$
stretching completely $M_N$.
Let $\varepsilon\in(0,1/10)$.

\begin{itemize}
\item Note that  Proposition~\ref{KBbd} ensures also that
\begin{equation}\label{FV/V}
N^{1 + \frac{\beta}{(2\beta+1)(\beta+1)}}\approx \frac{|FV|}{|V|}
\end{equation}
for any unstable curve $V\subset M_N$ such that 
$\mathcal FV\cap B_\varepsilon=\emptyset=\mathcal F^{N-1}V\cap B_\varepsilon$,
where we denote by $B_\varepsilon$ the set $\{\cos\varphi\ge \varepsilon N^{-\frac\beta{2\beta+1}}\}$.
\item Let $W_\varepsilon=W\cap \mathcal F^{-(N-1)}(B_\varepsilon)$.
Observe that $\mathcal F^{N-1}(W_\varepsilon)$ is a curve contained in $B_\varepsilon$. 
But, it follows from Proposition~\ref{prop3},
that the slope of the curve $\mathcal F^{N-1}(W_\varepsilon)$ (seen as a graph $r\mapsto\varphi$) has order $t_{N-1}^{-\beta-2}\approx N^{-\frac{\beta(\beta+2)}{(2\beta+1)(\beta+1)}}$. Gathering these two facts, we obtain that
\begin{equation}\nonumber
\frac{|\mathcal F^{N-1}W_\varepsilon|}{\varepsilon}\approx N^{\frac{\beta(\beta+2)}{(2\beta+1)(\beta+1)}}N^{-\frac\beta{2\beta+1}}=N^{\frac\beta{(\beta+1)(2\beta+1)}}
\end{equation}
and it then follows from the last estimate of Proposition~\ref{KBbd} that
\begin{equation}\nonumber
\frac{|FW_\varepsilon|}{\varepsilon}\approx N^{-\frac \beta{2\beta+1}}
\frac{|\mathcal F^{N-1}W_\varepsilon|}{\varepsilon}\approx 
N^{-\frac{\beta^2}{(\beta+1)(2\beta+1)}}\approx {\mathfrak l}_N^\pm.
\end{equation}
This ensures the existence of $\varepsilon_0$ such that if $\varepsilon<\varepsilon_0$, then 
\begin{equation}\nonumber
|F(W_\varepsilon)|_q< \frac{c'_1}{10}{\mathfrak l}_N^\pm,
\end{equation}
for any $W$ streching completely $M_N$ (from $\partial M_{N-1}$ to $\partial M_{N+1}$),
where we write $|\cdot|_q$ for the $q$-length.
\item Assume here only that $\mathcal FW$ intersects the set $ B_\varepsilon$. Then there exists a stable curve $\tilde W_\varepsilon$ contained in $M_N$, joining $W$ and $\mathbf s_0\cap\partial M_N$ such that $\mathcal F\tilde W_\varepsilon$ is contained in $B_\varepsilon$.
By time reversibility, it follows from the previous point that if $\varepsilon<\varepsilon_0$, then
\begin{equation}\nonumber
|\tilde W_\varepsilon|_q< \frac{c'_1}{10}{\mathfrak l}_N^\pm.
\end{equation}
\end{itemize}
 Assume now that $\varepsilon<\varepsilon_0$ and that $d_q(W,\mathbf s_0\cap\partial M_N)> \frac {c'_1}{10}{\mathfrak l}_N^\pm$ (writing $d_q$ for the distance in terms of the coordinate $q$). Then the curve $V_\varepsilon=W\setminus W_\varepsilon$ satisfies~\eqref{FV/V} and we infer that
\begin{equation}\nonumber
|W|>|V_\varepsilon|\approx N^{-1 - \frac{\beta}{(2\beta+1)(\beta+1)}}|FV_\varepsilon|\ge 
N^{-1 - \frac{\beta}{(2\beta+1)(\beta+1)}}\frac{9}{10}{\mathfrak l}_N^\pm\approx 
N^{-1-\frac\beta{2\beta+1}}.
\end{equation}
This combined with~\eqref{FW/W} ensures that these curves $W$ satisfy $|W|\approx N^{-1-\frac\beta{2\beta+1}}$.\\
To conclude, we observe that, for any $N$ large enough, we can consider a family of such curves $W$ (for example made of parallel pre-collisional vectors, i.e. such that
$d\varphi/dr=\max K(r)$) covering the set of points $(q,v)\in M_N$ with position
at distance at least $\frac{8}{10}c'_1{\mathfrak l}_N^\pm$ to $\mathbf s_0\cap\partial M_N$.
Hence 
\begin{equation}\nonumber
N^{-1-\frac \beta{2\beta+1}}{\mathfrak l}^\pm_N=\mathcal O(\mu(M_N)).
\end{equation}
This combined with~\eqref{dommuMN} ensures that $\mu(M_N)\approx N^{-\frac{2\beta+1}{\beta+1}}$.
\end{proof}
The following technical fact will be useful in the proof of the one-step expansion.
 \begin{lemma}\label{rk0}
There exists an integer $K'$ such that, if $W$ is an unstable curve  crossing completely $M_{N},...,M_{K'N }$, then $W$ has length of order $N^{-\frac{\beta}{2\beta+1}}$.

 \end{lemma}
\begin{proof}

It follows from~Proposition~\ref{prop3} that there exist constants $c''_1,c''_2$ such that
$$x\in M_N\Rightarrow t_1/ N^{\frac \beta{(2\beta+1)(\beta+1)}}\in[c''_1;c''_2].$$
This implies that for $x,y\in M_N$, we have
$$t_1(x)-t_1(y)
\leq (c''_2 -c''_1)N^{\frac \beta{(2\beta+1)(\beta+1)}}$$
Let $K'$ be an integer such that $c''_2<c''_1 K' $. 
If $W$ crosses totally $M_{N},...,M_{ K'N }$, 
 then $W$ crosses the region $\mathcal F^{-1}(D_N)$ 
 with $D_N$ to be the set of points $x=(q,\vel)\in\mathcal F(M)$ with $q$ corresponding to $c''_2 N^{\frac \beta{(2\beta+1)(\beta+1)}}<t<c''_1 (K'N)^{\frac \beta{(2\beta+1)(\beta+1)}}$. 
 But, recall that the slope of $\mathcal F(W)$ has then order $ t^{-\beta-2}\approx N^{-\frac {\beta(\beta+2)}{(2\beta+1)(\beta+1)}}$ (due to~\eqref{Noooo} since $\kappa\circ\mathcal F\approx t^{-\beta-2}$,
$\cos(\varphi\circ\mathcal F)=\mathcal O( t^{-\beta-1})$ and $\tau\approx t$). Thus the variation of $\varphi$ along $\mathcal F(W)$ has order $ t^{-\beta-2}t\approx t^{-\beta-1}\approx N^{-\frac\beta{2\beta+1}}$.
We conclude that there exists a curve $W'$ contained in $W$ such that $|\mathcal FW'|\approx|\mathcal FW|$      and such that $|\cos(\varphi(\mathcal FW'))|\approx  t^{-\beta-1}$. Thus $|\cF W'|_p\approx t^{-\beta-1} t=t^{-\beta}$ and, using the first estimate of~\eqref{P3}, it follows that
 $$ |W|\ge|W'|\approx |W'|_p\approx N^{-\frac\beta{(\beta+1)(2\beta+1)}} |\mathcal F W'|_p\approx N^{-\frac\beta{(\beta+1)(2\beta+1)}} t^{-\beta} \approx N^{-\frac\beta{2\beta+1}}. $$
But it follows from Lemma~\ref{dimMN} that $|W|=\mathcal O(N.N^{-1-\frac{\beta}{2\beta+1}})$.
So we conclude that $|W|\approx N^{-\frac\beta{2\beta+1}}.$
\end{proof}

\subsection{Decay of correlations for the induced system}\label{sec:decaycorr}

Decorrelation estimates will follow from the next crucial estimate,
via the notion of standard pair and standard family.
 \begin{lemma}\label{onestep}(\textbf{One-step expansion estimate})  Let $W$ be a short
 unstable curve in $M$ for the induced map $F$, and $\{W_i\}$ be the collection of
components in $W$ such that $\{FW_i\}$ are
homogeneous unstable curves.  Then
\begin{equation}\label{step1}
\liminf_{\delta_0 \to 0} \sup_{W \colon |W| < \delta_0} \sum_{i \geq 1} \frac{1}{\Lambda_i} < 1,
\end{equation}
where the supremum is taken over unstable curves $W \subset M$
such that $|W|<\delta_0$ and $\Lambda_i := \frac{|FW_i|}{|W_i|}$ denotes the local expansion factor of the connected component $W_i$ under the map $F$.
\end{lemma}
\begin{proof}
 Let $W \subset M$ be an unstable curve such that $|W|<\delta_0$. 
 Define $$W_{N,m,k} = W \cap M_N\cap \mathcal{F}^{-(N-1)} \mathbf{H}_k \cap \mathcal{F}^{-1} \mathbf{H}_m.$$
 Let us write $B_N$ for the set of couples $(m,k)$ such that
 $W_{N,m,k} $ is not empty. 
 It follows from Item~(3) of Proposition~\ref{prop3} that
 elements of $B_N$ should satisfy $k\ge m_N$ and $m\ge m_N$, 
 with $m_N\approx N^{\frac \beta {\vartheta(2\beta+1)}}$ (recall that $\vartheta$ has been chosen before~\eqref{bbH0}). 
 Furthermore, as in~\cite{CM05}, we observe that one of the two following monotonicity properties holds: 
 \begin{itemize}
  \item either: for all $(m,k)\in B_N$ such that $(m+1,k')$ is in $B_N$ for some $k'$, then $k'\ge k$;
 \item or: for all $(m,k)\in B_N$ such that $(m+1,k')$ is in $B_N$ for some $k'$, then $k'\le k$.
 \end{itemize}
 Thus, for any $n\in\mathbb Z$, there
 exists at most a couple $(k,m)\in B_N$ such that
 $k+m=n$ (or, resp. $k-m=n$) and this couple satisfies 
 $$n\ge 2m_N\quad\mbox{and}\quad  (k/m_N)^\vartheta(m/m_N)^\vartheta\ge  (\max(k,m)/m_N)^\vartheta\ge (|n|/(2m_N))^\vartheta\, ,$$
 or resp.
 $$
 (k/m_N)^\vartheta(m/m_N)^\vartheta\ge  (\max(k,m)/m_N)^\vartheta 
 \ge ((m_N+|n|)/m_N)^\vartheta.$$
 Furthermore Proposition \ref{KBbd} implies that the expansion factor satisfies:
 \begin{align*}
 \Lambda_{N,m,k}^{-1} := \frac{|W_{N,m,k}|}{|FW_{N,m,k}|} &\approx N^{-1 - \frac{\beta}{(2\beta+1)(\beta+1)}} (k/m_N)^{-\vartheta}(m/m_N)^{-\vartheta}\, .
 \end{align*}
 Thus, for all integers $n_1\le n_2$, 
 \begin{align*}
 \sum_{N = n_1}^{n_2} \sum_{(m,k) :W_{N,m,k}\ne\emptyset} 
 &\frac{1}{\Lambda_{N,m,k}}=\mathcal O\left(\sum_{N= n_1}^{n_2}\sum_{n'\ge m_N}N^{-1 - \frac{\beta}{(2\beta+1)(\beta+1)}} (n'/m_N)^{-\vartheta}\right)\\
 &=\mathcal O\left(  \sum_{N= n_1}^{n_2}N^{-1 - \frac{\beta}{(2\beta+1)(\beta+1)}}  m_N\right) =\mathcal O\left(\sum_{N= n_1}^{n_2}N^{-1-\frac{\beta}{2\beta+1}\left(\frac 1{\beta+1}-\frac 1\vartheta\right)}\right)\\
 &\le c''n_1^{-\frac{\beta}{2\beta+1}\left(\frac 1{\beta+1}-\frac 1\vartheta\right)}\, .
 \end{align*}
Assume there exists $n_1$ such that $W$ crosses $M_N$ for all $N \in\{ n_1,...,K'n_1\}$ (with $K'$ of Lemma~\ref{rk0}) and that $W$ does not cross fully $M_{n_1-1}$. Then
it follows from Lemma~\ref{rk0}
that 
$$\delta_0>|W|\ge c'  n_1^{-\frac\beta{2\beta+1}}\, ,$$
for some $c'>0$.
Combining the two above estimates, in that case, we obtain that
$$
 \sum_{N = n_1}^{n_2} \sum_{(m,k) :W_{N,m,k}\ne\emptyset} 
 \Lambda_{n,m,k}^{-1}
\le c''n_1^{-\frac{\beta}{2\beta+1}\left(\frac 1{\beta+1}-\frac 1\vartheta \right)}\le c''(\delta_0/c')^{\frac 1{\beta+1}-\frac 1\vartheta}
$$
which is strictly smaller than 1 provided $\delta_0$
is small enough  (since $\vartheta>\beta+1$).\\
Assume now that $W$ intersects $M_N$ for all $N \in\{ n_1,...,n_2\}$ and no other, with $n_2<K'n_1$.
Then
$$\sum_{N = n_1}^{n_2} \sum_{(m,k) :W_{N,m,k}\ne\emptyset} 
\frac{1}{\Lambda_{N,m,k}}\ll n_1^{-\frac{\beta}{2\beta+1}\left(\frac 1{\beta+1}-\frac 1\vartheta\right)}
$$
which is $<1$ as soon as $n_1$ is larger than some $N_0$.\\
If $n_1\le N_0$ and still $n_2<K'n_1$, then $n_2\le K' N_0$, this case can be treated as stable curves in $M_1$.
\end{proof}
The notion of a {\bf standard pair} was studied  by Chernov and Dolgopyat in \cite{CD}. A standard pair is a couple $(W,\nu_{\mathcal W})$ where $W$ is an homogeneous unstable curve and where $\nu_{\mathcal W}$ is a probability measure on $W$ with smooth density with respect to the length $m_W$. A {\bf standard family} $(\mathcal G,\nu_{\mathcal G})$ is then a family of standard pairs
 $\{(W_\aaa, \nu_\aaa)\,:\, \aaa\in \cA\}$ where $\cA\subset [0,1]$ and where $\cW:=\{W_\aaa\}_{\aaa\in\cA}$ is a measurable family of
unstable curves
endowed with a probability Borel measure $\lambda_{\cG}$ (factor measure) on $\cA$, which defines a measure $\nu_{\cG}$
 on $M$  by
\beq\label{cGnm}
   \nu_{\cG}(B):=\int_{\cA} \nu_{\aaa}(B\cap
   W_{\aaa})\,
   d\lambda_{\cG}(\aaa)\hspace{1cm}
   \eeq
for all  measurable sets $B\subset M$. The probability measures
$\nu_\aaa$ appear as the conditional probability of $\nu_{\mathcal G}$ conditioned to be in the unstable curve $W_\aaa$. 
We then consider the function $\cZ
$ defined on the set of standard families by the following formula
\beq\label{cZ}
\cZ
(\cG):=
\int_{\cA}|W_{\aaa}|^{-1
}\,d\lambda_{\cG}(\aaa). \eeq
For any unstable curve $W\in \cW$, any $x\in W$,  and any $k\geq 1$,
let $W^{k}(x)$ be the largest smooth homogeneous unstable curve contained in $F^k \cW$ that contains $F^k x$ and such that $F^{-i}(W^{k}(x))$ is homogeneous for all $i=1,...,k$. 
We define $r_k(x)$ as the minimal distance between $F^k x$ and the closest end points of $W^k(x)$, measured along $W^k(x)$.
We next state the Growth Lemma, proved in \cite{CZ09}, which is key in the analysis of hyperbolic systems with singularities. It expresses the fact that the expansion of unstable curves dominates the cutting by singular curves, in a uniform fashion for all sequences. The reason behind this is that unstable curves expand at a uniform exponential rate, whereas the cuts accumulate at only a finite number of singular points. The following Growth Lemma can be derived directly from Lemma \ref{onestep} — the One-step Expansion Estimate, see \cite{CM}, \cite{CZ09} for details (using also the distortion bound established in Appendix~\ref{distortion}). 
\begin{lemma}[Growth Lemmas]\label{growthlemma}
Let $\cG=(\cW, \nu)$ be a standard family. Then, there exists $C_0>0$, $C_1>0$ and $\theta_0\in(0,1)$ such that, for any $\eps>0$ and any $k\ge 0$,
$$\nu_{\mathcal G}(r_k<\eps)\le C_0
\eps \cZ(F^k\cG)\leq C_1(\theta_0^{k}\cZ
(\cG)+C_2)\eps\, .
$$
\end{lemma}
 The induced map
behaves as
a dispersing billiard system, which is uniformly hyperbolic and also enjoys exponential mixing speed on H\"older observables. Indeed the induced system satisfies all conditions in
 \cite{CZ09}.

For the reader's convenience we record explicitly how the structural conditions of the Chernov--Zhang framework---i.e.\ the assumptions (H1)--(H5) in the formulation of \cite{C99} (see also \cite[Ch.~4 and 7]{CM}), together with the one-step expansion and growth hypotheses of \cite{CZ09}---are verified for the induced map $(M,F,\mu_M)$:
\begin{itemize}
\item[\textbf{(H1)}] \emph{Hyperbolicity and invariant cones.} The stable/unstable cone fields $\mathcal C^{s}_x,\mathcal C^{u}_x$ are defined in Definition~\ref{defn:1}; their strict invariance under $\cF$ is \eqref{Noooo}, while the \emph{uniform} expansion of $F$ (the expansion factor on $M_N$ is $\gtrsim N^{1+\frac{\beta}{(2\beta+1)(\beta+1)}}$ for $N$ large, and $\ge\Lambda_{\min}>1$ uniformly) is Proposition~\ref{KBbd}. Thus $F$ is uniformly hyperbolic even though the underlying $\cF$ is not.
\item[\textbf{(H2)}] \emph{Singularities.} The singular sets $\mathcal S_{\pm1}$ are countable unions of compact $C^{2}$ curves, see \eqref{:sing} and the sequences $\{\mathbf s^{\pm}_n\}$ accumulating only at $x_\asy$; the alignment (monotonicity) of these curves used to bound complexity is the property exploited in the proof of Lemma~\ref{onestep}, and the dimensions of the cells $M_N$ are quantified in Lemma~\ref{dimMN}.
\item[\textbf{(H3)}] \emph{Regularity of unstable curves.} Unstable curves lie in the bounded cones of Definition~\ref{defn:1} (bounded slope) and have uniformly bounded curvature, by the bounds \eqref{takck} on $\tau,\cK$ and the wavefront relations \eqref{cBt}; the homogeneity strips $\bH_k$ (with exponent $\vartheta>\beta+1$) control the H\"older oscillation of $\log\cos\varphi$, see \eqref{controlcos}.
\item[\textbf{(H4)}] \emph{Distortion bounds.} The Jacobian of $F$ along homogeneous unstable curves is dynamically $\tfrac1{\vartheta+1}$-H\"older, Proposition~\ref{distortionbound}.
\item[\textbf{(H5)}] \emph{Absolute continuity.} The stable holonomy of $F$ is absolutely continuous with Jacobian bounded and dynamically H\"older, exactly as for any dispersing billiard map, cf.~\cite[Ch.~5]{CM}.
\item[\textbf{(OS)}] \emph{One-step expansion.} The summability $\liminf_{\delta_0\to0}\sup_{|W|<\delta_0}\sum_i\Lambda_i^{-1}<1$ is Lemma~\ref{onestep}; it is precisely here that the cusp order $\beta$ enters, through the expansion rate $\Lambda_N\approx N^{1+\frac{\beta}{(2\beta+1)(\beta+1)}}$ on $M_N$ and the requirement $\vartheta>\beta+1$ on the homogeneity exponent.
\item[\textbf{(GL)}] \emph{Growth lemma and equidistribution.} The growth lemma for standard families is Lemma~\ref{growthlemma}, and the equidistribution of proper families is Lemma~\ref{equidistribution} (\cite[Thm~2]{CZ09}); both belong to the standard-pair framework of \cite{CZ09}.
\end{itemize}
Conditions (H1)--(H5), (OS) and (GL) together yield, via \cite[Thm~3]{CZ09}, the exponential decay of correlations \eqref{correT} for the induced system.

For $\gamma\in (0,1)$, let $\cH_{\gamma}^\pm$ be the set of bounded functions $f$ such that the following quantity is finite
$$[f]_{\gamma}^\pm:= \sup_{ (x, y)\in\mathcal E^\pm}\frac{|f(x)-f(y)|}{\gamma^{s_{\pm}(x,y)}}<\infty,
$$
where $\mathcal E^+$ (resp. $\mathcal E^-$) is the set of couples $(x,y)$ lying on a same unstable (resp. stable) manifold and where $s_{+}(x,y)$ (resp. $s_{-}(x,y)$) is the forward (resp. backward) separation time, that is the infimum of integers $k\ge -1$ such that
$F^{k+1}(x)$ and $F^{k+1}(y)$ (resp. $F^{-k-1}(x)$ and $F^{-k-1}(y)$) do not belong to the same connected component of $M\setminus \mathcal S_{- 1}$ (resp. of $M\setminus \mathcal S_{1}$).
For every $f\in\mathcal H_\gamma^\pm$, we set
\beq \label{defCgamma}
\|f\|_{\gamma}^\pm:=\|f\|_{\infty}+[f]^\pm_{\gamma}.
\eeq

For a fixed large constant $C_\text{prop}>0$, any standard family $\cG$ with $\cZ(\cG)<C_\text{prop}$ will be called a {\it{proper}} family. The following was proved in \cite[Theorem 2]{CZ09}. 

 \begin{lemma}[Equidistribution]\label{equidistribution}
If $\cG=(\cW,\nu)$ is a proper family, then for any $g\in \cH_{\gamma}^+$ with $\gamma\in (0,1)$ and $k\ge 0$,
\beq\label{properg}
\left|\int_Mg\circ F^k\, d\nu-\int_Mg\, d\mu_M\right|\leq C \|g\|_{\gamma}^+\vartheta^k.
\eeq
\end{lemma}
For future reference, we record here a bound from \cite[Theorem 3]{CZ09}. For $f\in \cH_{\gamma}^+, g\in \cH_{\gamma}^-$, and any integer $k$, the correlations of $f$ and $g\circ   F^k$ satisfy:
\beq\label{correT}\left|\text{Cov}_{\mu_M}(f,g\circ   F^k)\right|=\left|\int_Mf\cdot g\circ   F^k\, d\mu_M-\int_Mf\, d\mu_M\int_Mg\, d\mu_M\right|\leq C \|f\|_{\gamma}^+ \|g\|_{\gamma}^-\vartheta^k\eeq
where $C>0$ and $\vartheta\in (0,1)$ are constants. 
Observe that \eqref{correT} applies to observables of the form
$f=f_0\cdot(f_1\circ F^{-1})\cdots (f_k\circ F^{-m})$ and $g=g_0\cdot (g_1\circ F)\cdots (g_k\circ F^\ell)$, with the $f_j$'s belonging to $\mathcal H_\gamma^+$
and the $g_j$'s belonging to $\mathcal H_\gamma^-$.

\subsection{Proof of Theorem~\ref{LLNfails}}\label{sec44}
In this Subsection, we provide the proof of Theorem~\ref{LLNfails}. 
Equation~\eqref{correT} ensures in particular the mixing and so the ergodicity of $(M,F,\mu_M)$.
This combined with Remark~\ref{rktower}
implies the recurrent ergodicity of $(\cM,\cF,\mu)$.
The Hopf ergodic theorem then ensures that, for any integrable functions $f$ and $g$ such that $\int_{\mathcal M}g\, d\mu\ne 0$, the following convergence holds true almost everywhere: $\lim_{n\rightarrow +\infty}\frac{S_nf}{S_ng}=\frac{\int_{\mathcal M}f\, d\mu}{\int_{\mathcal M}g\, d\mu}$. 
But this holds also true for $g=\mathbf 1$. Using the fact that $|\mu(f)|<\|f\|_1$, and $\mu(g)=\infty$, we obtain $$\lim_{n\to\infty} \frac{S_nf}{S_ng}=0$$ 
Since $S_ng=n$, it follows that $$\lim_{n\to\infty} \frac{S_nf}{n}=0$$as we claimed. This finishes the proof.
This ends the proof of Theorem~\ref{LLNfails}.
\subsection{Correlations for the Return time function}

In this Subsection, we prove an important property for the return time function. 
\begin{lemma}\label{deco1}
There exists $\theta_0>0$ such that 
$$\sup_{k\ge 1}
\mu_M(R>n^{\frac 1\alpha},R\circ F^k>n^{\frac 1\alpha})=\mathcal O\left(  n^{-1-\theta_0}\right)\, ,$$ 
as $n\rightarrow +\infty$.
\end{lemma}

\begin{proof}
We first define $$A_{n,j}=\{x\in M\,|\, R\circ F^j(x)> n^{1/\alpha}\}=F^{-j}(A_{n,0})\, .$$
For $k\ge 1$, we write
$$
\mu_M(A_{n,0}\cap A_{n,k})=\mu_M(A_{n,-1}\cap A_{n,k-1}).$$
Furthermore note that
$$A_{n,-1}=\cup_{m=\lfloor n^{1/\alpha}\rfloor+1}^{\ff}F(M_m),$$ and $F(M_m)$ is a union of two strips, of width $\le cm^{-1-\frac\beta{2\beta+1}}=\mathcal O\left(n^{-\frac{(\beta+1)(3\beta+1)}{\beta(2\beta+1)}}\right)$ (in the stable direction), made of a union of unstable curves $\{W_{\bar\alpha}\}$, such that 
$$|W_{\bar\alpha}|\approx |{\mathbf s}_m^\pm|\approx m^{-\frac{\beta^2}{(2\beta+1)(\beta+1)}}=\mathcal O\left( n^{-\frac {\beta}{2\beta+1}}\right)\, .$$ 
Recall that $\mu_M(A_{n,0})\approx n^{-1}$ (due to Proposition~\ref{Mm2}). 
We consider the associated standard family $\cG_a=(\cW_a^{(n)},\nu_a^{(n)})$ with $\nu_a^{(n)}:=\mu|_{A_{n,-1}}/\mu(A_{n,0})$; with factor measure $\lambda_{\cG_a^{(n)}}\approx n \mu(M_m)\mbox{Leb}
$ and with conditional measure $\nu_W\approx \mbox{length}/|{\mathbf s}_m^\pm|$ if $W\in \mathcal W_a^{(n)}$ is contained in $F(M_m)$. 
Recall that $\mu(M_m)=\mathcal O\left(m^{-\frac{2\beta+1}{\beta+1}}\right)$ (see Lemma~\ref{dimMN}). 
Thus
\begin{align*}
\cZ
(\cG_a^{(n)})&
=\int_{\cA}|W_{\bar\alpha}|^{-1}
\, d\lambda_{\cG_a^{(n)}}(\bar\alpha)\\
&
=\mathcal O\left( n 
\sum_{m= n^{1/\alpha}}^{\ff}m^{\frac{\beta^2}{(2\beta+1)(\beta+1)}} m^{-\frac{2\beta+1}{\beta+1}}
\right)
=\mathcal O\left( n \, n^{-\frac {\beta+1}{2\beta+1}}
\right)\, .
\end{align*}
Now observe that $A_{n,k-1}\subset \{r_{k-1}<\eps'_n\}$, with
$\eps'_n$ the maximal unstable diameter of $M_m$ for $m\ge n^{1/\alpha}$. Thus, due to Lemma~\ref{dimMN}, $\eps'_n=\mathcal O\left(n^{-\frac{(\beta+1)(3\beta+1)}{\beta(2\beta+1)}}\right)$.
It thus follows from Lemma~\ref{growthlemma} that
\begin{align*}
\mu_M\left(A_{n,0}\cap A_{n,k}\right)&=\mu_M\left(A_{n,-1}\cap A_{n,k-1}\right)\\
&=\mu_M\left(A_{n,-1}\cap \{ r_{k-1}<\eps'_n\}\right)\\
&=\mu_M(A_{n,0})\nu_a^{(n)}(r_{k-1}<\eps'_n)\\
&=\mathcal O\left( n^{-1}(\mathcal Z(\mathcal G_a^{(n)})+1)\eps'_n \right)\\
&=\mathcal O\left( n^{-\frac{\beta+1}{2\beta+1}} \, n^{-\frac{(\beta+1)(3\beta+1)}{\beta(2\beta+1)}}\right)=\mathcal O\left(n^{-\frac{(\beta+1)(4\beta+1)}{\beta(2\beta+1)}}\right)\, .\end{align*}
\end{proof}

\part{Proofs of limit theorems} \label{part3}

In this part, we present a general strategy to establish limit theorems and apply it to prove Theorems~\ref{thm:LLN},~\ref{thm:TCL},~\ref{thm:TCLflow} and~\ref{thm:TCLflow2}.
This strategy is based on induction on a finite measure set $M$.

In Section~\ref{secF}, we prove that the CSP and the FCLT for a map preserving an infinite measure follow from functional limit theorems for the probability preserving induced dynamical system. More precisely, we will see that the CSP and the FCLT in infinite measure follow from a (simple or joint) limit theorem for the induced system.

In Section~\ref{proofformap}, we prove, in a general context, that the above mentioned (simple or joint)
FCLT for the induced system appears as a consequence of a (simple or joint) limit theorem involving the point process of return times to the finite measure set $M$. We prove Theorems~\ref{thm:LLN},~\ref{thm:TCL} as an application of our general results.

In Section~\ref{proofTCLflow}, we establish a general limit theorems for flows as a consequence of the previous results. Again, we apply our general results to prove Theorems~\ref{thm:TCLflow} and~\ref{thm:TCLflow2}.

\section{General results in infinite measure based on induction and point processes}\label{secF}
Throughout this section, we consider the following general context based on induction, coupling and point processes.
In subsections~\ref{secF1} and~\ref{secF2}, we will
establish a result of Convergence to a Stable Process (CSP) and a Functional Central Limit Theorem (FCLT) in infinite measure via limit theorems established for the induced map on a set of finite measure.
Limit theorems in infinite measure using induction have been proved by Thomine in~\cite{Thomine1,Thomine2} (see also~\cite{PT1,PT2}). We present here a general approach to establish {\bf functional limit theorems}, via induction and point processes. 
Whereas a convergence in distribution of a single point process will be enough for the Convergence to a Stable Distribution (CSP), the convergence of double processes will be needed for the Central Limit Theorem. This convergence of double processes is a weaker condition than those appearing in~\cite{Thomine1,Thomine2} (see e.g.~\cite[Theorem 1.7]{Thomine1}).
In subsection~\ref{sec:F3}, we will give a criteria, based on the convergence of a family of point processes, to ensure the convergence of the double processes.  In this section, we will
make the following assumption.
\begin{hypothesis}\label{hypinduction}
Assume $(\cM,\cF,\mu)$ is a conservative infinite measure preserving dynamical system. Let $M\subset\cM$ be such that $\mu(M)\in (0,+\infty)$. Consider the induced system
$(M,F,\mu_M)$ corresponding to the induction on $M$, that is given by the first return map $F$ to $M$: $$F(x)=\cF^{R(x)}(x),\quad\mbox{with }R(x):=\inf\{n\ge1\, :\, F^n(x)\in M\}\, ,$$
and where $\mu_M$ is the restriction of $\mu/\mu(M)$ to the set of measurable sets contained in $M$.
\end{hypothesis}
For any $f:\cM\rightarrow\mathbb R$, we write $(S_nf)_n$ for the sequence of Birkhoff sums with respect to the initial system, i.e.
$$
S_nf:=\sum_{k=0}^{n-1}f\circ \cF^k,\quad S_0f=0\, .
$$
We set
\begin{equation}\label{defR}\mathcal R(u):=
\sum_{k=0}^{\lfloor u
\rfloor-1} R\circ F^k\, .
\end{equation}
\begin{hypothesis}\label{hypDelta}
Let $\alpha\in(0,1)$, $c_0>0$ and let $\Delta$ be an $\alpha$-stable
L\'evy process (with independent stationary increments) such that $\Delta(1)$ has characteristic function $$\Psi_{\alpha,c_0}:t\mapsto\exp(-|c_0t|^\alpha\Gamma(1-\alpha)e^{-i\frac{\pi\alpha}2\sign(t)}).$$
Furthermore, we define $\Delta^{-1}$ as the generalized inverse of $\Delta$ defined as follows
$$
\Delta^{-1}(t)=\inf\{s>0\, :\, \Delta(s)\ge t\}\, .
$$
\end{hypothesis}
Recall that $\Delta(s)$ has the same distribution as $s^{\frac 1\alpha}\Delta(1)$.
The proofs of Sections~\ref{secF1} and~\ref{secF2} use the
fact that $\Delta$ is c\`adl\`ag and strictly increasing, and
so that $\Delta^{-1}$ is continuous.
Below we will write $\bI_A$ for the indicator function of a set $A$.

\subsection{Convergence to a Stable Process}\label{secF1}
\begin{proposition}\label{corocvgceloi}
Assume Hypotheses~\ref{hypinduction} and~\ref{hypDelta} and that
\begin{equation}\label{cvgcedistrib}
\left(\varepsilon^{\frac 1\alpha}\mathcal R( s/\varepsilon)\right)_{s\geq 0}\Longrightarrow\left(\Delta(s)\right)_{s\geq 0}\, ,\quad\mbox{as }\varepsilon\rightarrow 0
\end{equation}
in the sense of finite dimensional distributions (with respect to $\mu_M$). Then, 
for all $T>0$,  
\beq \label{cvgeSn}
\left(S_{\lfloor nt\rfloor}\bI_{M}/n^\alpha\right)_{t\in[0,T]}\Longrightarrow
(\Delta^{-1}(t))_{t\in[0;T]}\quad \mbox{as }n\rightarrow +\infty 
\eeq
with respect to $\mu_M$, for the uniform metric.
\end{proposition}
\begin{proof}[Proof of Proposition~\ref{corocvgceloi}]
We start by making some comments on $\Delta$ and $\Delta^{-1}$. 
Recall that the process $s\mapsto \Delta(s)$ is almost surely strictly increasing, c\`adl\`ag (with discontinuities) and unbounded. Furthermore, by definition of the generalized inverse, for any $t,s>0$, 
$$
\Delta^{-1}(t)\le s\quad\Leftrightarrow\quad t\le \Delta(s)\, .
$$
Moreover $\Delta^{-1}$ is continuous, increasing (but not strictly increasing).

Observe that $n\mapsto S_n\bI_M$ corresponds to the inverse of $\mathcal R$ since, on $M$,
$$
S_n\bI_M\le m \quad\Leftrightarrow\quad \mathcal R(m)\ge n\, .
$$
Due to \cite[Theorem 1]{Zweimuller}, it is enough to prove the convergence in distribution of $((S_{\lfloor nt\rfloor}\bI _M/n^\alpha)_{t\in[0;T]})_n$ to 
    $(\Delta^{-1}(t))_{t\in[0;T]}$ with respect to the probability measure $\mu_M$. 
Let $J\in\mathbb N^*$, $0<s_1<...<s_J$ and $0<t_1<...<t_J$. 
We observe that
\begin{align*}
\mu_M\left(\forall j=1,...,J,\ S_{\lfloor nt_j\rfloor}\bI_{M}\le n^\alpha s_j\right)&=\mu_M\left(\forall j=1,...,J,\ \mathcal  R(\lfloor n^\alpha s_j\rfloor)\ge \lfloor nt_j\rfloor\right)\\
&=\mu_M\left(\forall j=1,...,J,\ t_j\lfloor nt_j\rfloor^{-1} \, \mathcal R( n^\alpha s_j)\ge t_j\right)
\end{align*}
which converges to $\mathbb P(\forall j=1,...,J,\ \Delta(s_j)\ge t_j)=\mathbb P(\forall j=1,...,J,\ \Delta^{-1}(t_j)\le s_j)$.
This implies the convergence of the finite dimensional
distributions of the sequence of continuous processes 
$\left(t\mapsto  \mathcal S_n(t)=n^{-\alpha}(S_{\lfloor nt\rfloor}\bI_M+(nt-\lfloor nt\rfloor) \bI_M\circ \mathcal F^{\lfloor nt \rfloor})\right)_n$.

Observe that, for all $n\ge 1$, $t\mapsto  \mathcal S_n(t)$ is increasing. Thus, for 
all $T>0$, the sequence of processes $\left(\left(\mathcal S_n(t)\right)_{t\in[0;T]}\right)_{n\ge 1}$ is tight for the $M_1$-metric. 
We conclude that it converges to $(\Delta^{-1}(t))_{t\in[0;T]}$ for the $M_1$-metric. 
Since the limit process $\Delta^{-1}$  is continuous,
we conclude that the convergence in distribution holds also in $\mathcal C^0([0;T])$. 
\end{proof}
\begin{corollary}\label{coro2cvgceloi}
Assume Assumptions of Proposition~\ref{corocvgceloi} and that $(\cM,\cF,\mu)$ is ergodic, then, 
 due to the Hopf ergodic theorem, for any $f\in L^1(\mu)$, 
\beq 
 \left(S_{\lfloor nt\rfloor}f/n^\alpha\right)_{t\in[0,T]}\Longrightarrow \left(\frac{\int_{\cM}f\, d\mu}{\mu(M)}\, \Delta^{-1}(t)\right)_{t\in[0;T]}\, ,
\eeq
\end{corollary}
\begin{proof}
The proof is based on the combination of Proposition~\ref{corocvgceloi} (convergence in distribution of $((S_{\lfloor nt\rfloor}\bI_M/n^\alpha)_t)_n$ to a continuous process vanishing at 0) and of the Hopf ergodic theorem.
We will use the decomposition 
$\frac{S_{\lfloor nt\rfloor}f}{n^\alpha}= 
\frac{S_{\lfloor nt\rfloor}f}{S_{\lfloor nt\rfloor}\bI_M}\frac{S_{\lfloor nt\rfloor}\bI_M}{n^\alpha}$. 
It follows from the Hopf ergodic theorem that 
$(\frac{S_{m}f}{S_{m}\bI_M})_{m\ge 1}$ converges almost everywhere to $\frac{\int_{\mathcal M}f\, d\mu}{\mu(M)}$. 
This implies that
$$
\sup_{t\in [\frac{\log n}n,T]}\left\vert\frac{S_{\lfloor nt\rfloor}f}{S_{\lfloor nt\rfloor}\bI_M}-\frac{\int_{\mathcal M}f\, d\mu}{\mu(M)}\right\vert\longrightarrow 0,\,\,\ \ \mu-a.e.\, ,
$$
and that
$$
\sup_{t\in [0,\frac{\log n}n]}\left\vert\frac{S_{\lfloor nt\rfloor}f}{n^\alpha}\right\vert\le \frac{S_{\log n}|f|}{n^\alpha}=\frac{S_{\log n}|f|}{(\log n)^\alpha}\frac{(\log n)^\alpha}{n^\alpha}\longrightarrow 0,\,\,\,\ \ \mu-a.e.
$$
In particular
$$\sup_{t\in[0,\frac{\log n}n]}\left|\frac{S_{\lfloor nt\rfloor}\bI_M}{n^\alpha}\right|\longrightarrow 0\ \ \mu-a.e\, .$$
Furthermore, it follows from Proposition~\ref{corocvgceloi} that
$$\sup_{t\in[0,T]}\left|\frac{S_{\lfloor nt\rfloor}\bI_M}{n^\alpha}\right|\Longrightarrow \sup_{t\in[0,T]}\left|\Delta^{-1}(t)\right|\, .$$
Therefore, gathering our estimates on $[0,\frac{\log n}n]$ and on
$[\frac{\log n}n,T]$, we obtain that
$$
\sup_{t\in [0,T]}\left\vert\frac{S_{\lfloor nt\rfloor}f}{n^\alpha}-\frac{\int_{\mathcal M}f\, d\mu}{\mu(M)} \frac{S_{\lfloor nt\rfloor}\bI_M}{n^\alpha}\right\vert\Longrightarrow 0\, ,
$$
and we conclude using Proposition~\ref{corocvgceloi}.
\end{proof}
\begin{remark}
Note that, if the measure $\mu$ were finite and ergodic, then $R$ would be integrable and the Birkhoff Theorem would ensure the almost sure convergence of $(\varepsilon\mathcal R(s/\varepsilon))_{s\in[0;T]}$ to $(s\mathbb E_{\mu_M}[R])_{s\in[0;T]} = (\frac{s}{\mu(M)})_{s\in[0;T]}$, instead of the convergence in distribution \eqref{cvgcedistrib} that arises in the infinite-measure case. In that case, Corollary \ref{coro2cvgceloi} would amount to a functional Law of Large Numbers, namely, $\left( S_{\lfloor nt \rfloor} f / n \right)_{t \in [0,T]} \longrightarrow \left( t \int_\cM f \, d\mu \right)_{t \in [0,T]}$ almost surely (relative to the sup norm).
\end{remark}
\subsection{Functional Central Limit Theorem}\label{secF2}
Recall that we have set 
$\mathcal R(u):=
\sum_{k=0}^{\lfloor u
\rfloor-1} R\circ F^k$. For $f:\cM\rightarrow\mathbb R$, we also set
$$
\mathcal S(u):=\int_0^u\widetilde f\circ F^{\lfloor s\rfloor}\, ds=\sum_{k=0}^{\lfloor u\rfloor-1}
\widetilde f\circ F^k+(u-\lfloor u\rfloor)
\widetilde f\circ F^{\lfloor u\rfloor}\, ,$$
with $\widetilde f:M\rightarrow \mathbb R$
given by $\widetilde f(x):=\sum_{k=0}^{R(x)-1}f\circ \cF^k(x)$.
\begin{proposition}\label{prop:CLT}
Let $\gamma>0$. 
Assume Hypotheses~\ref{hypinduction} (Induction context) and~\ref{hypDelta} (hypothesis on $\alpha$ and $\Delta$). 
Let $f:\cM \rightarrow \mathbb R$ be a measurable function.
Assume that the
following joint convergence in distribution holds true
with respect to $\mu_M$:
\begin{equation}\label{cvgcedistrib2}
\left(\left(\varepsilon^{\frac 1\gamma}\mathcal S(t/\varepsilon)\right)_{t\geq 0},\left(
\varepsilon^{\frac 1\alpha}\mathcal R(s/\varepsilon
)\right)_{s\geq 0}\right)\stackrel{distrib}{\longrightarrow}\left(\left(B_t\right)_{t\geq 0},\left(\Delta(s)\right)_{s\geq 0}\right)\ 
\quad\mbox{as }\varepsilon\rightarrow 0\, ,
\end{equation}
where $B$ is a continuous process independent of $\Delta$, 
for respectively the uniform metric on $\mathcal C([0,T])$ and the Skorohod  $M_1$-metric on $\mathcal D([0;T])$ for any $T>0$.
Assume moreover that 
$\left(\sup_{k=0,...,N-1} \frac{g_f\circ F^k }{N^{\frac 1\gamma}}\right)_{N\geq 1}$ converges in probability (with respect to $\mu_M$) to 0, where we set $$g_f(x):=\sup_{\ell=0,...,R(x)-1}\left|\sum_{k=0}^{\ell}f( \cF^k(x))\right|.$$
Then
\begin{equation*}
\left(\frac{\sum_{k=0}^{\lfloor nt\rfloor-1}f\circ\mathcal  F^k}{n^{\frac \alpha \gamma}}
\right)_{t\geq 0}
\stackrel{distrib}{\longrightarrow} \left(B_{\Delta^{-1}(t)} \right)_{t\geq 0}\, \quad\mbox{as }n\rightarrow +\infty\, ,
\end{equation*}
with respect to $\mu_M$, for the uniform metric on $[0;T]$ for all $T>0$.

If, moreover, $(\cM,\cF,\mu)$ is ergodic, then this convergence in distribution holds in the strong sense with respect to $\mu$, i.e.,
with respect to any $\mu$-absolutely continuous probability measure. 
\end{proposition}
We point out that~\cite[Theorem 1.7]{Thomine1}
obtained an analogous conclusion, by using a stronger criterion, which implies the crucial assumption~\eqref{cvgcedistrib2} of Proposition~\ref{prop:CLT}. 
\begin{remark}
This result is motivated by the study of the case when $\int_{\mathcal M}f\, d\mu=0$.
In this case, it can be useful to recall that, 
if $M$ is a sweep-out set for $\cF$, that is, $\bigcup_{k\in\N} \cF^{-k} M = \cM$ mod $\mu$, then $\int_\cM f \, d\mu = \int_M \widetilde f \, d\mu_M = 0$. 
This equality of integrals holds true if the system $(\mathcal M,\mathcal F,\mu)$ is recurrent ergodic
\end{remark}

\begin{remark}
 The additional condition on $g_f$ is satisfied if there exist $C,\eta_0>0$ and $h\in L^{\gamma+\eta_0}(\mu_M)$ such that $g_f\le C\left|\widetilde f\right|+h$. Indeed, this follows from the tightness of $((\varepsilon^{\frac 1\gamma}\cS(t/\varepsilon))_t)_\varepsilon$ and from the Markov inequality that ensures that
$$
\mu_M\left(\sup_{k=0,...,N-1} h\circ F^k>N^{\frac 1\gamma}\eta\right)
\le N\mu_M\left( h>N^{\frac 1\gamma}\eta\right)\le N\frac{\Vert h\Vert_{L^{\gamma+\eta_0}(\mu_M)}^{\gamma+\eta_0}}{(\eta N^{\frac 1\gamma})^{\gamma+\eta_0}}=\mathcal O(N^{-\frac{\eta_0}\gamma})\, .
$$
In particular, this condition is satisfied if $f$ is supported on $M$.
\end{remark}
\begin{remark}
If the measure $\mu$ was finite and ergodic, then  $(\varepsilon\mathcal R(s/\varepsilon))_{s\geq 0}$ would converge almost surely to $(s\mathbb E_{\mu_M}[R])$, or equivalently $(\frac{s}{\mu(M)})_{s\geq 0}$ as $\varepsilon\rightarrow 0$. In this case, a standard time change argument would be enough to infer the Limit Theorem for $((S_{\lfloor nt\rfloor}f/n^{\frac 1\gamma})_{t\geq 0})_{n\ge 1}$ from the convergence in distribution of $\left((\varepsilon^{\frac 1\gamma}\mathcal S(t/\varepsilon))_{t\geq 0}\right)_\varepsilon$. The difficulty in the infinite measure case lies in the fact that
the best we can hope for the families of processes $\left(\left(\varepsilon^{\frac 1\gamma}\mathcal S(t/\varepsilon)\right)_t\right)_\varepsilon$ and $\left(\left(
\varepsilon^{\frac 1\alpha}\mathcal R(s/\varepsilon
)\right)_s\right)_\varepsilon$ is a convergence in distribution (and not an almost everywhere convergence). But having these two separate convergences in distribution would not be enough to conclude (since the limit processes are not deterministic).
\end{remark}
\begin{proof}[Proof of Proposition~\ref{prop:CLT}]
Due to \cite[Theorem 1]{Zweimuller}, when $(\cM,\cF,\mu)$ is ergodic, the convergence with respect to $\mu_M$ will imply the convergence with respect to any probability measure absolutely continuous with respect to $\mu$.
So, from now on, we study the convergence with respect to $\mu_M$.

We first explain why, up to replacing $f$ by $\widetilde f\bI_M$, we can restrict our study to the case where $f$ is supported on $M$.  
Note first that 
changing $f$ into $\widetilde f\bI _M$ does not change the process $\mathcal S$.
Furthermore, on $M$,
$$ S_n(\widetilde f\bI_M-f)(x)=h_f\circ \cF^{n}(x),\quad \mbox{with }h_f(y):=\sum_{k=0}^{R(y)-1}f\circ\cF^{k}(y)\bI_{\cM\setminus M}(y)\, ,$$
up to extend the definition of $R$ to $\cM$ ($R$ is the first visit time to $M$). We infer that, on $M$,
\begin{equation}\label{:excnegl}
\sup_{t\in[0;T]}\left|\frac{\sum_{k=0}^{\lfloor nt\rfloor-1}(f-\widetilde f\bI_M)\circ \cF^k}{n^{\frac \alpha 2}}\right|=   \sup_{t\in[0;T]}\frac{|h_f|\circ \cF^{\lfloor nt\rfloor}}{n^{\frac \alpha 2}}
\end{equation}
and
\begin{equation*}
    \sup_{\ell=0,...,R(x)-1}\left|\sum_{k=\ell}^{R(x)-1}f\circ  \cF^k(x)\right|\le 2g_f(x)\, .
\end{equation*}
Thus, for all $K>0$ and $\eta>0$,
\begin{align*}
& \limsup_{n\rightarrow +\infty}\mu_M\left(\sup_{k=0,...,\lfloor nT\rfloor} \frac{|h_f|\circ\cF^k }{n^{\frac\alpha 2}}>\eta\right)\\
 & \quad\quad\le  \limsup_{n\rightarrow +\infty}\mu_M\left(\mathcal R(K n^{\alpha})<nT\right)+ \limsup_{n\rightarrow +\infty}\mu_M\left(
\sup_{k=0,...,K n^\alpha} \frac{2g_f\circ F^k }{n^{\frac\alpha 2}}>\eta \right)\\
& \quad\quad\le \mathbb P(\Delta(K)\le T^{\frac 1\alpha})
 +0\, .
\end{align*}
Therefore~\eqref{:excnegl} converges in probability to 0 as $n\rightarrow +\infty$. 
Thus it is enough to prove the result for $f=\widetilde f \bI_M$.

We assume from now on that $f=\widetilde f \bI_M$ (i.e. that $f$ is supported on $M$), so that $\widetilde f=f_{|M}$. 
Let $T'>0$ and $T>0$.
According to the Skorohod
representation theorem, there exists 
a sequence of processes
$\left((\mathfrak a_n(t),\mathfrak b_n(t))_{t\in[0;T]}\right)_{n\ge 1}$ converging (everywhere) to a process $(\mathcal A(t),\mathcal B_t)_{t\in[0;T]}$
(with respect respectively to the M1-metric on $\mathcal D([0,T])$ and to the uniform metric on $\mathcal C([0;T])$) and such that
\begin{itemize}
\item $(\mathcal A(t),\mathcal B_t)_{t\in[0;T]}$ has the same distribution as $(\Delta(t),B_t)_{t\in[0;T]}$,
\item for all $n\ge 1$, 
$(\mathfrak a_n(t),\mathfrak b_n(t))_{t\in[0;T]}$ has the same distribution as
$\left(n^{-1}\mathcal R(n^\alpha t), n^{-\frac\alpha \gamma}\mathcal S(n^\alpha t)\right)_{t\in[0;T]}$.
\end{itemize}
Since $(\mathfrak a_n)_n$ converges to $\mathcal A$ with respect to the $M1$ metric, there exist two sequences of couples of  
continuous maps 
$\left(u\mapsto(t_{\infty,n}(u), x_{\infty,n}(u))\right)_{n\geq 1}$
 and 
$(u\mapsto(t_n(u), x_n(u)))_{n\geq 1}$ 
such that 
\begin{itemize}
\item $t_n$ and $t_{\infty,n}$ are non decreasing,
\item 
$t_{\infty,n}(0)=x_{\infty,n}(0)=0$ and $t_{\infty,n}(1)=T$, $x_{\infty,n}(1)=\mathcal A(T)$,  $x_{\infty,n}(u)\in[\mathcal A(t_{\infty,n}(u)-),\mathcal A(t_{\infty,n}(u))]$ 
\item $t_n(0)=x_n(0)=0$ and $t_n(1)=T$, $x_n(1)=\mathfrak a_n(T)$
such that $x_n(u)\in[ \mathfrak a_n(t_n(u)-), \mathfrak a_n(t_n(u))]$.
\item 
$(t_n(u)-t_{\infty,n}(u),x_n(u)-x_{\infty,n}(u))_{n\geq 1}$ converges to  $(0,0)$, uniformly
in $u\in[0;1]$, as $n$ goes to infinity.
\end{itemize}
Let $\Omega_{T,T'}$ be the event $\{\mathcal A(T)> T'\}$.
Let us prove that, on $\Omega_{T,T'}$,
$(\mathfrak b_n\circ \mathfrak a_n^{-1})_{n\geq 1}
$ converges uniformly on $[0,T']$ to 
$\mathcal B_{\mathcal A^{-1}}=(\mathcal B_{\mathcal A^{-1}(s)})_{s\in [0,T']}$. 
Since $\mathcal A$ is strictly increasing, 
$t_{\infty,n}(u)=
\mathcal A^{-1}(x_{\infty,n}(u))$. 
If $\mathcal A(T)>T'$, then, there exists $n_0$ such that, for every integer $n\ge n_0$, $x_n(1)=\mathfrak a_n(T)>T'$
and so, for all $n\ge n_0$ and all $s\in[0;T']$,  $x_n^{-1}(s)=u\in[0;1]$ is the smallest $u>0$ such that $x_n(u)=s$. 
Thus, on $\Omega_{T,T'}$, 
\begin{align}
\lim_{n\rightarrow +\infty}\sup_{s\in[0;T']}\left| x_{\infty,n}(x_n^{-1}(s))-s\right|&= \lim_{n\rightarrow +\infty}\sup_{s\in[0;T']}\left| x_{\infty,n}(x_n^{-1}(s))-x_{n}(x_n^{-1}(s))\right|\nonumber\\
&\le \lim_{n\rightarrow +\infty}\left\Vert x_{\infty,n}-x_n\right\Vert_\infty=0\, .\label{zzz1}
\end{align}
But we also know that $t_n(x_n^{-1}(s))-t_{\infty,n}(x_n^{-1}(s))$ converges to 0 uniformly in $s\in[0, T']$, which can be rewritten as follows since $\mathcal A$ is strictly increasing:
\beq\label{zzz2}
\lim_{n\rightarrow +\infty}\sup_{s\in[0,T']}\left|t_n(x_n^{-1}(s))-\mathcal A^{-1}(x_{\infty,n}(x_n^{-1}(s)))\right| =0\quad\mbox{on }\Omega_{T,T'}\, .
\eeq
Since $\mathcal A^{-1}$ is continuous, it follows from \eqref{zzz1} and \eqref{zzz2} that
$$
\lim_{n\rightarrow +\infty}\sup_{s\in[0;T']}\left|t_n(x_n^{-1}(s))-\mathcal A^{-1}(s)\right| =0\, \quad\mbox{on }\Omega_{T,T'}.
$$
Since $(\mathfrak b_n)_n$ is a sequence of continuous functions converging uniformly  on $[0;T]$ to the continuous function $\mathcal B$, we conclude that
\beq
\lim_{n\rightarrow +\infty}\sup_{[0;T']}
\left|\mathfrak b_n(t_n\circ x_n^{-1})-\mathcal B_{\mathcal A^{-1}}\right|=0\quad \mbox{on}\quad \Omega_{T,T'}\, .\label{zzz2b}
\eeq
But, on $\Omega_{T,T'}$, for all $s\in[0;T']$,  $$\mathfrak a_n(t_n(x_n^{-1}(s))-)\le s=x_n(x_n^{-1}(s))\le \mathfrak a_n(t_n(x_n^{-1}(s))),$$ and so $|t_n(x_n^{-1}(s))-\mathfrak a_n^{-1}(s)|\le n^{-\alpha}$. Thus, on $\Omega_{T,T'}$, 
\beq
\sup_{[0;T']}\left|\mathfrak b_n(t_n\circ x_n^{-1})-\mathfrak b_n\circ \mathfrak a_n^{-1}\right|\le \sup_{t,t'\in[0;T]:|t-t'|\le n^{-\alpha}}|\mathfrak b_n(t)-\mathfrak b_n(t')|\, ,\label{zzz2c}
\eeq
which converges to 0 as $n\rightarrow +\infty$. 
Hence, it follows from~\eqref{zzz2b} and~\eqref{zzz2c} that, on the event $\Omega_{T,T'}$, $\mathfrak b_n\circ \mathfrak a_n^{-1}$ converges uniformly on $[0;T']$ of  to $ \mathcal B_{\mathcal A^{-1}}$ and so, we will infer 
the announced convergence in distribution since
$\mathcal B_{\mathcal A^{-1}}$ has the same distribution as $B_{\Delta^{-1}}$
and since $\mathfrak b_n\circ \mathfrak a_n^{-1}$ has the same distribution as $$\widetilde S_{n}:t\mapsto n^{-\frac \alpha \gamma}\sum_{k=0}^{\mathcal R^{-1}(nt)-1}
\widetilde
f\circ  F^k= n^{-\frac\alpha\gamma}\sum_{k=0}^{\lceil nt\rceil-1}f\circ\mathcal F^k\, ,$$
with $\mathcal R^{-1}(s):=\inf\{u > 0\, :\,  \mathcal R(u) \ge  s\}$ equal to 1 plus the number of returns of $(\mathcal F^k)_k$ to $M$ before time $\lceil s\rceil$.
Indeed, for any $h:\mathcal D([0,T'])\rightarrow \mathbb R$ bounded continuous with respect to the uniform norm, 
\begin{align*}
\limsup_{n\rightarrow +\infty}\left|\mathbb E_{\mu_M}\left[h\left(\widetilde S_n\right)\right]-
\mathbb E\left[h(B_{\Delta^{-1}})\right]\right|
&\le 2\Vert h\Vert_\infty\left(1-\mathbb P(\Omega_{T,T'})\right)\\
&\quad+\limsup_{n\rightarrow +\infty}\left|\mathbb E_{\mu_M}\left[\bI_{\Omega_{T,T'}}\left(h\left(\mathfrak b_n\circ \mathfrak a_n^{-1}\right)-h(\mathcal B_{\cA^{-1}})\right)\right]\right|\\
&\le 
 2\Vert h\Vert_\infty\mathbb P\left(\Delta_T\le T'\right)\, ,
\end{align*}
which goes to 0
as $T\rightarrow +\infty$.
Thus we have proved the convergence in distribution of $(n^{-\frac \alpha \gamma}\sum_{k=0}^{\lfloor nt\rfloor-1} f\circ \mathcal F^k)_{t\in[0;T]}$ to $(B_{\Delta^{-1}(t)})_{t\in[0;T]}$, with respect to $\mu_M$, as $n\rightarrow +\infty$.
\end{proof}

\subsection{Joint limit theorem via point processes}\label{sec:F3}
A classical strategy to prove a convergence result to the $\alpha$-stable L\'evy process $\Delta$ consists in proving the convergence of a family of point processes to a Poisson point process. In the context of dynamical systems, this has been enlighted in \cite{TyranKaminska1,TyranKaminska2} where necessary and sufficient conditions have been proved.
In particular, the convergence in distribution \eqref{cvgcedistrib} of the family of processes $\left(\left(\varepsilon^{\frac 1\alpha}\mathcal R(s/\varepsilon)\right)_s\right)_\varepsilon$ to $(\Delta(s))_s$ required in Proposition~\ref{corocvgceloi} follows from the convergence to a Poisson process of the family $(\mathcal N_\varepsilon)_\varepsilon$ of point processes given by:
\beq \label{eq:DefNeps}
\mathcal N_\varepsilon=\sum_{j\ge 1}\delta_{\left(j\varepsilon,
\varepsilon^
{\frac 1\alpha}R\circ F^j\right)}\mbox{ on }[0,+\infty)\times(0,+\infty)\, .
\eeq
\begin{proposition} \label{doublecvgce2}
Let $\alpha\in(0,1)$ and $\gamma>0$. 
Assume Hypotheses~\ref{hypinduction} (Induction context) and~\ref{hypDelta} (hypothesis on $\alpha$ and $\Delta$). Let $\widetilde f:M\rightarrow \mathbb R$ be a measurable function.
Assume that $\mu_M(R>N)\sim c_0^\alpha N^{-\alpha}$ as $N\rightarrow +\infty$ and that 
\beq\label{cvdistribBP}
\left(\left(  \left(\varepsilon^{\frac 1\gamma}
\mathcal S(t/\varepsilon)\right)_{t\in[0,T]}, \mathcal N_\varepsilon \right)\right)_\varepsilon\stackrel{distribution}{\longrightarrow}
(B,\mathcal P)\eeq
with respect to $\mu_M$, for respectively the uniform topology and the vague convergence topology\footnote{The vague convergence of a sequence of measures $(\mathfrak m_n)_n$ on a topological space $V$ to $\mathfrak m$ means that, for every $h:V\rightarrow \mathbb R$ continuous with compact support, $\lim_{n\rightarrow +\infty}\int_Vh\, d\mathfrak m_n=\int_Vh\, d\mathfrak m$. It is equivalent to the convergence of 
$\mathfrak m_n(E)$ to $\mathfrak m(E)$ for every relatively compact set $E\subset V$ such that $\mathfrak m(\partial E)=0$.}, as $\varepsilon\rightarrow 0$, 
where $B$ is a continuous process and where
$\mathcal P$ is a Poisson process
on $[0,+\infty)\times(0,+\infty)$ with density $(u,z)\mapsto c_0^\alpha\alpha z^{-\alpha-1}\bI_{\{z>0\}} $, independent one from the other. Then Assumption \eqref{cvgcedistrib2} (joint convergence assumption) of Proposition~\ref{prop:CLT} holds true, with $B$ and $\Delta$
independent one from the other (with $\Delta$ as in hypothesis~\ref{hypDelta}).\\
As a consequence, the FCLT stated in Proposition~\ref{prop:CLT} holds true for $f:=\widetilde f\bI_M$. 
\end{proposition}
\begin{proof}
By assumption, we already know that the first coordinate of~\eqref{cvgcedistrib2} converges
in distribution. Let us prove the convergence of the second coordinate knowing the first one. The tightness of this family of conditioned processes will follow from the monotonicity of the processes. Thus it is enough to prove the convergence of the conditional finite dimensional distributions.
Let $J$ be a positive integer, $s_0=0<s_1<...<s_J$ and  $t=(t_1,...,t_J)\in\mathbb R^J$. 
Set $Z_\varepsilon:=\left(\varepsilon^{\frac 1\gamma}
\mathcal S(
t/\varepsilon
)
\right)_{t\in[0,T]}$ and let $U$ be an open set in $\mathcal C^0([0;T])$ such that $\mathbb P(B\in \partial U)=0$ and such that $\mathbb P(B\in U)\ne 0$. 
A crucial fact is that
$$
V_\varepsilon:=\left(\varepsilon^{\frac 1\alpha}\sum_{k=0} ^{\lfloor s_j/\varepsilon\rfloor-1}R\circ F^k\right)_{j=1,...,J}=\left(\int_{[0,s_j)\times \mathbb (0,+\infty)}w\, 
d {\mathcal N}_\varepsilon(u,w)\right)_{j=1,...,J}\, .$$
We will prove that, for every $\delta>0$, the family of random variables
$$
W_\varepsilon(\delta):=\left(\varepsilon^{\frac 1\alpha}\sum_{k=0} ^{\lfloor s_j/\varepsilon\rfloor-1}R\circ F^k \cdot \bI_{\{\varepsilon^{\frac 1\alpha}R\circ F^k> \delta\}}\right)_{j=1,...,J}=\left(\int_{[0,s_j)\times(\delta,+\infty)}w\, 
d {\mathcal N}_\varepsilon(u,w)\right)_{j=1,...,J}$$
converges in distribution (with respect to $\mu_M$ knowing $\{Z_\varepsilon\in U\}$), as $\varepsilon\rightarrow 0$, to
$$W(\delta):=\left(\int_{[0,s_j)\times(\delta,+\infty)}w\, 
d {\mathcal P}(u,w)\right)_{j=1,...,J}\, ,
$$
i.e.
\begin{equation}\label{AAA3}
\forall \delta>0,\quad \forall J\in\mathbb N^*,\quad \forall t\in\mathbb R^J,\quad  
\lim_{\varepsilon\rightarrow 0}
\left|\mathbb E_{\mu_M}[e^{i\langle t,W_\varepsilon(\delta)\rangle}|Z_\varepsilon\in U]
- \mathbb  E[e^{i\langle t,W(\delta)\rangle}]\right|=0.
\end{equation}
Moreover, we will prove that, for every $\gamma_0>0$, we have
\begin{equation}\label{PPP2}
\lim_{\delta\rightarrow 0}\limsup_{\varepsilon\rightarrow 0}
   \mu_M\left(|T_\varepsilon(\delta)|>\gamma_0|Z_\varepsilon\in U\right)=0 \, ,
\end{equation}
with
\beq  \label{Noo}
T_\varepsilon (\delta):=V_\varepsilon-W_\varepsilon(\delta)=\left(\varepsilon^{\frac 1\alpha}\sum_{k=0}^{\lfloor s_j/\varepsilon\rfloor-1} R\circ F ^k\bI_{\{\varepsilon^{\frac 1\alpha} R\circ F^k\le \delta\}}\right)_{j=1,...,J}\, .
\eeq
We will also prove that
\begin{equation}\label{PPPP3}
\lim_{\delta\rightarrow 0}\left|{\mathbb E}_{\mu_M}\left[e^{i\langle t,W (\delta)\rangle}
\right]-\prod_{j=1}^J\Psi_{\alpha,c_0}\left((s_j-s_{j-1})^{\frac 1\alpha}\sum_{k=j}^Jt_k\right)\right|=0 \, 
\end{equation}
recalling that $\Psi_{\alpha,c_0}$ given by
$\Psi_{\alpha,c_0}(t):=e^{-|c_0t|^\alpha\Gamma(1-\alpha)e^{-\frac{i\pi\alpha}2\sign(t)}}$ is
the characteristic function of $\Delta(1)$.
Assuming for the moment~\eqref{AAA3},~\eqref{PPP2} and~\eqref{PPPP3},
let us see how we can prove~\eqref{cvgcedistrib2}. Let $t=(t_1,...,t_J)\in\mathbb R^J$. 
Let $\eta>0$. Due to \eqref{PPP2} and \eqref{PPPP3}, we consider $\delta>0$ and $\varepsilon_0>0$ such that
\begin{equation}\label{AAA1}
\forall \varepsilon\in(0,\varepsilon_0],\quad
\mu_M\left(|\langle t,T_\varepsilon(\delta)\rangle|>\tfrac{\eta}6\Big|Z_\varepsilon\in U\right)<\tfrac\eta 6
\end{equation}
and
\begin{equation}\label{AAA2}
\left|{\mathbb E}_{\mu_M}\left[e^{i\langle t,W (\delta)\rangle}
\right]-\prod_{j=1}^J\Psi_{\alpha,c_0}\left((s_j-s_{j-1})^{\frac 1\alpha}\sum_{k=j}^Jt_k\right)\right|<\tfrac\eta 2 \, .
\end{equation}
Therefore, we have
\begin{eqnarray*}
 && \limsup_{\varepsilon\rightarrow 0}\displaystyle\left|\mathbb E_{\mu_M} [e^{i\langle t,V_\varepsilon\rangle}|Z_\varepsilon\in U] -\prod_{j=1}^J\Psi_{\alpha,c_0}\left((s_j-s_{j-1})^{\frac 1\alpha}\sum_{k=j}^Jt_k\right)\right| \\
&\le&\tfrac\eta 2+\limsup_{\varepsilon\rightarrow 0}\left|\mathbb E_{\mu_M}[e^{i\langle t,V_\varepsilon\rangle}|Z_\varepsilon\in U]-\mathbb E [e^{i\langle t,W(\delta)\rangle}]\right|\ \ \ \mbox{due to \eqref{AAA2}}\\
&\le& \tfrac{\eta} 2+\limsup_{\varepsilon\rightarrow 0}\left|\mathbb E_{\mu_M}[e^{i\langle t,(W_\varepsilon(\delta)+T_\varepsilon(\delta))\rangle}|Z_\varepsilon\in U]-\mathbb E
[e^{i\langle t,W(\delta)\rangle}]\right|\ \ \ \mbox{due to \eqref{Noo}}\\
&\le& \tfrac{\eta} 2+\limsup_{\varepsilon\rightarrow 0}\left|\mathbb E_{\mu_M}[e^{i\langle t,(W_\varepsilon(\delta)+T_\varepsilon(\delta))\rangle}
-e^{i\langle t,W_\varepsilon(\delta)\rangle}|Z_\varepsilon\in U]\right|\ \ \ \mbox{due to \eqref{AAA3}}\\
&\le& \eta\ \ \ \mbox{due to \eqref{AAA1}}.
\end{eqnarray*}
This inequality being true for any $\eta>0$, we conclude that the above $\limsup$
is null.\\
Therefore, we have proved that
\beq
\lim_{\varepsilon\rightarrow 0}\mathbb E_{\mu_M} [e^{i\langle t,V_\varepsilon\rangle}|Z_\varepsilon\in U] =\prod_{j=1}^J\Psi_{\alpha,c_0}\left((s_j-s_{j-1})^{\frac 1\alpha}\sum_{k=j}^Jt_k\right) \, .
\eeq
Recall that $\Delta$ is $1/\alpha$-stable with independent increments and that $\Psi_{c,\alpha_0}$ is the characteristic function of 
$\Delta(1)$. Thus
\beq
\mathbb E\left[e^{i\sum_{j=1}^Jt_j\Delta(s_j)}\right]
=\mathbb E\left[ \prod_{j=1}^J e^{i\sum_{k=j}^Jt_k(\Delta(s_j)-\Delta(s_{j-1}))}\right]
=\prod_{j=1}^J\Psi_{\alpha,c_0}\left((s_j-s_{j-1})^{\frac 1\alpha}\sum_{k=j}^Jt_k\right)\, .
\eeq
This ends the proof of the convergence of the finite dimensional distributions of
$(V_\varepsilon|Z_\varepsilon\in U)_{\varepsilon}$ to those of
$\Delta$
as $\varepsilon$ goes to 0.
Since $(Z_\varepsilon)_\varepsilon$ converges in distribution
to $B$ and since $(V_\varepsilon)_\varepsilon$ is tight in $M_1$, 
we conclude the convergence in distribution of the pair of processes.
\end{proof}

\begin{proof}[Proof of \eqref{AAA3} (Conditional convergence in distribution of $(W_\varepsilon(\delta))_\varepsilon$)]
Let $K>\delta$.
It follows from~\eqref{cvdistribBP} that
$(\mathcal N_\varepsilon|Z_\varepsilon\in U)_{\varepsilon}$ converges in distribution, as $\varepsilon\rightarrow 0$, to $(\mathcal P|B\in U)$, i.e. to $\mathcal P$.
Note that 
$$G:\mathfrak m\mapsto\left(\int_{(0,s_j)\times (\delta,K]}w\, d\mathfrak m(u,w)\right)_{j=1,...,J}$$
is continuous 
(for the vague topology on the set of non negative Radon measures on $[0,+\infty)\times(0,+\infty)$) 
at $\mathfrak m=\mathcal P$ except if $\mathcal P(\partial((0,s_j)\times(\delta,K)))>0$,
which holds true with null probability. 
It follows from the continuous mapping theorem 
that $$(G(\mathcal N_\varepsilon )=W_\varepsilon(\delta)-W_\varepsilon(K)|Z_\varepsilon\in U)_{\varepsilon>0}$$ converges in distribution to $G(\mathcal P)=W(\delta)-W(K)$. 
Thus, for any $t\in\mathbb R^J$,
\beq\label{diffWdeltaK}
\lim_{\varepsilon\rightarrow 0}\mathbb E_{\mu_M}[e^{i\langle t,W_\varepsilon(\delta)-W_\varepsilon(K)\rangle}
|Z_\varepsilon\in U]=\mathbb E[e^{i\langle t,W(\delta)-W(K)\rangle}]
\eeq
and, for any  $\varepsilon>0$,
\begin{align}
\nonumber
\left|\mathbb E_{\mu_M}[e^{i\langle t,W_\varepsilon(\delta)\rangle}|Z_\varepsilon\in U]-\mathbb E[e^{i\langle t,W(\delta)\rangle}]\right|&\le\left|\mathbb E_{\mu_M}[e^{i\langle t,W_\varepsilon(\delta)-W_\varepsilon(K)\rangle}
|Z_\varepsilon\in U]-\mathbb E[e^{i\langle t,W(\delta)-W(K)\rangle}]\right|\\
&\quad+\mu_M(W_\varepsilon(K)\ne 0|Z_\varepsilon\in U)+\mathbb P(W(K)\ne 0)\, .\label{GGG1}
\end{align}
But $W(K)$ converges in probability\footnote{Recall that, when the limit is constant, convergence in distribution is equivalent to convergence in probability.} to $0$ as $K\rightarrow +\infty$. Indeed, $|W(K)|_\infty=(W(K))_J$ and
\begin{align}\nonumber
   \mathbb  E\left[\exp\left(it(W(K))_J\right)\right]
&=   \mathbb  E\left[\exp\left(it\int_{(0,s_J)\times (K,+\infty)} w \, d{\mathcal P}(u,w)\right)\right]\\
     &= \exp\left(c_0^\alpha s_J\int_K^\infty\alpha\frac{e^{itx}-1}{x^{\alpha+1}}\, dx\right)\longrightarrow 1\, ,\label{GGG2}
\end{align}
as $K\rightarrow +\infty$ (see~\cite[Proposition 3.6]{Resnick}). 
Furthermore, there exists $\varepsilon_1>0$ such that
\begin{equation} \label{GGG3}
 \sup_{\varepsilon\in(0,\varepsilon_1)} \mu_M\left(
   W_\varepsilon(K)\neq0|Z_\varepsilon\in U\right) \longrightarrow 0 ,
  \quad   \mbox{as }K\rightarrow\infty \, .
 \end{equation}
Indeed, consider $\varepsilon_1>0$ such that
$$\forall \varepsilon\in(0,\varepsilon_1)\, ,\quad
\mu_M(Z_\varepsilon\in U)>\mathbb P(B\in U)/2\, .
$$
For any $\varepsilon\in(0,\varepsilon_1)$, we have
\begin{eqnarray*}
  \mu_M\left(
    W_\varepsilon(K)
    \neq0|Z_\varepsilon\in U\right)
  &\leq &  \frac{\mu_M\left(\exists k\le s_J/\varepsilon\, : \, \varepsilon^{\frac 1\alpha}R\circ F^k>K\right)}{\mu(Z_\varepsilon\in U)} \\
   &\leq&  2\frac{\sum_{k\le s_J/\varepsilon} \mu_M\left(R>\varepsilon^{-\frac 1\alpha} K\right)}{\mathbb P(B\in U)} \\
    &\leq&  2\frac{ (s_J/\varepsilon)C\varepsilon K^{-\alpha}}{\mathbb P(B\in U)}\\ 
     &\leq&  2\frac{ s_JC K^{-\alpha}}{\mathbb P(B\in U)} 
\, ,
\end{eqnarray*}
implying
\eqref{GGG3}, setting $C:=\sup_{y>0}y^\alpha\mu_M(R>y)$.
This combined with \eqref{diffWdeltaK}, \eqref{GGG1} and~\eqref{GGG2} ends the proof of the convergence in distribution 
of $(W_\varepsilon(\delta)|Z_\varepsilon\in U)_\varepsilon$ to $W(\delta)$.
\end{proof}
\begin{proof}[Proof of \eqref{PPP2}]
Since $\lim_{\varepsilon\rightarrow 0}\mu_M(Z_\varepsilon\in U)=\mathbb P(B\in U)>0$, 
it is enough to prove that $\lim_{\delta\rightarrow 0}\limsup_{\varepsilon\rightarrow 0}\mathbb E_{\mu_M}[(T_\varepsilon(\delta))_J]=0$. This follows from the following computation
\begin{eqnarray*}
\mathbb E_{\mu_M}[(T_\varepsilon(\delta))_J]&\le&
s_J/\varepsilon \, \mathbb E_{\mu_M}\left[\varepsilon^{\frac 1\alpha}R\bI_{\{\varepsilon^{\frac 1\alpha}R\le \delta\}}\right]\\
&\le&s_J/\varepsilon
      \int_0^{\delta}\mu_M\left(\delta\ge \varepsilon^{\frac 1\alpha}R>z\right)\, dz\\
&\le&s_J/\varepsilon
      \int_0^{\delta}\mu_M\left(R>\varepsilon^{-\frac 1\alpha}z\right)\, dz\\
&\le& Cs_J/\varepsilon 
\int_0^{\delta}\varepsilon z^{-\alpha}\, dz\\
&\le& Cs_J\frac{\delta^{1-\alpha}}{1-\alpha}\longrightarrow 0\, ,\quad\mbox{as }\ \delta\rightarrow 0\, ,
\end{eqnarray*}
where we used again the notation $C:=\sup_{y>0}y^\alpha\mu_M(R>y)$.
\end{proof}

\begin{proof}[Proof of \eqref{PPPP3}]
Due to the exponential formula, we have
\begin{align*}
\mathbb E&\left[\prod_{j=1}^Je^{it_j\int_{(0,s_j)\times(\delta,+\infty)}w\, d\mathcal P(u,w)  }\right]=
\mathbb E\left[\exp\left(i\int_{(0,+\infty)\times(\delta,+\infty)}\left(
\sum_{j=1}^Jt_j\bI_{\{u\le s_j\}}\right)w\, d\mathcal P(u,w)\right)\right]\\
&=
\mathbb E\left[\exp\left(i\sum_{j=1}^J\left(\sum_{k=j}^Jt_k\right)
\int_{(s_{j-1},s_j]\times(\delta,+\infty)}w\, d\mathcal P(u,w)\right)\right]\\
&=
\mathbb E\left[\exp\left(i\sum_{j=1}^J\left(\sum_{k=j}^Jt_k\right)
\int_{(\delta,+\infty)}w\, d\mathcal P_j(w)\right)\right]\, ,
\end{align*}
where $\mathcal P_j$ are independent Poisson processes with respective intensities 
$z\mapsto (s_j-s_{j-1})c_0^\alpha \alpha z^{-\alpha-1}\bI_{\{z>0\}}$. 
Therefore, using e.g.~\cite[Proposition 3.6]{Resnick}, it follows that
\begin{align*}
\mathbb E&\left[\prod_{j=1}^Je^{it_j\int_{(0,s_j)\times(\delta,+\infty)}w\, d\mathcal P(u,w)  }\right]=\prod_{j=1}^J\exp\left(c_0^\alpha(s_j-s_{j-1})\int_\delta^{+\infty}\frac{e^{i\sum_{k=j}^Jt_kz}-1}{z^{\alpha+1}}
\alpha\, dz\right)\, .
\end{align*}
Due to \cite[p. 568]{Feller2}, we have
$$ \lim_{\delta\rightarrow 0}\int_\delta^{+\infty}\frac{e^{itz}-1}{z^{\alpha+1}}
\alpha\, dz=-|t|^\alpha\Gamma(1-\alpha)e^{-\frac{i\pi\alpha}2\sign(t)}\quad 
\, . $$
So
$\lim_{\delta\rightarrow 0}
\mathbb E\left[\prod_{j=1}^Je^{it_j\int_{(0,s_j)\times(\delta,+\infty)}w\, d\mathcal P(u,w)  }\right]=
\prod_{j=1}^J\Psi_{\alpha,c_0}\left((s_j-s_{j-1})^{\frac 1\alpha}\sum_{k=j}^Jt_k\right)$.
\end{proof}

\section{Limit theorems for the map}\label{proofformap}
Still as in the context of Section~\ref{secF} and with the notations therein, assuming in particular Hypothesis~\ref{hypinduction}, we will make the following additional assumptions that are satisfied in the context of billiards with infinite cusps considered in the present article.
\begin{hypothesis}\label{hypinductionbis}
Let $\alpha\in(0,1)$. 
We assume that $(M,F,\mu_M)$ is invertible and that
$M$ admits a finite partition in subsets $M^{[1]},...,M^{[I_0]}$, and we set $\mathcal I(x)=i$ if $x\in M^{[i]}$. We assume furthermore that
\begin{itemize}
    \item[(A)] For any $i=1,...,I_0$, there exists $c_i>0$ such that $\lim_{N\rightarrow +\infty}N^{\alpha}\mu_M(M^{[i]}\cap\{R>N\})= c_i^\alpha $ as $N\rightarrow +\infty$, and we set $c_0^\alpha:=\sum_{i=1}^{\mathcal I}c_i^\alpha$,
    \item[(B)]  For all $K_0>0$, there exists $\eta>0$ such that $\sup_{1\le k\le K_0\log N}\mu_M(R>N,\ R\circ F^k>N)=\mathcal O(N^{-\alpha-\eta})$,
    \item[(C)]  Let $\gamma\in(0;1)$, 
there exist two subsets $\mathcal E^\pm$ of $ M^2$, and constants $C>1$ and $\theta\in(0,1)$ such that for
any $f\in \cH^-_{\gamma}$, any $g\in \cH^+_\gamma$
and any $n \geq
 1$,
\beq  \label{expboundgene} 	
\left|\int_Mf\circ F^n.g\, d\mu_M-\int_Mf\, d\mu_M\int_Mg\, d\mu_M\right|\leq C \|f\|_{\gamma}^-\|g\|_{\gamma}^+ \theta^{n}\, .
\eeq 
where we set
$$\Vert g\Vert_\gamma^\pm:=\Vert g\Vert_\infty+[g]_\gamma^\pm,\quad [g]_\gamma^\pm:=\sup_{(x,y)\in \mathcal E^\pm}\frac{|g(x)-g(y)|}{\gamma^{\mathbf{s}_\pm(x,y)}}\, ,$$
where $\mathbf s_{\pm}:\mathcal E^{\pm}\rightarrow \mathbb N\cup\{+\infty\}$ define two separation times satisfying
$$(x,y)\in \mathcal E^{\pm}\quad\Rightarrow\quad
(F^{\mp k}(y), F^{\mp k}(x))\in\mathcal E^{\pm}\mbox{ and } 
\mathbf{s}_\pm(F^{\mp k}(x),F^{\mp k}(y))=\mathbf{s}_\pm(x,y)+k\, ,
$$
and where
\beq
(x,y)\in\mathcal E^-\quad \Rightarrow\quad \left[R(x)=R(y)\mbox{ and }\cI(x)=\cI(y)\right]\, ,
\eeq
\beq
(x,y)\in\mathcal E^+\quad\Rightarrow \quad \left[R(F^{-1}(x))=R(F^{-1}(y))\mbox{ and }\cI(F^{-1}(x))=\cI(F^{-1}(y))\right]\, .
\eeq
\end{itemize}
\end{hypothesis}
\begin{remark}\label{rk:hypsatisfied}
In the case of the billiard with infinite cusps, we consider a partition so that $M^{[i]}\supset M^{(i)}$ for all $i=1,...,\cJ$,
so that the set of configurations making excursions in the $i$-th cusp is contained in the set $M^{[i]}$.  
Condition~(A) is proved in Proposition~\ref{Mm2} with $c_i^\alpha=\frac{4 I_{1/\alpha}^\alpha}{\beta\mu(M)}\mathfrak c_i^{\frac 1{\beta+1}}$ 
and $c_0^\alpha=\frac{4 \mathfrak CI_{1/\alpha}^\alpha}{\beta\mu(M)}$ , Condition~(C) is proved in \eqref{correT}.
Condition~(B) is proved in Lemma~\ref{deco1}.
\end{remark}

Recall that 
$$\mathcal N_\varepsilon:=\sum_{j\ge 1}\delta_{\left(j\varepsilon,
\varepsilon^
{\frac 1\alpha}R\circ F^j\right)}\quad\mbox{ on }[0,+\infty)\times(0,+\infty)\, .$$
We define the marked point process
$${\mathcal N}^{\cI}_\varepsilon:=\sum_{j\ge 1}\delta_{\left(j\varepsilon,
\varepsilon^
{\frac 1\alpha}R\circ F^j,\cI\circ F^j\right)}\quad\mbox{ on }[0,+\infty)\times(0,+\infty)\times\{1,...,I_0\}\, .$$
Observe that $\mathcal N_\varepsilon$ is the image measure of $\mathcal N_\varepsilon^{\cI}$ by the projection on the two first coordinates.  Thus the convergence in distribution of  $(\mathcal N_\varepsilon^{\cI})_{\varepsilon>0}$ implies the convergence in distribution of $(\mathcal N_\varepsilon)_{\varepsilon>0}$.

\subsection{Convergence to Stable Processes  for maps and proof of Theorem~\ref{thm:LLN}}\label{proofLLN}
We start by stating general results under our general hypotheses and we will end this subsection with the proof of Theorem~\ref{thm:LLN}.
\subsubsection{A general result of Convergence to Stable Processes  for maps}
\begin{proposition}\label{propcvPP}
Let $\alpha\in(0,1)$. Assume Hypotheses~\ref{hypinduction},~\ref{hypDelta} and~\ref{hypinductionbis}. 
Then the family of point processes $(\mathcal N_\varepsilon^\cI)_\varepsilon$ converges in distribution (with respect to $\mu_M$), as $\varepsilon\rightarrow 0$, to a Poisson point process with intensity of density $(t,y,i)\mapsto c_i^\alpha\alpha |y|^{-\alpha-1}$. 
As a consequence, 
\begin{equation}\label{cvgcedistrib0}
\left(\varepsilon^{\frac 1\alpha}\mathcal R(\lfloor s/\varepsilon\rfloor)\right)_{s\geq 0}=\left(\varepsilon^{\frac 1\alpha}\sum_{k=0}^{\lfloor s/\varepsilon\rfloor-1} R\circ F^k\right)_{s\geq 0}\stackrel{distrib}{\longrightarrow}\left(\Delta(s)\right)_{s\geq 0}\, , \quad\mbox{as }\varepsilon\rightarrow 0\, ,
\end{equation}
with respect to $\mu_M$ in the sense of finite dimensional distribution and 
for the $M_1$-metric on $\mathcal D([0;T])$ for all $T>0$, where $\Delta$ 
is as in Hypothesis~\ref{hypDelta}.\\
In particular, the CSP for $(\mathcal M,\mathcal F,\mu)$ and $\bI_M$ stated in Proposition~\ref{corocvgceloi} holds true.     
\end{proposition}
\begin{proof}
Let $a>0$. We set $A_\varepsilon=\{R>a\varepsilon^{-\frac 1\alpha}\}$. It follows from Hypothesis~\ref{hypinductionbis} that $\mu_M(A_\varepsilon)\sim c_0^\alpha a^{-\alpha}\varepsilon$ as $\varepsilon\rightarrow 0$. 
Setting $\tau_{A_\varepsilon}(x):=\inf\{m\ge 1:\,  F^m(x)\in A_\varepsilon\}$, due \cite[Theorem 2.1]{BenoitSoaz} and
to \cite[Corollary 4.10]{JPZ}, it is enough to prove the existence of $p_\varepsilon=o(\varepsilon^{-1})$ such that
\begin{equation}\label{control1}
\mu_M(\tau_{A_\varepsilon} \le p_\varepsilon
|A_\varepsilon) = o(1)\, ,
\end{equation}
and such that,
for every $c>0$, every $K \ge 1$, and every choice of intervals $W_1=[a_1,b_1), ..., W_K=[a_K,b_K) \subset(0,+\infty)$,
\begin{equation}\label{control2}
\sup
_{A,B\subset A_\varepsilon\, :\, A\in\mathcal G_\varepsilon, B\in\sigma( \cup_{j\ge p_\varepsilon} F^{-j} (\mathcal G_\varepsilon))}
|\mu_M(B \cap A) - \mu_M(A)\mu_M(B)|= o(\varepsilon)\, ,
\end{equation}
with $\mathcal G_\varepsilon$ 
 the $\sigma$-algebra of $\{R\ge c\varepsilon^{-\frac 1\alpha}\}$ generated by the sets of the form 
$M^{[i']}\cap \{\varepsilon^{\frac 1\alpha}R\in W_i\}$.\\
Considering $\theta$ of Condition~(C) of Hypothesis~\ref{hypinductionbis},
we 
set $p_\varepsilon:=\lceil-K_0\log\varepsilon\rceil$, with $K_0>\frac 1{-\log\theta} $.
With this choice, \eqref{control1} will follow from the following estimate~:
$$
\mu_M(\tau_{A_\varepsilon} \le p_\varepsilon
|A_\varepsilon) \le 
\sum_{k=1}^{p_\varepsilon}
\frac{\mu_M(A_\varepsilon\cap F^{-k}(A_\varepsilon))}{\mu_M(A_\varepsilon)}=\mathcal O\left(
p_\varepsilon\frac{\varepsilon^{1+\frac\eta\alpha}}{\varepsilon}\right)=o(1)\, ,
$$
for $\eta>0$ as in Condition (B) of Hypothesis~\ref{hypinductionbis}. 
It remains to prove~\eqref{control2}.
Let $A \in\mathcal G_\varepsilon$ and $B\in\sigma( \cup_{j\ge p_\varepsilon} F^{-j} (\mathcal G_\varepsilon))$.
It follows from Condition~(C) of Hypothesis~\ref{hypinductionbis} that
\begin{align*}
\left|\mu_M(B \cap A) - \mu_M(A)\mu_M(B)\right|
\le C\Vert \bI_{B}\circ F^{-p_\varepsilon}\Vert_\gamma^- \Vert \bI_{A}\circ F^{-1}\Vert_\gamma^+ \theta^{p_\varepsilon-1}\le C \theta^{p_\varepsilon-1}\, ,
\end{align*}
where we used the fact that $[\bI_A\circ F^{-1}]_\gamma^+=0$ (since it is a measurable function of $(R,\cI)\circ F^{-1}$) and $[\bI_B\circ F^{-p_\varepsilon}]_\gamma^-=0$ 
(since $\bI_B\circ F^{-p_\varepsilon}$ is a measurable function of $((R,\cI)\circ F^k,\ k\ge 0)$).\\
This ends the proof of the convergence in distribution of
the family of point processes. 
Finally, \eqref{cvgcedistrib0} then follows from \cite[Theorem 1.3]{TyranKaminska2} or from Proposition~\ref{doublecvgce2} with $\gamma=1$ and $f=\bI_M$ (the tightness in $M_1$ being automatic since we study increasing processes).
\end{proof}

\subsubsection{Proof of Theorem~\ref{thm:LLN}}
Theorem~\ref{thm:LLN} follows from the next corollary.
\begin{corollary}\label{COROLLN}
Assume Hypotheses~\ref{hypinduction},~\ref{hypDelta}  and~\ref{hypinductionbis} and that $(\cM,\cF,\mu)$ is recurrent ergodic. Then, for any $f\in L^1(\cM,\mu)$,
 $$\left(\frac{S_{\lfloor nt\rfloor} f}{n^\alpha}\right)_{t\in[0;T]}\Longrightarrow  \frac{\int_{\cM}f\, d\mu}{\mu(M)} \, (\Delta^{-1}(t))_{t\in[0;T]} \, ,$$
 with respect to $\mu$ and for the uniform metric.
\end{corollary}
\begin{proof}
The result for $f=\bI_M$ follows from Proposition~\ref{propcvPP}. The result for any  $f\in L^1(\cM,\mu)$
then follows from the Hopf ratio ergodic theorem.
\end{proof}
Let us say a few words about the above limit.
\begin{remark}\label{remchgtinvDelta}
It follows from the definition of the generalized inverse that
$$
\Delta^{-1}(t)/\mu(M)=(\Delta(\mu(M)\cdot))^{-1}(t)\, .
$$
But, since $\Delta$ is $\alpha$-stable, the process $\Delta(\mu(M)\cdot)$ has the same distribution as
$\widetilde\Delta:=\mu(M)^{\frac 1\alpha}\Delta$ which is a $\alpha$-stable L\'evy process such that $\widetilde\Delta(1)$ has characteristic function 
$\Psi_{\alpha,c_0(\mu(M))^{\frac 1\alpha}}$.
\end{remark}

\begin{proof}[Proof of Theorem~\ref{thm:LLN}]
Recall that the ergodicity of $(\cM,\cF,\mu)$ has been proved in Theorem~\ref{LLNfails}. 
Theorem~\ref{thm:LLN} follows from Corollary~\ref{COROLLN} due to Remarks~\ref{rk:hypsatisfied} and~\ref{remchgtinvDelta}. 
It follows from the previous remark that $\Delta^{-1}/\mu(M)$ has the same distribution
as $\widetilde \Delta^{-1}=\Delta_1^{-1}$ since $c_0\mu(M)^{1/\alpha}=c_1$ (see Remark~\ref{rk:hypsatisfied} and~\eqref{:c1} where $c_1$ is defined). This ends the proof of Theorem~\ref{thm:LLN}.
\end{proof}

\subsection{FCLT for maps and proof of
Theorem~\ref{thm:TCL}}\label{proofTCLmap}
\subsubsection{General FCLT for maps}
Again, the proof of Theorem~\ref{thm:TCL} is postponed at the end of this subsection and will follow from general results stated under our general hypotheses.

Consider a function $f:\cM\rightarrow\mathbb R$ centered 
with respect to $\mu$.
We write $\widetilde f:M\rightarrow \mathbb R$ given by $$\widetilde f(x)=\sum_{k=0}^{R(x)-1}f\circ \mathcal F^k(x).$$ Note that, if $f$ is null outside $M$, we simply have 
$\widetilde f(x)=f(x)$.
Due to Propositions~\ref{prop:CLT} and~\ref{doublecvgce2}, the
Central Limit Theorem for $\left(\left(S_{\lfloor nt\rfloor}f\right)_t\right)_n$ will follow
from the following result. 
\begin{proposition}\label{doublecvgce}
Assume Hypotheses~\ref{hypinduction},~\ref{hypDelta}  and~\ref{hypinductionbis} and that $\Vert \widetilde f\circ F^{-1}\Vert_\gamma^++\Vert \widetilde f\Vert_\gamma^-<\infty$.
Then the family of couples of processes
$$\left(\left(  \left(\varepsilon^{\frac 12}
\int_0^{t/\varepsilon}\widetilde f\circ F ^
{\lfloor s\rfloor}\, ds
\right)_{t\in[0,T]}, \mathcal N^{\mathcal I}_\varepsilon \right)\right)_\varepsilon$$
converges in distribution with respect to $\mu_M$ (for respectively the uniform topology and the vague convergence topology), as $\varepsilon\rightarrow 0$, to $(B,\mathcal P)$
where $B$ and $\cP$ are two processes independent one from  the other, $B$ being a Brownian motion with variance $\sigma_{\widetilde f}^2:=\sum_{k\in\mathbb Z}\mathbb E_{\mu_M}\left[\widetilde f.\widetilde f\circ F^k\right]$, and
$\mathcal P$ being a Poisson process on $[0,+\infty)\times\mathbb R^*\times\{1,...,I_0\}$ with intensity of density $(t,y,i)\mapsto c_i^\alpha\alpha |y|^{-\alpha-1}\bI_{\{y>0\}}$.\\
In particular, Proposition~\ref{doublecvgce2} applies and the FCLT stated in 
Proposition~\ref{prop:CLT} holds true for $(\mathcal M,f,\mu)$ and $f:=\widetilde f\bI_M$
with $c_0^\alpha:=\sum_{i=1}^{I_0}c_i^\alpha$.
\end{proposition}
Proposition~\ref{doublecvgce} states a result of joint
convergence in distribution for a couple made of two families of processes. Of course, since the two families of processes are not independent, this result 
is not a direct consequence of the two separate convergence results, but its proof will use these separate convergence results. 
Let us recall that the convergence in distribution of $(\mathcal N_\varepsilon)_\varepsilon$ has been proved in Proposition~\ref{propcvPP}. 
Before proving Proposition~\ref{doublecvgce}, let us prove the convergence of the first family of processes.
\begin{proposition}\label{standardCLT}
Assume that $(M,F,\mu_M)$
is invertible and satisfies Condition (C) of Hypothesis~\ref{hypinductionbis} and that $\Vert \widetilde f\circ F^{-1}\Vert_\gamma^++\Vert \widetilde f\Vert_\gamma^-<\infty$. 
Then the family of processes
$$\left( \left(\varepsilon^{\frac 12}
\int_0^{t/\varepsilon}\widetilde f\circ F ^{\lfloor s\rfloor}\, ds
\right)_{t\in[0,T]}\right)_\varepsilon$$
converges in distribution (with respect to $\mu_M$) for the uniform metric, as $\varepsilon\rightarrow 0$, to a Brownian motion with variance $\sigma^2_{\widetilde f}=\sum_{k\in\mathbb Z}\mathbb E_{\mu_M}[\widetilde f.\widetilde f\circ F^k]$.
\end{proposition}
\begin{proof}
This is a classical result. 
For example, we can apply \cite[Theorem 2.1.3]{FP_ESAIM_02} with $F(x,\omega):=\widetilde f(\omega)$. Indeed the multiple decorrelation assumption of \cite{FP_ESAIM_02} follows from Condition~(C) of Hypothesis~\ref{hypinductionbis} and from
the comment thereafter. 
For any natural integers $n,m,m'$, $k_1,...,k_m$ and $\ell_1,...,\ell_{m'}$, 
\begin{align*}
&\left|\int_M\prod_{i=1}^m \widetilde f\circ F^{k_i+n}.\prod_{j=1}^{m'}\widetilde f\circ F^{\ell_j}\, d\mu_M-\int_M\prod_{i=1}^m \widetilde f\circ F^{k_i}\, d\mu_M\int_M\prod_{j=1}^{m'}\widetilde f\circ F^{\ell_j}\, d\mu_M\right|\\
&\quad\quad\leq C 
\left\|\prod_{i=1}^m \widetilde f\circ F^{k_i}\right\|_{\gamma}^-\left\|\prod_{j=1}^{m'}\widetilde f\circ F^{\ell_j-\max_{j'}\ell_{j'}-1}\right\|_{\gamma}^+ \theta^{n-\max_{j'}\ell_{j'}-1}\\
&\quad\quad\leq C 
 \prod_{i=1}^m\left\| \widetilde f\right\|_{\gamma}^-.\prod_{j=1}^{m'}\left\|\widetilde f\circ F^{-1}\right\|_{\gamma}^+ \theta^{n-\max_{j'}\ell_{j'}-1}\, .
\end{align*}

\end{proof}
\begin{proof}[Proof of Proposition~\ref{doublecvgce}]
The tightness follows from the
separate convergences in distribution of
$\left(\mathcal N_\varepsilon\right)_\varepsilon$
(stated in Proposition~\ref{propcvPP})
and of $\left(\left(\sqrt{\varepsilon}
\int_0^{t/\varepsilon}\widetilde f\circ 
F^{\lfloor s\rfloor}\, ds\right)_{t\in[0;T]}\right)_\varepsilon$ (stated in Proposition~\ref{standardCLT}). 
To identify the limit, let us prove that,
for any  $D=\bigcup_{j=1}^J[a_j,b_j]\times (c_j,d_j)\times\{i_j\}$, $0=x_0\le x_1< ...<x_J$ 
and, for any $t_1,...,t_J\in\mathbb R$ and any positive integer $m$,
\begin{equation}\label{(I)} 
\boxed{\mathbb E_{\mu_M}\left[(\mathcal N^\cI_\varepsilon(D))^m\prod_{j=1}^J e^{it_j\sqrt{\varepsilon}
\int_{x_{j-1}/\varepsilon}^{x_j/\varepsilon} \widetilde f\circ F
^{\lfloor s\rfloor}\, ds}\right] \sim \mathbb E\left[(\mathcal P(D))^m\right]\prod_{j=1}^J\mathbb E\left[ e^{-\frac {\sigma^2t_j^2(x_j-x_{j-1})^2}{2}}\right] \, , }\end{equation}
as $\varepsilon\rightarrow 0$. 
It follows from the 
definition of $\mathcal N_\varepsilon$ given in~\eqref{eq:DefNeps} and from the above definition of $D$
that
$$\mathcal N_\varepsilon(D)=\sum_{j=1}^J\sum_{n=
\lceil a_j/\varepsilon\rceil}^{
\lfloor b_j/\varepsilon\rfloor}\bI_{\{(R,\cI)\circ F^n\in(c_j\varepsilon^{-\frac 1\alpha};d_j\varepsilon^{-\frac 1\alpha})\times\{i_j\}\}},$$ we are led to the study of sums of
\begin{equation*}
I_{\varepsilon,n_1,...,n_m}^{(j_1,...,j_m)}:=\mathbb E_{\mu_M}\left[
\prod_{k=1}^m\bI_{\{F^{n_k}\in M^{[i_{j_k}]},R\circ F^{n_{k}}\in(c_{j_k}\varepsilon^{-\frac 1\alpha};d_{j_k}\varepsilon^{-\frac 1\alpha})\}}\prod_{j=1}^J e^{it_j\sqrt{\varepsilon}
\int_{x_{j-1}/\varepsilon}^{x_j/\varepsilon} \widetilde f\circ F
^{\lfloor s\rfloor}\, ds}\right]
\, ,
\end{equation*}
over $j_1,...,j_m\in\{1,...,J\}$ and $a_{j_k}/\varepsilon\le n_k\le b_{j_k}/\varepsilon$.
We rewrite
the sum $\sum_{j=1}^J\int_{x_{j-1}/\varepsilon}^{x_j/\varepsilon} t_j\widetilde f\circ F^{\lfloor s\rfloor}\, ds$
in a sum of $J':=J+m$ terms of the form
$$t'_j\int_{x'_{j-1}/\varepsilon}^{x'_j/\varepsilon} \widetilde f\circ F^{\lfloor s\rfloor}\, ds\, , $$
where $x'_0=x_0$ and where $x'_j/\varepsilon$ is the $j$-th smallest
element of $\{x_1/\varepsilon,...,x_J/\varepsilon\}\cup\{n_1,...,n_m\}$.
Set $K_\varepsilon:=(\log\varepsilon)^2$. Introducing gaps, we will approximate $I_{\varepsilon,n_1,...,n_m}$ by
\begin{equation*}
I^{*,(j_1,...,j_m)}_{\varepsilon,n_1,...,n_m}:=\mathbb E_{\mu_M}\left[
\prod_{k=1}^m\bI_{\{F^{n_k}\in M^{[i_{j_k}]},R\circ F^{n_k}\in(c_{j_k}\varepsilon^{-\frac 1\alpha};d_{j_k}\varepsilon^{-\frac 1\alpha})\}}\prod_{j=1}^{J'} e^{it'_j\sqrt{\varepsilon}\int_{x^-_{j}/\varepsilon}^{x^+_j/\varepsilon} \widetilde f\circ F^{\lfloor s\rfloor}\, ds}\right]
\, ,
\end{equation*}
with $x^-_{j}/\varepsilon=x'_{j-1}/\varepsilon+K_\varepsilon$
and $x^+_{j}/\varepsilon=x'_{j}/\varepsilon-K_\varepsilon$
and with convention $\int_a^b\, ...ds=0$ if $b<a$. 
Observe that
\begin{align*}
&\left| I_{\varepsilon,n_1,...,n_m}^{(j_1,...,j_m)}-I^{*,(j_1,...,j_m)}_{\varepsilon,n_1,...,n_m}  \right|\le J'\Vert \widetilde f\Vert_\infty\max_j|t_j|\sqrt{\varepsilon}2K_\varepsilon\\
&\quad\quad\quad\quad\quad \quad\quad \mu_M\left(
\forall k=1,..., m,\ 
R\circ F^{n_k}\in(c_{j_k}\varepsilon^{-\frac 1\alpha};d_{j_k}\varepsilon^{-\frac 1\alpha})
\right)\, .
\end{align*}
Observe that, for $0<a<b$,
\beq
\left[ e^{it\int_{a}^{b}\widetilde f\circ F^{\lfloor s\rfloor}\, ds}\right]_\gamma^- \le |t|\int_{a}^{b}\left[\widetilde f\right]_\gamma^-\gamma^{\lfloor s\rfloor}\, ds\le |t|\frac{\gamma^{\lfloor a\rfloor}}{1-\gamma}\left[\widetilde f\right]_\gamma^-\, 
\eeq
and
$$
\left[ e^{it\int_{-b}^{-a}\widetilde f\circ F^{\lfloor s\rfloor}\, ds}\right]_\gamma^+ 
\le |t|\int_{a}^{b}\left[\widetilde f\right]_\gamma^+\gamma^{\lfloor s\rfloor}\, ds\le |t|\frac{\gamma^{\lfloor a\rfloor}}{1-\gamma}\left[\widetilde f\right]_\gamma^+\, .$$
Furthermore $[\bI_{B} ]_\gamma^-=0$ for all $B\in\sigma((R,\cI)\circ F^k,\, k\ge 0)$ and  $[\bI_{D} ]_\gamma^+=0$ for all $D\in\sigma((R,\mathcal I)\circ F^{-k},\, k\ge 1)$. As usual for H\"older norms, 
$\Vert gh\Vert_\gamma^\pm\le \Vert g\Vert_\gamma^\pm\Vert h\Vert_\gamma^\pm$. 
Thus, applying several times Condition (C) of Hypothesis~\ref{hypinductionbis}, 
setting $x_{-1}^+=-\infty$ and
$x_{J'+1}^-=+\infty$, we obtain that
\begin{align*}
&I^{*,(j_1,...,j_m)}_{\varepsilon,n_1,...,n_m}=
\prod_{j=1}^{J'}\mu_M\left(\forall k\mbox{ s.t. } x_{j-1}^+/\varepsilon<n_k<x_j^-/\varepsilon,\ F^{n_k}\in M^{[i_{j_k}]},\ R\circ F^{n_k}\in]c_{j_k}\varepsilon^{-\frac 1\alpha};d_{j_k}\varepsilon^{-\frac 1\alpha}[\right)\\
&\quad\quad\quad 
\mathbb E_{\mu_M}\left(
e^{it_j\sqrt{\varepsilon}\int_{x^-_{j}/\varepsilon}^{x^+_j/\varepsilon} \widetilde f\circ F^{\lfloor s\rfloor}\, ds}\right)
+\mathcal O\left(\theta^{K_\varepsilon}\right)
\\
&=
\mu_M\left(\forall k=1,...,m,\ F^{n_k}\in M^{[i_{j_k}]},\ R\circ F^{n_k}\in]c_{j_k}\varepsilon^{-\frac 1\alpha};d_{j_k}\varepsilon^{-\frac 1\alpha}[\right)\mathbb E_{\mu_M}\left(
\prod_{j=1}^{J'}e^{it_j\sqrt{\varepsilon}\int_{x^-_{j}/\varepsilon}^{x^+_j/\varepsilon} \widetilde f\circ F^{\lfloor s\rfloor}\, ds}\right) \\
&\quad\quad\quad +\cO\left(\theta^{K_\varepsilon}\right)\\
&=
\mu_M\left(\forall k=1,...,m,\ F^{n_k}\in M^{[i_{j_k}]},\ R\circ F^{n_k}\in(c_{j_k}\varepsilon^{-\frac 1\alpha};d_{j_k}\varepsilon^{-\frac 1\alpha})\right)\\
&\quad
\left(\mathbb E_{\mu_M}\left(
\prod_{j=1}^Je^{it_j\sqrt{\varepsilon}\int_{x_{j-1}/\varepsilon}^{x_j/\varepsilon} \widetilde f\circ F^{\lfloor s\rfloor}\, ds}\right)
+\mathcal O(\sqrt{\varepsilon}K_\varepsilon)\right)+\mathcal O\left(\theta^{K_\varepsilon}\right)\, .
\end{align*}
We infer by summation over  $j_1,...,j_m\in\{1,...,J\}$ and over $n_1,...,n_m$ such that $a_{j_k}/\varepsilon\le n_k\le b_{j_k}/\varepsilon$ that
\begin{align*}
\mathbb E_{\mu_M}&\left((\mathcal N_\varepsilon^\cI(D))^m\prod_{j=1}^J e^{it_j\sqrt{\varepsilon}\int_{x_{j-1}/\varepsilon}^{x_j/\varepsilon} \widetilde f\circ F^{\lfloor s\rfloor}\, ds}\right)\\
&=\mathbb E_{\mu_M}\left((\mathcal N^\cI_\varepsilon(D))^m\right)
\left(\mathbb E_{\mu_M}\left(
\prod_{j=1}^Je^{it_j\sqrt{\varepsilon}\sum_{\ell=x_{j-1}/\varepsilon}^{x_j/\varepsilon} \widetilde f\circ F^{\ell}}\right)
+\mathcal O(\sqrt{\varepsilon}K_\varepsilon)\right)+\mathcal O\left(\varepsilon^{-m}\theta^{K_\varepsilon}\right)\, .
\end{align*}
To conclude the proof of \eqref{(I)}, it remains:
\begin{itemize}
\item to notice that  $\left(\left(\sqrt{\varepsilon}\sum_{k=0}^{\lfloor t/\varepsilon\rfloor}\widetilde f\circ F^k\right)_t\right)_\varepsilon$ converges to a Brownian motion of variance $\sigma^2_{\widetilde f}$, in the sense of finite dimensional distributions (this follows from Proposition~\ref{standardCLT}),
\item to prove that $\mathbb E_{\mu_M}\left((\mathcal N^\cI_\varepsilon(D))^m\right)\sim\mathbb E[(\mathcal P(D))^m]$. 
\end{itemize}
This last estimate will use Proposition~\ref{propcvPP} which ensures that
$(\mathcal N^\cI_\varepsilon)_\varepsilon$ converges to $\mathcal P$ as $\varepsilon\rightarrow 0$. To conclude the convergence of the moments via uniform integrability, it is thus enough to prove that, for any positive integer $m$, $\limsup_{\varepsilon\rightarrow 0}\Vert \mathcal N^\cI_\varepsilon(D)\Vert_{L^m(\mu_M)}<\infty$.
Observe that
$$
\Vert \mathcal N^\cI_\varepsilon(D)\Vert_{L^m(\mu_M)}^m\le m !
\sum_{p=1}^m
\sum_{0\le k_1<...<k_p\le \max_i b_i/\varepsilon}\mu_M\left( \forall j=1,...,p,\quad 
    R\circ F^{k_j}>\min_ic_i \varepsilon^{-\frac 1\alpha}    \right)
\, .$$
But, due to Condition (C) of Hypothesis~\ref{hypinductionbis}, 
if 
there exists $j'=1,...,p-1$ such that $k_{j'+1}>k_{j'}+K_\varepsilon$, then 
\begin{align*}
\mu_M&\left( \forall j,\  
    R\circ F^{k_j}>\min_ic_i \varepsilon^{-\frac 1\alpha}    \right)\\
    &=\mu_M\left(\forall j\le j', \ R\circ F^{k_j}>\min_ic_i \varepsilon^{-\frac 1\alpha}\right)\mu_M\left(\forall j> j', \ R\circ F^{k_j}>\min_ic_i \varepsilon^{-\frac 1\alpha}\right)+\mathcal O\left(\theta^{K_\varepsilon}\right)\, .
\end{align*}
Let $L$ be the number of such $j'$'s and, for all $l=1,...,L$, we set
$j_l$ for the $l$-th smallest $j'=1,...,p-1$ such that $k_{j'+1}>k_{j'}+K_\varepsilon$. We also set $j_0:=0$, $j_{L+1}:=p$. Applying successively the above estimate with $j_1,...,j_L$, we obtain
\begin{align}
\mu_M&\left( \forall j,\  
    R\circ F^{k_j}>\min_ic_i \varepsilon^{-\frac 1\alpha}    \right)\\
    &=\prod_{l=1}^{L+1} \mu_M\left( \forall j=j_{l-1}+1,...,j_l,\  
    R\circ F^{k_j}>\min_ic_i \varepsilon^{-\frac 1\alpha}    \right)
+\mathcal O(\theta^{K_\varepsilon})\, .
\end{align}
It follows from Conditions (B) and (C) of Hypothesis~\ref{hypinductionbis} that
$$ 
 \mu_M\left( \forall j=j_{l-1}+1,...,j_l,\  
    R\circ F^{k_j}>\min_ic_i \varepsilon^{-\frac 1\alpha}    \right)\le C(\varepsilon/\min_ic_i^\alpha)^{1+\frac\eta\alpha\bI_{j_l>1+j_{l-1}}}\, .
$$
But, for fixed $L$ and fixed pairwise distinct indices $j_1,...,j_L\in\{1,...,p-1\}$, the number of $(k_1,...,k_p)$ such that $0\le k_1<...<k_p\le\max_i b_i/\varepsilon$ and  $k_{j+1}-k_{j}\le K_\varepsilon$ for all $j\in\{1,...,p-1\}\setminus\{j_1,...,j_L\}$ is bounded by $(\max_i b_i/ \varepsilon)^{L+1} K_\varepsilon^{p-1-L}$.
We conclude that
$$
\Vert \mathcal N^\cI_\varepsilon(D)\Vert_m^m\le \mathcal O\left(
\varepsilon^{-p}\theta^{K_\varepsilon}+\sum_{L=0}^{p-1} (\varepsilon \max_ib_i/\varepsilon)^{L+1}  (K_\varepsilon)^{p-1-L} \varepsilon^{\frac\eta\alpha(\bI_{\{L<p-1\}})}
\right)=\mathcal O(1)\, .
$$
As explained above, this ends the proof of the convergence of \eqref{(I)}. To conclude, we use the tightness and consider a limit point $(\mathfrak A,\mathfrak B)$ of the couple of processes for the convergence in distribution. Then
there exists $\varepsilon_k\rightarrow 0$ such that
$$
(\mathcal N^\cI_{\varepsilon_k}(D))^m\prod_{j=1}^J e^{it_j\sqrt{\varepsilon_k}\int_{x_{j-1}/\varepsilon_k}^{x_j/\varepsilon_k} \widetilde f\circ F^{\lfloor s\rfloor}\, ds}$$
converges in distribution, as $k\rightarrow +\infty$, to $(\mathfrak B(D))^m\prod_{j=1}^J e^{it_j(\mathfrak A_{x_j}-\mathfrak A_{x_{j-1}})}$. Thus, the uniform integrability of $((\mathcal N^\cI_\varepsilon(D))^m)_{\varepsilon}$ ensures that
$$
\lim_{k\rightarrow +\infty}
\mathbb E_{\mu_M}\left[ (\mathcal N^\cI_{\varepsilon_k}(D))^m\prod_{j=1}^J e^{it_j\sqrt{\varepsilon_k}\int_{x_{j-1}/\varepsilon_k}^{x_j/\varepsilon_k} \widetilde f\circ F^{\lfloor s\rfloor}\, ds}\right]=\mathbb E\left[(\mathfrak B(D))^m\prod_{j=1}^J e^{it_j(\mathfrak A_{x_j}-\mathfrak A_{x_{j-1}})}\right]\, .$$
But \eqref{(I)} implies that
$$
\mathbb E\left[(\mathfrak B(D))^m\prod_{j=1}^J e^{it_j(\mathfrak A_{x_j}-\mathfrak A_{x_{j-1}})}\right]=\mathbb E[(\mathcal P(D))^m]
\mathbb E\left[\prod_{j=1}^J e^{it_j(B_{x_j}-B_{x_{j-1}})}\right]
$$
In particular
$$
\mathbb E\left[e^{it\mathfrak B(D)}\prod_{j=1}^J e^{it_j(\mathfrak A_{x_j}-\mathfrak A_{x_{j-1}})}\right]=\mathbb E[e^{it\mathcal P(D))}]
\mathbb E\left[\prod_{j=1}^J e^{it_j(B_{x_j}-B_{x_{j-1}})}\right]
$$
from which we conclude that the only limit in distribution is $(B,\mathcal P)$.
The consequence about the convergence of the FCLT for $f=\widetilde f\bI_M$
follows from Proposition~\ref{doublecvgce2}.
\end{proof}
\subsubsection{proof of
Theorem~\ref{thm:TCL}}

Here we provide the proof of Theorem~\ref{thm:TCL}.
Up to increase the value of $\widehat t$ used in Section~\ref{Sec:3} to define the set $M$, we can assume without any loss of generality that $f$ is supported in $M$. Furthermore, Hypotheses~\ref{hypinduction} and~\ref{hypinductionbis} are satisfied (see Remark~\ref{rk:hypsatisfied} and~\eqref{:c1}). It then follows from Propositions~\ref{doublecvgce} and~\ref{doublecvgce2} that~\eqref{cvgcedistrib2} holds and so we end the proof by applying Proposition~\ref{prop:CLT}.
The limit is then given by $(B_{\Delta^{-1}(t)})_t$
where $B$ is a Brownian motion with variance
$\sigma^2=\sigma^2_f/\mu(M)$, where
$\sigma^2_f:=\sum_{k\in\mathbb Z}\int_{M}
f.f\circ F^k\, d\mu$ and where $\Delta$ is an independent $\alpha$-stable L\'evy process with characteristic distribution at time $1$ given by $\Psi_{\alpha,c_0}$
(see Remark~\ref{rk:hypsatisfied} and~\eqref{:c1}) with 
$$c_0=(4\mathfrak C/(\beta \mu(M)))^{\frac 1\alpha}I_{1/\alpha}=c_1/\mu(M)^{\frac 1\alpha}.$$
To conclude that the limit has the announced form, we observe that $B_{\Delta^{-1}(t)}=B_{ \mu(M)\Delta^{-1}(t)/\mu(M)}$, and $B_{\mu(M)\cdot}$
is a Brownian motion of variance $\sigma^2_f$
whereas it has been seen in Remark~\ref{remchgtinvDelta} than $\Delta^{-1}(\cdot)/\mu(M)$ has same distribution as 
$\Delta_1^{-1}$ with characteristic function at 1 given by $\Psi_{\alpha,c_1}$.

\section{FCLT for suspension flows}\label{proofTCLflow}
We will extend the previous general results to flows using the suspension flows structure recalled in the next hypothesis.
\begin{hypothesis}[Suspension flow]\label{hypspecial}
Let $(\mathcal M_0,(\Phi^t)_t,\nu)$ be a measure preserving 
dynamical system in continuous time that can be represented by the suspension flow over $(M,F,\mu_{|M})$ a finite measure preserving dynamical system with $M\subset \cM_0$ and with roof function $\widetilde\tau$ corresponding to the first return time of the flow $\Phi$ in $M$, which is given by
\beq
\widetilde\tau(x):=\inf\{t>0\, :\, \Phi^t(x)\in M\}\, .
\eeq
This means that $\nu$ is the image measure by $\Pi:(x,s)\mapsto \Phi^s(x)$ of the restriction of the measure $\mu_{|M}\otimes {\mbox Leb}$ on $\cM_1':=\{(x,s)\in M\times[0,+\infty)\, :\, s<\widetilde\tau(x)\}$, where $\mu_{|M}=\mu(\cdot\cap M)$.\\
Let $\mu_M:=\mu_{|M}/\mu(M)$ be the probability measure proportional to $\mu_{|M}$.
\end{hypothesis}

\begin{remark}\label{remhyp}
Hypothesis~\ref{hypspecial} holds true with
$\widetilde\tau(x)=\sum_{k=0}^{R(x)-1}\tau\circ T^k(x)$ if Hypothesis~\ref{hypinduction} holds true and if 
$\cM\subset \cM_0$ and if
$(\mathcal M_0,(\Phi^t)_t,\nu)$ can be represented as the suspension flow over $(\cM,\cF,\mu)$ with roof function $\tau:\cM\rightarrow ]0,+\infty[$ the first return time of $\Phi$ to $\mathcal M$.
\end{remark}
We will obtain different results depending on whether the measure $\nu$
is finite or infinite, which is equivalent to the integrability or non integrability of $\widetilde \tau$ with respect to $\mu_M$ (or equivalently to the integrability or non-integrability of $\tau$ with respect to $\mu$, in the context of Remark~\ref{remhyp}).

\subsection{FCLT for flows in finite measure and proof of Theorem~\ref{thm:TCLflow}}\label{proofstandardTCLflow}
\subsubsection{General FCLT for flows in finite measure}
The proof of Theorem~\ref{thm:TCLflow} will follow from the following theorem valid in a general context.
\begin{theorem}[finite measure case]\label{thmflowfinitemeasure}
Assume Hypothesis~\ref{hypspecial} with $\nu$ finite, that $(M,F,\mu_M)$ is invertible and ergodic and satisfies Condition~(C) of Hypothesis~\ref{hypinductionbis}. 
Let $g:\mathcal M_0\rightarrow\mathbb R$ be a $\nu$-centered bounded function such that the function $U_ g:M\rightarrow\mathbb R$ given by $U_g(x)= \int_0^{\widetilde\tau(x)}g(\Phi^s(x))\, ds$ satisfies $\Vert U_ g\circ F^{-1}\Vert_\gamma^++\Vert U_ g\Vert_\gamma^-<\infty$ and is such that $ U_{|g|}$ is uniformly bounded. 
Then
\begin{equation}\label{limiteflow}
\forall T>0, \quad \left(\frac{\int_0^{nt} g\circ \Phi^s\, ds}{n^{\frac 1{2}}}\right)_{t\in[0;T]}
\Longrightarrow(\sigma
B_{t\mu(M)/\nu(\mathcal M_0)
})_{t\in[0;T]}\, ,
\end{equation}
in the uniform topology, 
with respect to $\nu$, as $n\rightarrow +\infty$, 
where $B$ is a standard Brownian motion
and where $\sigma^2:=\sum_{k\in\mathbb Z}\mathbb E_{\mu_M}[U_ g.U_ g\circ F^k]$.
\end{theorem}
\begin{remark}\label{remsigmagflow}
Observe that the asymptotic variance of the limit Brownian motion in the right hand side of~\eqref{limiteflow} is given by
\beq\label{sigmag_flow}
\boldsymbol{\sigma}^2_g:=\sigma^2\frac{\mu(M)}{\nu(\mathcal M_0)}= \frac 1{\nu(\mathcal M_0)}\sum_{k\in\mathbb Z}\int_M U_g. U_g\circ F^k\, d\mu\, .
\eeq
\end{remark}

\begin{remark}[Non-degeneracy of $\boldsymbol\sigma_g$]
In both flow theorems the asymptotic variance equals, up to a positive multiplicative constant, the Green--Kubo variance $\sigma_{U_g}^2$ of the induced observable $U_g$ for $(M,F,\mu_M)$ obtained by integrating $g$ along the flow over a return excursion (cf.~\eqref{sigmag_flow} for the finite-area case and the discussion after Theorem~\ref{thm:TCLflow2} for the infinite-area one). Hence, as in the previous remark for the map, $\boldsymbol\sigma_g\ge 0$ and
$$
\boldsymbol\sigma_g=0\quad\Longleftrightarrow\quad \text{$U_g$ is an $L^2(\mu_M)$ coboundary for $F$.}
$$
Equivalently, $\boldsymbol\sigma_g=0$ iff $g$ is a measurable coboundary for the flow, i.e.\ the additive functional $t\mapsto\int_0^{t}g\circ\Phi^s\, ds$ stays bounded in $L^2$ (heuristically, $g=Xh$ is a derivative of some $h$ along the flow direction $X$). Thus the Brownian limits in Theorems~\ref{thm:TCLflow} and~\ref{thm:TCLflow2} are non-degenerate for all $g$ outside this non-generic class.
\end{remark}

\begin{proof}[Proof of Theorem~\ref{thmflowfinitemeasure}]
Due to \cite[Theorem 1]{Zweimuller}, it is enough to prove the
convergence in distribution with respect to the probability measure $P$ on $\cM_1'$
given by $dP(x,s)=d\mu_M(x)ds/\widetilde\tau(x)$. 
Since $\nu$ is finite, $\widetilde\tau$ is $\mu_M$-integrable and 
\beq
\int_{M}\widetilde\tau\, d\mu_M=
\frac 1{\mu(M)}\int_{M}\widetilde\tau\, d\mu= \frac{\nu(\mathcal M_0)}{\mu(M)}\, .
\eeq
It follows from Proposition~\ref{standardCLT} that
\begin{equation}
\left(n^{-\frac 12}\sum_{k=0}^{\lfloor nt\rfloor-1} U_ g\circ F^k\right)_t\stackrel{distrib}{\longrightarrow}\left(\sigma W_t\right)_{t}\, ,
\end{equation}
with respect to $\mu_M$ which is the image measure of $P$ by $(x,s)\mapsto x$, and
where $W$ is a standard Brownian motion and  
$\sigma^2:=\sum_{k\in\mathbb Z}\mathbb E_{\mu_M}[ U_g. U_g\circ F^k]$. 
Using the ergodicity of $(M,F,\mu_M)$ and the standard time change theorem \cite[section 14]{Billingsley}, we infer that
\begin{equation*}
\left(\frac{\int_0^{nt} g\circ \Phi^s\, ds}{n^{\frac 1{2}}}\right)_{t\in[0;T]}
=n^{-\frac 12}\left(\sum_{k=0}^{{\mathcal T}^{-1}(nt)}
 U_g\circ F^k+\mathcal O(1)\right)_{t\in[0;T]}\stackrel{distrib}{\longrightarrow}(\sigma B_{t/\int_{M}\widetilde\tau\, d\mu_M})_{t\in[0;T]}
\end{equation*}
with respect to $P$, 
with $\mathcal T(t)=\sum_{k=0}^{\lfloor t\rfloor-1}\widetilde{\tau}\circ \mathcal F^k$, since $\mathcal T(t)/t$
converges $\mu_M$-almost surely to $\int_{M}\widetilde\tau\, d\mu_M= \nu(\mathcal M_0)/\mu(M)$. As announced, this ends the proof of Theorem~\ref{thmflowfinitemeasure} due to~\cite[Theorem 1]{Zweimuller}.
\end{proof}

\subsubsection{Proof of Theorem~\ref{thm:TCLflow}}\label{:proofThm5}
We consider here the case of the billiard flow in the domain with infinite cusps when $\beta>1$, so that the billiard domain has finite area. 
Recall that the unit billiard flow $(\Phi^t)_t$
preserves the Lebesgue measure $\nu$ on $\mathcal M_0:=Q\times {\mathbb S}^1$, which is thus finite (since $\beta>1$). 
The ergodicity of this flow follows from the ergodicity of the billiard map. This implies the first part of Theorem~\ref{thm:TCLflow}. 
Let $g:Q\times\mathbb S^1\rightarrow\mathbb R$ be a compactly supported H\"older function centered with respect to $\nu$. Then $ U_g:x\mapsto \int_0^{\widetilde\tau(x)}g(\Phi^s(x))\, ds$ is in $\mathcal H_\gamma^-$ and $U_g\circ F^{-1}\in\mathcal H_\gamma^+$ for some $\gamma\in(0,1)$. 
To check this, we note that since $g$ has compact support, there exists $K>0$ such that $ U_g(x)= \int_0^{\min(K,\widetilde\tau(x))}g(\Phi^s(x))\, ds$.
Thus $U_g$ is bounded, and for $x,y$ on the same unstable (resp.\ stable) manifold of $F$ the orbit arcs $(\Phi^s x)_{0\le s\le K}$ and $(\Phi^s y)_{0\le s\le K}$ stay $\mathcal O(d_e(x,y))$-close, no singularity separating them up to the separation time; since $g$ is H\"older, $U_g$ is therefore H\"older along stable/unstable manifolds. As the Euclidean distance is in turn bounded by a fixed power of $\gamma_0^{s_\pm(\cdot,\cdot)}$ along these manifolds of the uniformly hyperbolic map $F$, this yields $U_g\in\mathcal H_\gamma^-$ and $U_g\circ F^{-1}\in\mathcal H_\gamma^+$ for some $\gamma\in(0,1)$; cf.~\cite[Ch.~7]{CM}.
We conclude the proof of Theorem~\ref{thm:TCLflow} with the use of
Theorem~\ref{thmflowfinitemeasure}.

\subsection{FCLT for flows in infinite measure and proof of Theorem~\ref{thm:TCLflow2}}\label{proofTCLnonstandardflow}  
\subsubsection{General FCLT for flows in infinite measure}
We state now a continuous time version of Proposition~\ref{prop:CLT} in the context of suspension
flow considered in Hypothesis~\ref{hypspecial}. We
set $\mathcal R_{\Phi}(u):=\sum_{k=0}^{\lfloor u\rfloor-1}\widetilde \tau\circ F^k$ and, for $f:\mathcal M_0\rightarrow\mathbb R$, set also
$\mathcal S_\Phi(u)=\int_0^{u}U_f\circ F^{\lfloor s\rfloor}\, ds$,   
with $U_ f(x):=\int_{0}^{\widetilde\tau(x)}f\circ\Phi^s(x)\, ds$
for any $x\in M$.
\begin{proposition}\label{prop:CLTflow2}
Assume Hypothesis~\ref{hypspecial}. 
Let $\beta>0$ and $\gamma>0$. 
Let $f:\cM_0 \rightarrow \mathbb R$ be a measurable function.
Assume that the
following joint convergence in distribution holds true
with respect to $\mu_M$:
\begin{equation}\label{cvgcedistrib2flow}
\left(\left(\mathfrak S_\varepsilon(t):=\varepsilon^{\frac 1\gamma}\mathcal S_\Phi(t/\varepsilon)\right)_t,\left(\mathfrak R_\varepsilon(s):=
\varepsilon^{\frac 1\beta}\mathcal R_\Phi(s/\varepsilon
)\right)_s\right)\stackrel{distrib}{\longrightarrow}\left(\left(B_t\right)_t,\left(\Delta_\Phi(s)\right)_s\right)\ ,
\quad\mbox{as }\varepsilon\rightarrow 0\, ,
\end{equation}
where $B$ and $\Delta_\Phi$ are two independent processes (respectively a continuous process and a strictly increasing process)
for respectively the uniform metric on $\mathcal C([0,T])$ and the Skorohod  $M1$-metric on $\mathcal D([0;T])$ for any $T>0$.
Assume moreover that $\left(\sup_{k=0,...,N} \frac{g_{f,\Phi}\circ F^k }{N^{\frac 1\gamma}}\right)_N$ converges in probability (with respect to $\mu_M$) to 0, where we set $g_{f,\Phi}(x):=\sup_{t\in[0,\widetilde\tau(x)]}\left|\int_{0}^{t}f( \Phi^ s(x))\, ds\right|$.
Then
\begin{equation*}
\left(\varepsilon^{\frac\beta \gamma}\int_{0}^{t/\varepsilon}f\circ\Phi^s\, ds
\right)_t\stackrel{distrib}{\longrightarrow}\left(B_{\Delta_\Phi^{-1}(t)}\right)_t\, ,\quad\mbox{as }n\rightarrow +\infty\, ,
\end{equation*}
with respect to $\mu_M$.\\
If, moreover, $(\cM_0,(\Phi^t)_t,\nu)$ is ergodic, then this convergence in distribution holds 
with respect to any probability measure absolutely continuous with respect to $\nu$ (i.e. in the strong sense with respect to $\nu$).
\end{proposition}
\begin{proof}
As in the proof of Proposition~\ref{prop:CLT} with the notation $(\mathcal
M_0,\nu)$ therein,
the joint convergence implies the convergence in distribution, for the uniform norm on $[0;T]$ for all $T>0$ of, the family of processes
$\left(\mathfrak S_\varepsilon\circ\mathfrak R_\varepsilon^{-1}\right)$ to $B_{\Delta_\Phi^{-1}}$ with respect to $\mu_M$.
But
\beq\nonumber
\mathfrak R_\varepsilon^{-1}(t):=\inf\left\{u>0\, :\, \varepsilon^{\frac 1\beta}\sum_{k=0}^{\lfloor u/\varepsilon\rfloor-1}\widetilde \tau\circ F^k\ge t\right\}
\eeq
is equal to $\varepsilon \left(1+\mathfrak N(\varepsilon^{-\frac 1\beta}t)\right)$, where
$\mathfrak N(\varepsilon^{-\frac 1\beta}t)$ is the number of visits of the flow $(\Phi^t)_t$ to $M$ in the time interval $(0,\varepsilon^{-\frac 1\beta}t)$.
Recall that $\mu_M$ is the image measure of $P$ by the canonical projection $\Pi_M:\cM'_1\rightarrow M$ given by $\Pi_M(x,s)= x$, where $P$ is the probability measure on $\mathcal M'_1$ given by $dP(x,s)=ds\, d\mu_M(x)/\widetilde \tau(x)$.  
Thus we have proved the convergence in distribution with respect to $P$ of
\beq\nonumber
\left(\mathfrak S_\varepsilon(\mathfrak R_\varepsilon^{-1}(t))\circ\Pi_M=\varepsilon^{\frac 1\gamma}
\int_0^{1+\mathfrak N(\varepsilon^{-\frac 1\beta}t)\circ\Pi_M}U_f\circ F^{\lfloor s\rfloor}\circ\Pi_M\, ds\right)_{t\in[0;T]}
\eeq
to $(B_{\Delta_\Phi^{-1}(t)})_{t\in[0;T]}$. 
To conclude, we observe that, on $\mathcal M_0$,
\beq\nonumber
\sup_{t\in[0;T]}\left|\varepsilon^{\frac 1\gamma}
\int_0^{\varepsilon^{-\frac 1\beta}t} f\circ \Phi^s\circ\Pi\, ds-\mathfrak S_\varepsilon(\mathfrak R_\varepsilon^{-1}(t))\circ\Pi_M\right|\le 2 \varepsilon^{\frac 1\gamma}
\sup_{k=0,...,\mathfrak N(\varepsilon^{-\frac 1\beta}T)}g_{f,\Phi}\circ F^{k}\circ\Pi_M\, ,
\eeq
and that, for all $K>0$,
\begin{align*}
\limsup_{\varepsilon\rightarrow 0}
P&\left(\varepsilon^{\frac 1\gamma}\sup_{k=0,...,\mathfrak N(\varepsilon^{-\frac 1\beta}T)}g_{f,\Phi}\circ F^{k}\circ\Pi_M>\eta\right)\\
&\le\limsup_{\varepsilon\rightarrow 0}
\mu_M\left(\mathfrak N(\varepsilon^{-\frac 1\beta}T)>K/\varepsilon\right)+
\limsup_{\varepsilon\rightarrow 0}
\mu_M\left(\sup_{k=0,...,K\varepsilon^{-1}}\varepsilon^{-\frac 1\gamma}g_{f,\Phi}\circ F^{k}>\eta\right)\\
&\le\limsup_{\varepsilon\rightarrow 0}
\mu_M\left(\varepsilon^{\frac 1\beta}\sum_{k=0}^{\lfloor K/\varepsilon\rfloor}\widetilde\tau\circ F^k<T\right)+0\\
&\le \mathbb P\left(\Delta_{\Phi}(K)<T\right)
\, ,
\end{align*}
which goes to 0 as $K\rightarrow +\infty$. 
The last point of the proposition follows again from \cite[Theorem 1]{Zweimuller}.
\end{proof}
\begin{remark}\label{jointintegerpart}
Assume $0<\beta<1$. 
The joint convergence~\eqref{cvgcedistrib2flow} follows from the following one
\beq\nonumber
\left(\left(\mathfrak S_\varepsilon(t):=\varepsilon^{\frac 1\gamma}\mathcal S_\Phi(t/\varepsilon)\right)_t,\left(
\varepsilon^{\frac 1\beta}\sum_{k=0}^{\lfloor s/\varepsilon\rfloor-1}\lfloor\widetilde \tau\rfloor\circ F^k\right)_s\right)\stackrel{distrib}{\longrightarrow}\left(\left(B_t\right)_t,\left(\Delta_\Phi(s)\right)_s\right)
\quad\mbox{as }\varepsilon\longrightarrow 0\, .
\eeq
Indeed 
\beq\nonumber
\sup_{s\in[0;T]}\left|\varepsilon^{\frac 1\beta}\sum_{k=0}^{\lfloor s/\varepsilon\rfloor-1}(\widetilde\tau-\lfloor\widetilde \tau\rfloor)\circ F^k\right|\le T\varepsilon^{\frac 1\beta-1}\rightarrow 0
\quad\mbox{as }\varepsilon\rightarrow 0\, .
\eeq
Furthermore Propositions~\ref{doublecvgce2} and~\ref{doublecvgce} give a strategy to prove this joint distribution convergence.
\end{remark}
\begin{proposition}\label{prop:CLTflow2b}
Assume Hypothesis~\ref{hypspecial} and that $(M,\mu_M,F)$ is ergodic. 
Let $\beta\in(0,1)$. Assume that
\beq\label{RR}
\left(\mathfrak R_\varepsilon(t):=
\varepsilon^{\frac 1\beta}\mathcal R_\Phi(t/\varepsilon
)\right)_t\stackrel{distrib}{\longrightarrow}\left(\Delta_\Phi(t)\right)_t
\eeq
for the M1-metric on $\mathcal D([0;T])$ for any $T>0$, where $\Delta_\Phi$ is a $\beta$-stable increasing process. Let $f:\cM_0 \rightarrow \mathbb R$ be a $\nu$-integrable function. 
Then
\begin{equation*}
\left(\varepsilon^{\beta}\int_{0}^{t/\varepsilon}f\circ\Phi^s\, ds
\right)_t\stackrel{distrib}{\Longrightarrow}\left(\frac{\int_{\mathcal M_0}f\, d\nu}{\mu(M)}\Delta_\Phi^{-1}(t)\right)_t\, ,\quad\mbox{as }n\rightarrow +\infty\, ,
\end{equation*}
strongly in distribution with respect to $\nu$ (with respect to the uniform topology).
\end{proposition}
\begin{proof}
The ergodicity of the system
$(M,\mu_M,F)$ ensures, via the Birkhoff Ergodic Theorem, that
$(\mathfrak S_\varepsilon(t):=\varepsilon\mathcal S_\Phi(t/\varepsilon))_t$
converges almost surely, in the uniform topology, to $(t\int_M U_f\, d\mu_M)$ as $\varepsilon\rightarrow 0$. Since the limit is deterministic, this combined with~\eqref{RR} implies the
following joint convergence in distribution
with respect to $\mu_M$:
\begin{equation}\label{cvgcedistrib2flow0}
\left(\left(\mathfrak S_\varepsilon(t)\right)_t,\left(\mathfrak R_\varepsilon(s):=
\varepsilon^{\frac 1\beta}\mathcal R_\Phi(s/\varepsilon
)\right)_s\right)\stackrel{distrib}{\longrightarrow}\left(\left(t\int_MU_f\, d\mu_M \right)_t,\left(\Delta_\Phi(s)\right)_s\right)\ ,
\quad\mbox{as }\varepsilon\rightarrow 0\, ,
\end{equation}
Furthermore, the ergodicity of $(M,\mu_M,F)$ ensures that $\left(\sup_{k=0,...,N} \frac{g_{f,\Phi}\circ F^k }{N}\right)_N$ converges almost surely (with respect to $\mu_M$) to 0.
Thus, it follows from Proposition~\ref{prop:CLTflow2} with $\gamma=1$ that 
\begin{equation*}
\left(\varepsilon^{\beta}\int_{0}^{t/\varepsilon}f\circ\Phi^s\, ds
\right)_t\stackrel{distrib}{\Longrightarrow}\left(\int_M U_f\, d\mu_M{\Delta_\Phi^{-1}(t)}\right)_t\, ,\quad\mbox{as }n\rightarrow +\infty\, ,
\end{equation*}
strongly in distribution with respect to $\nu$.
We conclude by using the fact that 
$\int_MU_f\, d\mu_M=\int_{\mathcal M_0}f\, d\nu\, /\mu(M)$.
\end{proof}
\begin{lemma}\label{asymptildetau}
For the billiard flow with cusps, there exists 
a sequence of real numbers $(\varepsilon_N)_N$ such that $\lim_{N\rightarrow +\infty}\varepsilon_N=0$,
 for an excursion in the cusp parametrized by $f_{\mathfrak c}$, 
 such that we have 
\beq
\left|\widetilde \tau
- \frac{2}\beta  \left(\mathfrak cR/I_{1/\alpha}\right)^{\frac\alpha\beta}I_{\frac 1\beta-1}\right|\le  R^{\frac\alpha\beta}\varepsilon_R\, ,
\eeq
where $\alpha=\beta/(\beta+1)$ and 
$I_{\frac 1\beta-1}:=\int_{0}^{\frac \pi 2}(\sin\varphi)^
{\frac 1\beta-1}\, d\varphi$ (as defined in~\eqref{formulaI}).
\end{lemma}
\begin{proof}
Assume $R=N$ and let us use the notations of  Proposition~\ref{prop3}. We will also use intermediate results established in its proof given in Appendix~\ref{proofprop3}. Then
\begin{align}\label{decomptildetau}
\widetilde \tau&=\mathcal O(t_1+t_{N-1})
+\sum_{n=1}^{N-2}\tau_n 
\end{align}
But it follows from Item~(4) of Proposition~\ref{prop3} that $t_1\approx t_{N-1}\approx  N^{\frac{\beta}{(2\beta+1)(\beta+1)}}=o(N^{\frac 1{\beta+1}})=o(N^{\frac\alpha\beta})$.
Let us study $\sum_{n=1}^{N_2}\tau_n$, the study of $\sum_{n=N_2+1}^{N-2}\tau_n$ having the same behaviour. It follows from~\eqref{deftaun},~\eqref{Taylortnb} and \eqref{3.7} that 
\begin{align*}
\sum_{n=1}^{N_2}\tau_n
&=\mathfrak c\sum_{n=1}^{N_2}\frac{t_{n+1}^{-\beta}+t_n^{-\beta}}{\beta\sin v_n}\\
&=\mathfrak c\sum_{n=1}^{N_2}\frac{t_{n+1}^{-\beta-1}t_n\left(2+\mathcal O(
n^{-1})\right)}{\beta\sin v_n}\frac {v_{n+1}-v_n}{2\mathfrak ct_{n+1}^{-\beta-1}(1+\mathcal O(t_{n+1}^{-2(\beta+1)}))}\, .
\end{align*}
Thus
\begin{align} \label{sumtaun0}
\sum_{n=1}^{N_2}\tau_n
&= \sum_{n=1}^{N_2}\frac{
v_{n+1}-v_n
}{\beta t_n^{-1} \sin v_n}\left(1+\mathcal O\left(n^{-\frac{2(\beta+1)}{2\beta+1}}t_1^{-2\beta-2}+n^{-1}\right)\right)
\end{align}
due to Item~(4) of Proposition~\ref{prop3}.\\
Since $v_{n+1}-v_n= 2\alpha_{n+1}\le\frac{v_n}n$ (using~\eqref{3.7} and~\eqref{nalphan}) and $v_n=\mathcal O(\sin v_n)$, and so using Item~(4) of Proposition~\ref{prop3}, we infer that
\begin{align}\nonumber
\sum_{n=1}^{N_2}
\frac{v_{n+1}-v_n}{t_n^{-1}\sin v_n}&\left(n^{-\frac{2(\beta+1)}{2\beta+1}}t_1^{-2\beta-2}+n^{-1}\right)\\
\nonumber=&\cO\left(\sum_{n=1}^{N_2}
\frac{t_1n
^{\frac 1{2\beta+1}}}n\left(n^{-\frac{2(\beta+1)}{2\beta+1}}t_1^{-2\beta-2}
+n^{-1}\right)\right)\\
=&\cO\left( t_1^{-2\beta-1}+t_1\right)=\cO\left(N^{\frac\beta{(2\beta+1)(\beta+1)}}\right)=o\left(N^{\frac 1{\beta+1}}\right)\, . \label{sumtaun1}
\end{align}
But, due to Item (7) of Proposition~\ref{prop3}, $t_n^{-\beta}\sin v_n=I_{1/\alpha}^\alpha \left(\mathfrak cN)^{-\alpha}\left(1+\cO(n^{-1}\right)\right)$ 
and so
\begin{align}
\sum_{n=1}^{N_2}\frac{v_{n+1}-v_n}{ t_{n}^{-1}\sin v_n}
&=  (\mathfrak cN/I_{1/\alpha})^{\frac\alpha\beta}
\sum_{n=1}^{N_2}\frac{v_{n+1}-v_n}{ (\sin v_n)^{1-\frac 1\beta}}(1+\mathcal O(n^{-1}))\nonumber\\
& =  (\mathfrak cN/I_{1/\alpha})^{\frac\alpha\beta} (1+\mathcal O(\varepsilon_N))\int_0^{\frac \pi 2}\frac{dv}{(\sin v)^{1-\frac 1\beta}} \, ,\label{sumtaun2}
\end{align}
with $\lim_{N\rightarrow +\infty}\varepsilon_N=0$.
Combining~\eqref{sumtaun0},~\eqref{sumtaun1} and~\eqref{sumtaun2}, we infer that
\beq\label{000}
\sum_{n=1}^{N_2}\tau_n=\frac 1\beta(\mathfrak cN/I_{1/\alpha})^{\frac\alpha\beta}I_{\frac 1\beta-1}
(1+\mathcal O(\varepsilon_N))\, .
\eeq
By symmetry, we also conclude that 
\beq\label{000b}
\sum_{n=N_2+1}^{N-2}\tau_n=\frac 1\beta(\mathfrak cN/I_{1/\alpha})^{\frac\alpha\beta}I_{\frac 1\beta-1}
(1+\mathcal O(\varepsilon'_N))\, ,
\eeq
with $\lim_{N\rightarrow +\infty}\varepsilon'_N=0$. 
We conclude by combining~\eqref{decomptildetau},~\eqref{000} and~\eqref{000b}
\end{proof}
\subsubsection{Proof of 
Theorem~\ref{thm:TCLflow2}}\label{:proofThm6}
Recall that $\mathfrak c_i$ is the constant such that the $i$-th cusp can be represented by the graph of $f_{\mathfrak c_i}$. 
It has be seen in Remark~\ref{rk:hypsatisfied} that
the billiard map satisfies Hypotheses~\ref{hypinduction} and~\ref{hypinductionbis} with $c_i=c'\mathfrak c_i^{\frac 1\beta}$
with $c':=( 4/(\beta\mu(M)))^{\frac 1\alpha} I_{1/\alpha}$. 
We know that $U_f$ is in $\mathcal H_\gamma^-$ and that $U_f\circ F^{-1}\in\mathcal H_\gamma^+$ for some
$\gamma\in(0,1]$.
Thus, due to Proposition~\ref{doublecvgce},
\beq
\left(\left(\varepsilon^{\frac 12}\int_0^{t/\varepsilon}U_f\circ F^{\lfloor s\rfloor}\, ds\right)_t,\sum_{j\ge 1}\delta_{\left(j\varepsilon,\varepsilon^{\frac 1\alpha}R\circ F^j,\cI\circ F^j\right)}\right)\stackrel{distribution}\longrightarrow (B,\mathcal P)\, ,
\eeq
where $B$ is a Brownian motion with variance
$\sigma^2_{U_g}$ and where $\mathcal P$ is a Poisson process
on $[0,+\infty)\times\mathbb R^*\times\{1,...,\mathcal J\}$ with density $(u,z,i)\mapsto (c'/\mathfrak c_{i}^{1/\beta})^\alpha\alpha z^{-\alpha-1}
$, independent one from the other.
Due to Lemma~\ref{asymptildetau}, it follows that
\beq
\left(\left(\varepsilon^{\frac 12}\int_0^{t/\varepsilon}U_f\circ F^{\lfloor s\rfloor}\, ds\right)_t,\sum_{j\ge 1}\delta_{\left(j\varepsilon,\varepsilon^{\frac 1\beta}\widetilde\tau \circ F^j\right)}\right)\stackrel{distribution}\longrightarrow (B,\mathcal P'')\, ,
\eeq
where $\mathcal P''$ is the image of $\mathcal P$ by $\sum_{j}\delta_{(t_j,x_j,i_j)}\mapsto \sum_j\delta_{(t_j,a(\mathfrak c_{i_j}x_j)^{\alpha/\beta})}$, with $a:=\frac 2\beta I_{1/\alpha}^{-\frac\alpha\beta}I_{\frac 1\beta-1}$. Thus $\mathcal P''$ is a Poisson process
on $[0,+\infty)\times\mathbb R^*$ with density $(u,z)\mapsto c'_0 \beta z^{-\beta-1}
$, where
\beq
c'_0:=(c')^\alpha \sum_{i=1}^{\mathcal J}\mathfrak c_j a^{\beta}= \sum_{i=1}^{\mathcal J}2\mathfrak c_j(2/\beta)^{1+\beta}I^\beta_{\frac 1\beta-1}/\mu(M)\, ,
\eeq
recalling that $\cJ$ is the number of cusps of order $\beta$. 
Thus, it follows from Proposition~\ref{doublecvgce2} (applied to $\lfloor \widetilde \tau\rfloor$ instead of $R$) and Remark~\ref{jointintegerpart} that
Condition~\eqref{cvgcedistrib2flow} of Proposition~\ref{prop:CLTflow2} holds true with $\Delta_\Phi(1)$ having characteristic function $\Psi_{\beta,(2/\beta)^{\frac 1\alpha}(2\sum_{j=1}^{\mathcal J}\mathfrak c_j/\mu(M))^{\frac 1\beta} I_{\frac 1\beta-1}}$ and with $B$ of variance
$\mathbb E_{\mu_M}[U_f.U_f\circ F^k]=\sigma^2_{U_f}/\mu(M)$ with
\beq\label{formulasigmaginfini}
\sigma^2_{U_f}:=\sum_{k\in\mathbb Z}\int_MU_f.U_f\circ F^k\, d\mu\, .
\eeq
Thus
\begin{equation*}
\left(\varepsilon^{\frac\beta 2}\int_{0}^{t/\varepsilon}f\circ\Phi^s\, ds
\right)_t\stackrel{distrib}{\longrightarrow}(B_{\Delta_\Phi^{-1}(t)})_t\, ,\quad\mbox{as }n\rightarrow +\infty\, .
\end{equation*}
But $B_{\cdot/\mu(M)}$ is a Brownian motion of variance $\sigma^2_{U_g}$, and
it follows from Remark~\ref{remchgtinvDelta}
that $\Delta_2:=\Delta_\Phi(\mu(M)\cdot)$ is a $\beta$-stable
L\'evy process with characteristic function at time 1 given by $\Psi_{\beta,(2\sum_{j=1}^{\mathcal J}\mathfrak c_j)^{\frac 1\beta}(2/\beta)^{\frac 1\alpha} I_{\frac 1\beta-1}}$,
ending the proof 
of the last part of Theorem~\ref{thm:TCLflow2}.
For the first part, Proposition~\ref{prop:CLTflow2b} combined with Remark~\ref{remchgtinvDelta} ensures that
\begin{equation*}
\left(\varepsilon^{\beta}\int_{0}^{t/\varepsilon}f\circ\Phi^s\, ds
\right)_t\stackrel{distrib}{\Longrightarrow}\left(\frac{\int_{\mathcal M_0}f\, d\nu}{\mu(M)}\Delta_\Phi^{-1}(t)=\int_{\mathcal M_0}f\, d\nu\, \Delta_2^{-1}(t)\right)_t\, ,\quad\mbox{as }n\rightarrow +\infty\, ,
\end{equation*}
strongly in distribution with respect to $\nu$ (with respect to the uniform topology).
This ends the proof of Theorem~\ref{thm:TCLflow2}.

\part{Appendices} \label{appendices}
\appendix

\section{Proof of Proposition \ref{prop3}}\label{proofprop3}
For simplicity, we fix any $i \in \{1,...,\mathcal J\}$, consider a point $x\in M_N^{(i)}$ and drop the index $i$ in the following analysis. 
We recall that, by \textbf{(h1)}, the boundaries can be represented in a Cartesian coordinate system $(t, z)$ by the graphs $\Gamma_{\pm}$ given by
\beq\label{zsP1}
z = \pm f_{\mathfrak{c}}(t), \quad \text{with} \ f_{\mathfrak{c}}(t) = \frac{\mathfrak{c}}{\beta t^{\beta}}, \quad \forall t \in [\widehat{t}, \infty)\, ,
\eeq
with here $\mathfrak c=\mathfrak c_i$. 
We use the notations introduced before
Proposition~\ref{prop3} and set also
$$ s_n = 1/t_n.$$
For any $n$,
\beq \label{alphansn}
\alpha_n = \arctan({\mathfrak{c}} s_n^{\beta+1})= \mathfrak{c} t_n^{-\beta-1} - \tfrac{1}{3} \alpha_n^3 + \cO(\alpha_n^5), 
\eeq
i.e.
\begin{equation}\label{snfirstorder}
\mathfrak{c} t_n^{-\beta-1} = \tan (\alpha_n) = \alpha_n \left(1 + \alpha_n^2/3 + \cO(\alpha_n^5)\right).
\end{equation}
Consider the two following quantities
\beq\label{uwn}
u_n=\frac{\alpha_n}{\alpha_{n+1}},\,\,\,\,\,w_n=\frac{ v_n}{\alpha_{n}}\, .
\eeq
Note that $u_n\ge 1$ for all $n=1,...,N_2-1$.
\begin{itemize}
\item \underline{Preliminary estimates.}
One can check that $v_n>0$ defines an increasing sequence that slowly grows from about 0 to about $\pi/2$ for $n=1,\cdots, N_2-1$. 
At the $(n+1)$-th reflection in the cusp, the angle between the $t$-axis and the incident direction is
$v_{n}=\gamma_n+\alpha_n$ and $\gamma_{n+1}-\alpha_{n+1}$, which leads to
\beq\label{3.6}
\gamma_{n+1}=\gamma_{n}+\alpha_n+\alpha_{n+1},\eeq
i.e.
\beq\label{3.7}
v_{n+1}=v_n+2\alpha_{n+1}=v_n+2\mathfrak c t_{n+1}^{-\beta-1}-2\, \alpha_{n+1}^3/3+\cO(\alpha_n^5)
\, ,\eeq
where we used~\eqref{snfirstorder}. 
Since $\alpha_n$ is decreasing, this implies that 
\begin{equation}\label{nalphan}
2(n-1)\alpha_n\le v_n-v_1\le v_n<\frac \pi 2\, ,
\end{equation}
and so that \beq\label{firstests}
\alpha_n=\cO(n^{-1}),\quad s_n:=t_n^{-1}=\cO(\alpha_n^{\frac 1{\beta+1}})=\cO(n^{-\frac 1{\beta+1}}) \, .\eeq
Moreover $\alpha_n\le \alpha_1\le v_1=\alpha_1+\gamma_1$, so it follows from~\eqref{nalphan} that
\beq
\forall n=2,..., N_2,\quad 
\alpha_n\cot v_n\le  \frac{\alpha_n}{\sin v_n}\le \frac{v_n}{2(n-1)\sin v_n}\le \frac{\pi}{4(n-1)}\, ,\label{alphancotvn-1}
\eeq
and 
\beq 
\alpha_1\cot v_1\le \frac{v_1}{\sin v_1}\le \frac \pi 2.
\eeq
This implies that
\begin{equation}\label{alphancotv}\alpha_n \cot v_n\le\alpha_n/\sin v_n=\cO(n^{-1}),\,\,\,\,\alpha_n=\cO(n^{-1}\sin v_n )\, .
\end{equation}

Moreover
\beq\label{3.71}
 t_{n+1}=t_n+\tau_n \cos v_n,
\eeq
where $\tau_n$ is the free path between two collisions based at $(t_n,\pm f_{\mathfrak c}(t_n))$ and $(t_{n+1},\mp f_{\mathfrak c}(t_{n+1}))$:
\beq\label{deftaun}
\tau_n:=\tau(r_n,\varphi_n)=\mathfrak c\frac{t_{n+1}^{-\beta}+t_n^{-\beta}}{ \beta \sin v_n}\, .
\eeq
Note that $(\tau_n)_n$ is a strictly decreasing sequence (since $(t_n)_n$ and $(v_n)_n$ are increasing).
It follows from~\eqref{3.71} and~\eqref{deftaun} that
\beq\label{3.71bb}
 t_{n+1}-t_n=\frac{\mathfrak c}\beta (t_{n+1}^{-\beta}+t_n^{-\beta})\cot v_n\, .
\eeq

\item \underline{First Collision in the cusp for $N\ge 2$: Control of $v_1$ and $w_1$.}\\
The goal of this part is to establish bounds for the entrance velocity $v_0=v_1-2\alpha_1$ in the cusp in terms of $\alpha_1$ (i.e. of $t_1$). More precisely we will prove that $-\alpha_1\le v_0=\cO(\alpha_1)$, this will ensure that $v_1\approx \alpha_1$.\\
Since $\gamma_1\ge 0$, it follows that $v_0=\gamma_1-\alpha_1\ge -\alpha_1$
 (the limit case $v_0= -\alpha_1$, corresponds to the case where the first collision in the cusp is tangent to the boundary at $(t_1,\pm f_{\mathfrak c}(t_1))$). This gives the lower bound on $v_0$. 
 Now let us establish an upper bound on $v_0$. This will follow from Items~(a)--(c) below.
\begin{itemize}
\item[(a)] Let  $t'_0,t'_1$ be two positive real numbers such that the tangent line to $\Gamma_\mp$ at $(t'_0,\mp f_{\mathfrak c}(t'_0))$ contains $(t'_1,\pm f_{\mathfrak c}(t'_1))$.
It follows from~\eqref{3.71bb} that
$$t'_1-t'_0= \frac{\mathfrak c}\beta\cot(v'_0)({t'_0}^{-\beta}+{t'_1}^{-\beta})\, .$$
Since $\cot(v'_0)=|f'_{\mathfrak c}(t'_0)|^{-1}=(\mathfrak c(t'_0)^{-\beta-1})^{-1}$, 
it follows that
$$
t'_1-t'_0= 
\frac{(t'_0)^{-\beta}+{t'_1}^{-\beta}}{\beta (t'_0)^{-\beta-1}}
=\frac{t'_0}{\beta}(1+(t'_0/t'_1)^{\beta})\, .$$
Thus
\begin{equation}\label{t1-t0}
t'_1/t'_0-1= \frac 1{\beta}(1+(t'_0/t'_1)^{\beta})\, .
\end{equation}
We infer that $t'_0=\widetilde ct'_1$  where $\widetilde c$
is the unique solution in $(0,1)$ of
$$
{\widetilde c}^{-1}-1= \frac 1{\beta}(1+\widetilde c^{\beta})\, .
$$
\item[(b)] This ensures in particular that, if $t_1\ge \widehat t/\widetilde c$, then there exists $t'_0:=\widetilde ct_1\in[\widehat t,+\infty)$ such that
the tangent line to $\Gamma_\mp$ at $(t'_0,\mp f_{\mathfrak c}(t'_0))$ contains $(t_1,\pm f_{\mathfrak c}(t_1))$. 
For this choice of $t'_0$, let us write $\alpha'_0$ for the angle made by the tangent vector $(1,\mp f'_{\mathfrak c}(t'_0))$ to $\Gamma_\mp$ with $(1,0)$.
We then have
\begin{equation}\label{majov0}
v_0\le \alpha'_0:=
\arctan(\mathfrak c(t'_0)^{-\beta-1})
=\arctan(\widetilde c^{-\beta-1}\tan\alpha_1)\, .
\end{equation}
We conclude that $v_1\le 2\alpha_1+ \alpha'_0=\mathcal O(\alpha_1)$, if $t_1\ge \widehat t/\widetilde c$.
\item[(c)] 
If $t_1\in[\widehat t,\widehat t/\widetilde c)$, then $t'_0$ does not exist but $\alpha_1\ge \arctan(\mathfrak c (\widehat t/\widetilde c)^{-\beta-1})>0$ and so $v_1\le\frac\pi 2\le \frac \pi {2 \arctan(\mathfrak c (\widehat t/\widetilde c)^{-\beta-1})}\alpha_1$.
\item[(d)]It follows from the above Item (a) that for all $n\ge 2$, 
\begin{equation}\label{borntntn-1}
\widetilde c t_{n}\le t_{n-1}\le t_n
\end{equation}
Indeed we proved in Item (a) that $(t_n,\pm f_{\mathfrak c}(t_n))$ is on the tangent line to the graph of $\mp f_{\mathfrak c}$ at the point $(\widetilde ct_n,\mp f_{\mathfrak c}(\widetilde ct_n))$ and $f_{\mathfrak c}$ is convex. This implies in particular that
\beq\label{bound_un}
u_n\approx 1\, .
\eeq
\end{itemize}
It follows from Items (a)--(c) that 
\begin{equation}\label{w_1}
v_0=\mathcal O(\alpha_1),\quad 
\alpha_1\le v_1=\mathcal O(\alpha_1), \quad \mbox{and so}\quad \boxed{w_1\approx 1}\, .
\end{equation}
\item \underline{Definition and study of $D_n$.}
We define 
\beq\label{tnn1}
{\boxed{D_n=-\beta \frac{t_{n+1}-t_n}{t_n},\quad\mbox{i.e.}\quad\frac{t_{n+1}}{t_n}=1-\frac{D_n}{\beta}\, .}}
\eeq
It follows from~\eqref{3.71bb},~\eqref{borntntn-1},~\eqref{alphansn} and~\eqref{alphancotv} that
\begin{equation}\label{MajoDn}
-D_n=\frac{\mathfrak c\cot v_n(t_{n+1}^{-\beta}+t_n^{-\beta})}{t_n}\le 2\mathfrak c t_n^{-\beta-1}\cot v_n=\cO\left(\alpha_n\cot v_n\right)=\cO(n^{-1})\, .
\end{equation}
This combined with~\eqref{firstests} ensures the following estimate that we will use throughout the proof
\beq\label{boundDnan}
\boxed{|D_n|+\alpha_n=\cO(n^{-1})\, .}
\eeq
A Taylor expansion gives
\beq\label{Taylortnb}
\frac{t_{n+1}^{-\beta}+t_n^{-\beta}}{t_n^{-\beta}}=\left(
1-\frac{D_n}\beta\right)^{-\beta}+1
=2+
D_n+\frac{(\beta+1) D_n^2}\beta+\cO(D_n^3)\, .
\eeq
Combining~\eqref{MajoDn} with \eqref{Taylortnb}, we obtain
\begin{align}
\nonumber -D_n
 &=2\mathfrak c t_n^{-\beta-1} \cot v_n\left(1+D_n/2+\frac{(\beta+1)D_n^2}{2\beta}+\cO(D_n^3)\right)\, .\label{3.711}
 \end{align}
This also implies that
\beq\label{3.7111}
2\mathfrak c t_n^{-\beta-1} \cot v_n=\frac{-D_n(1+\cO(D_n^3))}{1+D_n/2+\frac{(\beta+1)D_n^2}{2\beta}}
=-D_n\left(1-\frac{D_n}2\right)+\mathcal O(D_n^{3})\, ,
\eeq
which combined with~\eqref{alphansn} and~\eqref{boundDnan} leads to
\beq\label{3.7111bis}
\boxed{ 2 \alpha_n \cot v_n=-D_n\left(1-\frac{D_n}2\right)+\mathcal O(n^{-3})\, .}
\eeq
\item \underline{First study of $u_n-1$.}
It follows from~\eqref{tnn1} and~\eqref{borntntn-1} that
$0\le -\frac{D_n}\beta=\frac{t_{n+1}-t_n}{t_n}\le\widetilde c^{-1}-1$. 
Taylor expansions ensure that
\begin{align}
    \ln t_{n+1}^{-\beta}-\ln t_n^{-\beta}&=-\beta\ln\left(1-\frac{D_n}{\beta}\right)=D_n+ \frac{D_n^2}{2\beta}+\cO(D_n^3)\, ,\label{difflnftn}
\end{align}
and that
\begin{align}
\nonumber\frac{t_{n+1}^{-\beta-1}}{t_n^{-\beta-1}}-1&=\left(1-\frac{D_n}\beta\right)^{-\beta-1}-1\\
&=\frac{\beta+1}{\beta}\left(D_n+\frac{(\beta+2)D_n^2}{2\beta}+\cO(D_n^3)\right)\, .\label{tnbeta+1}
\end{align}
This combined with~\eqref{alphansn} 
implies that, for $n\in [1,N_2-1]$,
\begin{align}\label{un-1}
   u_n^{-1}-1&= \frac{\alpha_{n+1}}{\alpha_n}-1 =\frac{t_{n+1}^{-\beta-1}}{t_{n}^{-\beta-1}}\(1-\frac{\alpha_{n+1}^2-\alpha_n^2}3+\cO(\alpha_n^4)\)-1=\mathcal O(n^{-1})\, .
\end{align}
Thus, using $u_n\approx 1$ proved in~\eqref{bound_un}, it follows that
\begin{equation}\label{diffalphan}
\alpha_{n+1}^2-\alpha_n^2=\alpha_n^2(u_n^{-2}-1)
=\cO(n^{-3})\, .
\end{equation}
Using this rough estimate and~\eqref{tnbeta+1} 
in the right hand side of~\eqref{un-1}, we obtain
\begin{align*}
   u_n^{-1}-1&=\frac{\beta+1}{\beta}D_n\left(1+\frac{(\beta+2)D_n}{2\beta}
   \right)+\cO(n^{-3})\, ,
\end{align*}
where we used again~\eqref{boundDnan}. 
This implies that 
\begin{align}
u_n-1&=-\frac{\beta+1}{\beta}\left(D_n
-\frac{D_n^2}{2}
+\cO(n^{-3})\right)
\, .
\label{A15}
\end{align}
This combined with~\eqref{3.7111bis} leads to
\begin{align}
 u_n-1&=\frac{\beta+1}{\beta}2\alpha_n\cot v_n
+\cO\left(
n^{-3}\right)\label{un-1ba1}\\
&=\frac{v_n}{w_n}\frac {2(\beta+1)}{\beta\tan v_n}
+\cO\left(
n^{-3}\right)\, ,\label{un-1ba}
\end{align}
where we use the fact that $\alpha_n=v_n/w_n$
by definition of $w_n$.

\item \underline{Study of $w_n=v_n/\alpha_n$.}
Recall that 
\begin{equation}\label{minown}
w_n\ge 2(n-1)
\end{equation}
(as seen in \eqref{nalphan}).
Moreover $w_1=\frac{\gamma_1+\alpha_1}{\alpha_1}\ge 1$. 
 Combining the previous estimate~\eqref{un-1ba} of $u_n-1$ with \eqref{3.7} and~\eqref{uwn} yields:
 \begin{align}
\nonumber w_{n+1}&=\frac{v_{n+1}-2\alpha_{n+1}}{\alpha_{n+1}}+2=\frac{v_n}{\alpha_{n+1}}+2=w_n u_n+2=w_n+2+w_n(u_n-1)\\
&\nonumber = w_n+2+\frac{2(\beta+1)}\beta\frac{v_n}{\tan v_n}+\cO(w_n/n^3)
\nonumber\\
&=w_n(1
+\mathcal O(n^{-3}))+2+ \frac{2\beta+2}\beta \cos v_n+\cO\left(v_n^2\cos v_n\right)\\
&=w_n(1
+\mathcal O(n^{-3}))+ \frac{4\beta+2}\beta +\cO\left(v_n^2
\right)
\, .
\label{wnEn00pre}
\end{align}
Thus $w_{n+1}\le w_n(1+\cO(n^{-2}))+\cO(1)$, which by induction leads to
\begin{align*}
w_n&\le \cO(1+w_1) \sum_{k=1}^n\prod_{j=k}^n(1+\cO(j^{-2}))
\le \cO(1)\sum_{k=1}^ne^{\sum_{j=k}^n\mathcal O(j^{-2})}=\cO(n)\, ,
\end{align*}
where we used the fact that $w_1=\cO(1)$ (established in~\eqref{w_1}). It follows from this combined with~\eqref{minown} that
\beq \label{vn&nalphan}
\boxed{w_n\approx n,\,\,\,\, v_n\approx n\alpha_n\, .}
\eeq
Thus~\eqref{wnEn00pre} becomes
\begin{equation}\label{wnEn00}
w_{n+1}=w_n 
+\frac{2(2\beta+1)}\beta +\cO\left(n^{-2}+v_n^2\right)\, .
\end{equation}
Moreover, since the sequence $(v_n)_{n=1,...,N_1}$ is increasing, $\sum_{k=1}^n v_k^2\le nv_n^2$ and so
 \begin{align}
 w_n&=w_1 +\frac{2(2\beta+1)}\beta n +\cO\left(1 + n v_n^2\right)
 \\
 &=\frac{2(2\beta+1)}\beta n
 \left(1+\cO\left( n^{-1}+v_n^2\right)\right)\, ,\label{wnEn1}
\end{align}
where we used again the fact that $w_1=\cO(1)$ (coming from~\eqref{w_1}).

\item \underline{Proof of Items~(2), (3) and of the first part of  (4).}
It follows from~\eqref{wnEn1} (and~\eqref{vn&nalphan}) that
\begin{align}
\nonumber 2\alpha_n\cot v_n&=\frac {2v_n}{w_n}\cot v_n=\frac{2+\cO(v_n^2)}{w_n}=\frac{\beta}{(2\beta+1) n}\left(1
+\cO\left(\frac{1}n+v_n^2
\right) \right)
\, .
\label{cotvn}
\end{align}
Combining this estimate with~\eqref{un-1ba1},
we obtain that
\begin{equation}\label{un-1c}
u_n-1=\frac {\beta+1}{(2\beta+1)n} \left(1+\cO\left(\frac{1}n+v_n^2 \right)\right)\, .
\end{equation}
Since $u_n\approx 1$ (due to~\eqref{bound_un}), writing 
$\log u_n=\log (1+u_n-1)=u_n-1+\cO((u_n-1)^2)$, it follows that
\beq\label{alphann=1}
\log \frac{\alpha_n}{\alpha_{n+1}}=\log u_n=\frac {\beta+1}{(2\beta+1)n}\left(1
+\cO\left(\frac{1}n+v_n^2
\right)\right)\, .
\eeq
In particular, for $k<n$,
\begin{equation}
\log \frac{\alpha_k}{\alpha_{n+1}}
=\sum_{m=k}^n\log \frac{\alpha_m}{\alpha_{m+1}}=\frac {\beta+1}{(2\beta+1)}\log (n/(k+1))+\cO\left(\frac
1{k+1}+v_n^2
\right)\, ,\label{ratioalpha}
\end{equation}
where we used the fact that $\sum_{m=k}^n\frac{1
} {m^2}\le \frac{
1
}k$
and that
\beq
\sum_{m=k}^n\frac{v_m^2}m\le v_n \sum_{m=k}^n\frac{v_m}m\approx  v_n \sum_{m=k}^n\alpha_m= v_n\frac{v_n-v_{k-1}}2\le \frac{v_n^2}2\, .
\eeq

Let $i\in\{1,2\}$.
It follows from  \eqref{vn&nalphan} and from the definition of $N_i$ that
\beq
\label{alphaNi}  
\alpha_{N_i}\approx v_{N_i}/N_i\approx N_i^{-1}\, .
\eeq
Applying \eqref{ratioalpha} to couples $(k,n)=(m,N_i)$,
we obtain that
\beq\label{alpham}
\alpha_{m}
\approx \alpha_{N_i}(N_i/m)^{\frac{\beta+1}{2\beta+1}}\approx m^{-\frac{\beta+1}{2\beta+1}}
N_i^{-\frac{\beta}{2\beta+1}}\, ,\quad m\in[1,N_i]\, .
\eeq
This, combined with~\eqref{un-1c}
proves Items (2) and the first part of (4), 
since $\alpha_n\approx t_n^{-(\beta+1)}$. Item (3) then follows from Item~(2) combined with~\eqref{defvn},~\eqref{w_1} and~\eqref{vn&nalphan}.

\item \underline{Proof of Item (1).}
Note that \eqref{alpham} applied to $(N_1,N_2)$ combined with \eqref{alphaNi} implies that $N_1\approx N_2$. Furthermore, 
due to \eqref{3.7} and \eqref{alphaNi},
\beq\nonumber
 N_2^{-1}
 \approx\alpha_{N_2}\le \frac{v_{N_2}-v_{N_1}}{2(N_2-N_1)}\le \alpha_{N_1}
 \approx N_1^{-1}\, .
\eeq
Thus $N_1\approx N_2\approx N_2-N_1$ and we have proved Item (1).

\item \underline{Proof of the first part of Item~(7): proof of~\eqref{controlHn}.} 
We study $H_n=s_n^{\beta}\sin v_n$. 
Combining~\eqref{3.7111} and \eqref{tnbeta+1}, we obtain
\beq\label{alphan+1vn}
2 \mathfrak c t_{n+1}^{-\beta-1}\cot v_n=-D_n-\frac{\beta+2}{2\beta} D_n^2+\cO(D_n^3) 
\eeq
and also, using~\eqref{alphansn} and~\eqref{boundDnan}.
\begin{align}\nonumber
2& \alpha_{n+1}\cot v_n
=\left(-D_n-\frac{\beta+2}{2\beta} D_n^2+\cO(D_n^3)\right)\left(1-\alpha_{n+1}^2/3+\cO(\alpha_n^4)\right)\\
&=-D_n-\frac{\beta+2}{2\beta}D_n^2+\cO(n^{-3})
\,.
\label{alphan+1vn1}  
\end{align}
Observe that
\begin{align}\nonumber
\frac{\sin v_{n+1}}{\sin v_n}-1&=\cot v_n(v_{n+1}-v_n) -\frac 12(v_{n+1}-v_n)^2+\cO\left(\cot v_n(v_{n+1}-v_n)^3\right)\\
&=2\alpha_{n+1}\cot v_n -2\alpha_{n+1}^2+\cO\left(\alpha_n^3\cot v_n\right)\nonumber
\\
&
=-D_n-\frac{\beta+2}{2\beta}D_n^2-2\alpha_{n+1}^2+\cO\left(n^{-3}\right)
\, ,\label{rationsinvn}
\end{align}
where we used~\eqref{alphan+1vn1},~\eqref{firstests} and~\eqref{alphancotv}. 
It follows from~\eqref{3.7} that $v_n\le v_{n+1}=v_n+2\alpha_{n+1}\le v_n+2\alpha_n\le 2v_n$ and so
using~\eqref{rationsinvn} combined with a Taylor expansion formula, we get 
\begin{align}
\ln & \frac{\sin v_{n+1}}{\sin v_n}
=\left(\frac{\sin v_{n+1}}{\sin v_n}-1\right)-\frac{1}{2}\left(\frac{\sin v_{n+1}}{\sin v_n}-1\right)^2+\cO\left(\left(\frac{\sin v_{n+1}}{\sin v_n}-1\right)^3\right)\nonumber\\
&
=-D_n-\frac{\beta+2}{2\beta}D_n^2-2\alpha_{n+1}^2-\frac{D_n^2}2
+\cO(n^{-3})\nonumber\\
&=-D_n-\frac{\beta+1}\beta 
D_n^2-2\alpha_{n+1}^2+\cO(n^{-3})\, ,\label{lnsinvn}
\end{align}
where we used~\eqref{boundDnan}.
Summing~\eqref{lnsinvn} with~\eqref{difflnftn}, we obtain for $n=1,\cdots, N_2$,
\begin{align}\label{eqA34}
&\ln \frac{\sin v_{n+1}}{\sin v_n}+\ln \frac{t_{n+1}^{-\beta}}{t_n^{-\beta}}=-
\frac{2\beta+1}{2\beta}
D_n^2-2\alpha_{n+1}^2+\cO(n^{-3})\, .
\end{align}
Recalling that $s_n=t_n^{-1}$, we take exponential from both sides of the above equation to get 
$$\frac{H_{n+1}}{H_n}=\frac{s_{n+1}^\beta\sin v_{n+1}}{s_n^\beta\sin v_n}= \exp{\left(-\frac{2\beta+1}{2\beta}D_n^2-2\alpha_{n+1}^2+\cO(
n^{-3})\right)}\, .$$ 
This leads to
\begin{align}
\label{sinvnftn}
A_n&:=\frac{H_{N_2}}{H_n}= \exp\left(-\sum_{m=n}^{N_2-1}\left(\frac{2\beta+1}{2\beta}D_m^2+2\alpha_{m+1}^2+\cO(m^{-3})\right)\right)\\
\nonumber &= \left(1-\sum_{m=n}^{N_2-1}\left(\frac{2\beta+1}{2\beta}D_m^2+2\alpha_{m+1}^2+\cO(m^{-3})\right)+\cO\left(\left(\sum_{m= n}^{N_2-1}m^{-2}\right)^2\right)\right)\\
&=1-\sum_{m=n}^{N_2-1}\left(\frac{2\beta+1}{2\beta}D_m^2+2\alpha_{m+1}^2\right)+\cO(n^{-2}-N_2^{-2})=1+\cO(n^{-1}-{N_2}^{-1})\, ,
\label{Anexpansion}
\end{align}
using~\eqref{boundDnan}. Thus we have proved the first part of~\eqref{controlHn}.
The second part of~\eqref{controlHn} will follow using~\eqref{HN2} combined with Item~(4).
\item \underline{Proof of the second part of Item~(7): proof of~\eqref{HN2}}.\\
We wish to establish a precise formula for $H_{N_2}$.
Recall from~\eqref{sinvnftn}
that
\begin{align}\label{Cnvnalpha}
 (H_{N_2}/ A_n)^{\frac{\beta+1}{\beta}}=\left(\frac{\sin v_n}{ t_n^\beta}\right)^{\frac{\beta+1}{\beta}}=\left(\sin v_n\right)^{\frac{\beta+1}{\beta}} \frac{\alpha_n}{\mathfrak c}
 \left(1+ \cO(\alpha_n^2)
 \right)\, ,
\end{align}
where we used~\eqref{alphansn} for the last estimate.
Using Items (2) and (3), we know that
\beq\label{error11}
\sum_{n=2}^{N_2}\left(\sin v_n\right)^{\frac{\beta+1}{\beta}} \alpha_n^3
=   \sum_{n=2}^{N_2}\cO\left(n^{-\frac{2\beta+2}{2\beta+1}} N^{-\frac{4\beta+1}{2\beta+1}}\right)\leq C   N^{-\frac{4\beta+1}{2\beta+1}}\, .
\eeq
Summing~\eqref{Cnvnalpha} combined with~\eqref{error11} and applying a series-integral comparison, we obtain that
\begin{align}
\nonumber \sum_{n=2}^{N_2}& 2( H_{N_2}/ A_n)^{\frac{\beta+1}{\beta}}
=\sum_{n=2}^{N_2}\left(
\sin v_n
\right)^{\frac{\beta+1}{\beta}} \frac{v_n-v_{n-1}}{\mathfrak c}+\cO(N^{-\frac{4\beta+1}{2\beta+1}})\\
&=\mathfrak c^{-1}\left(I_{1/\alpha}-\int_0^{v_1}  (\sin v)^{\frac{\beta+1}{\beta}} \, dv
-\int_{v_{N_2}}^{\frac\pi 2}  (\sin v)^{\frac{\beta+1}{\beta}} \, dv
\ -J_{N_2} \right)+\cO(N^{-\frac{4\beta+1}{2\beta+1}})
\, ,
\label{CNANCheck}
\end{align}
where 
\begin{equation}\label{I1aJ2}
I_{1/\alpha}:= \int_{0}^{\frac \pi 2} (\sin v)^{\frac{\beta+1}{\beta}} \, dv\quad\mbox{and}\quad    
J_{N_2}:=\sum_{n=2}^{N_2} \int_{v_{n-1}}^{v_{n}} \left((\sin v)^{\frac{\beta+1}{\beta}} - (\sin v_n)^{\frac{\beta+1}{\beta}} \right) \, dv \, .
\end{equation}
It follows from Item (3) that
\begin{equation}\label{intsmallv}
\int_0^{v_1}  (\sin v)^{\frac{\beta+1}{\beta}} \, dv=\mathcal O\left(v_1^{1+\frac{\beta+1}\beta}\right)=\mathcal O\left(v_1^{\frac{2\beta+1}\beta}\right)=\mathcal O\left(N^{-1}\right)\, .
\end{equation}
Furthermore, since $v_{N_2}\le \pi/2<v_{N_2}+2\alpha_{1+N_2}$,
\begin{equation}\label{intbigv}
\int_{v_{N_2}}^{\frac \pi 2}  (\sin v)^{\frac{\beta+1}{\beta}} \, dv
\le\frac \pi 2-v_{N_2}=\mathcal O\left(\alpha_{1+N_2}\right)=\mathcal O\left( N^{-1}\right)
\, ,
\end{equation}
using again Item (2) for the estimate of $\alpha_{1+N_2}$.
Combining~\eqref{CNANCheck} with~\eqref{intsmallv} and~\eqref{intbigv}, we obtain that
\begin{align}
 \sum_{n=2}^{N_2} 2( H_{N_2}/ A_n)^{\frac{\beta+1}{\beta}}
&=\mathfrak c^{-1}\left(I_{1/\alpha}-J_{N_2}\right) +\cO(N^{-1})\, .
\label{CNANCheck1}
\end{align}

\underline{Study of the left hand side of~\eqref{CNANCheck1}.\\}
It follows from~\eqref{Anexpansion} that
\begin{align}\label{An-b/b+1}
A_n^{-\frac{\beta+1}\beta}&=1+\frac{\beta+1}\beta \sum_{m=n}^{N_2-1}\left(\frac{2\beta+1}{2\beta}D_m^2+2\alpha_{m+1}^2\right)+\cO(n^{-2}-N_2^{-2})
\, ,
\end{align}
and so
\begin{align*}
\sum_{n=2}^{N_2}A_n^{-\frac{\beta+1}\beta}&=(N_2-1)+\frac{\beta+1}\beta
\sum_{m=2}^{N_2}(m-1)\left(\frac{2\beta+1}{2\beta}D_m^2+2\alpha_{m+1}^2\right)
+\cO(1)\, .
\end{align*}
But, since $m\alpha_m\le v_m$ and $\alpha_{m+1}=\cO(\alpha_m)$,
\beq
\sum_{m=2}^{N_2}m\alpha_{m+1}^2=\cO\left( \sum_{m=2}^{N_2}v_m\alpha_{m}\right)=\cO\left(\int_{v_1}^{v_{N_2}}v\, dv\right)=\cO(1)\, .
\eeq
Thus
\beq
\sum_{n=2}^{N_2} A_n^{-\frac{\beta+1}\beta}=N_2+\cO(1)+\frac{\beta+1}\beta\frac{2\beta+1}{2\beta}\sum_{m=2}^{N_2}m D_m^2\, .\label{sumAn_bis}
\eeq
Furthermore, it follows from~\eqref{wnEn1} that
\begin{align}
\frac{2(2\beta+1)}\beta\nonumber\sum_{m=2}^{N_2}m D_m^2
&=
\sum_{m=2}^{N_2}\left(w_m
+\cO\left(
1+mv_m^2\right)\right)D_m^2\\
&=\frac {2\beta}{2\beta+1}\log N_2+
\cO\left(1\right)\, ,\label{summDm^2} 
\end{align}
where we used the two following estimates
\beq 
\sum_{m=2}^{N_2}mv_m^2D_m^2=
\cO\left(\sum_{m=2}^{N_2}v_m^2/m\right)=
\cO\left(\sum_{m=2}^{N_2}v_m\alpha_m\right)=\cO\left(\int_{v_1}^{v_{N_2}}v\, dv\right)=\cO(1)\, ,
\eeq
and
\begin{align}
\nonumber\sum_{k=2}^{N_2}w_kD_k^2&= \left(\sum_{k=2}^{N_2}2v_k(\cot v_k)^2(2\alpha_k)\right)+\cO(1)\\
\nonumber&=2 \int_{v_1}^{v_{N_2}}v^{-1}\, dv +\cO(1)=2(\log v_{N_2}-\log v_1)+\cO(1)\\
&=\frac{2\beta}{2\beta+1}\log N+\cO(1)\, ,\label{sumwmDm^2}
\end{align}
using the definition of $w_k$,~\eqref{3.7111bis}, Item (3)  and
the following computation
\begin{align*}
\int_{v_{k-1}}^{v_k}(v_k(\cot v_k)^2-v^{-1})\, dv&\le \sup_{v\in[v_{k-1};v_k]}\left[\frac{(v_k\cot v_k)^2
     -1}{v_k}+v_k^{-1}-v^{-1}\right]2\alpha_k\\
     &= \cO\left(\frac{v_k^2}{v_k}+\frac{v_k-v_{k-1}}{v_{k-1}v_k}\right)2\alpha_k\\
     &= \cO\left(v_k\alpha_k+\frac{\alpha_k^2}{v_k^2}\right)=\cO\left(\alpha_k+k^{-2}\right)\, ,
\end{align*}
(since $v_k\cot v_k=1+\cO(v_k^2)$), 
and $\sum_{k=2}^{N_2}\alpha_k=\frac{v_{N_2}-v_1}2=\cO(1)$.

It follows from~\eqref{sumAn_bis}
and~\eqref{summDm^2} that
\beq
2\sum_{n=2}^{N_2} A_n^{-\frac{\beta+1}\beta}=N+
\frac{\beta+1}{2\beta+1}
\log N+\cO(1)\, .\label{sumAn_final}
\eeq
\underline{Study of the right hand side of~\eqref{CNANCheck1}.}
To estimate $J_{N_2}$ defined in~\eqref{I1aJ2}, we observe that
\begin{align*}
\left((\sin v)^{\frac{\beta+1}{\beta}} - (\sin v_n)^{\frac{\beta+1}{\beta}} \right) 
&=\tfrac{\beta+1}{\beta}(\sin v_n)^{\frac {\beta+1}\beta}\cot v_n (v-v_n)+\cO\left((\sin v_n)^{\frac{\beta+1}\beta}(\alpha_n/\sin v_n)^2\right)\\
&=\tfrac{\beta+1}{\beta}(\sin v_n)^{\frac {\beta+1}\beta}\left((v-v_n)\cot v_n+\cO\left(n^{-2}\right)\right)
\, ,
\end{align*}
and so
\begin{align*}
\int_{v_{n-1}}^{v_n}
\left((\sin v)^{\frac{\beta+1}{\beta}} - (\sin v_n)^{\frac{\beta+1}{\beta}} \right)\, dv
&=\tfrac{\beta+1}{\beta}(\sin v_n)^{\frac {\beta+1}\beta}\left(-\tfrac 12(v_{n-1}-v_n)^2\cot v_n+\cO\left(\alpha_nn^{-2}\right)\right)\\
&=\tfrac{\beta+1}{\beta}(\sin v_n)^{\frac {\beta+1}\beta}\left(-2\alpha_{n}^2\cot v_n+\cO\left(\alpha_nn^{-2}\right)\right)\\
&=\tfrac{\beta+1}{\beta}(\sin v_n)^{\frac {\beta+1}\beta}\alpha_{n}\left(D_n+\cO\left(n^{-2}\right)\right)
\, ,
\end{align*}
where we used the fact that $v_{n}-v_{n-1}=2\alpha_n$ combined with~\eqref{3.7111bis}. 
Thus, using~\eqref{alphansn} and the definition of $A_n$ given in~\eqref{sinvnftn} and again~\eqref{boundDnan}, we obtain that
\begin{align}
\nonumber \mathfrak c^{-1}J_{N_2}&=\mathfrak c^{-1}\tfrac{\beta+1}{\beta}\sum_{n=2}^{N_2}
(\sin v_n)^{\frac{\beta+1}{\beta}}\alpha_n\left(D_n+\cO(n^{-2})\right)\\
\nonumber&=\frac{\beta+1}{\beta}\sum_{n=2}^{N_2}(\sin v_n /t_n^\beta)^{\frac{\beta+1}{\beta}}\left(1
+\cO(\alpha_n^2)
\right)\left(D_n+\cO(n^{-2})\right)\\
\nonumber&=\frac{\beta+1}{\beta} H_{N_2}^{\frac{\beta+1}{\beta}}\sum_{n=2}^{N_2}A_n^{-\frac{\beta+1}{\beta}} \left(D_n
+\cO(n^{-2})\right)\\
&=\frac{\beta+1}{\beta} H_{N_2}^{\frac{\beta+1}{\beta}}\left(\cO(1)+\sum_{n=2}^{N_2}A_n^{-\frac{\beta+1}{\beta}} D_n
\right)\, ,\label{controlIN2}
\end{align}
since $A_n\approx 1$ due to Items~(3) and~(4). 
It follows from~\eqref{Anexpansion} and~\eqref{3.7111bis} that
\begin{align}
\sum_{n=2}^{N_2}A_n^{-\frac{\beta+1}{\beta}} D_n
&=\sum_{n=2}^{N_2}\left(D_n+\cO\left(n^{-2}\right)\right)=\cO(1)+\sum_{n=2}^{N_2}D_n\\
&=\cO(1)-2\sum_{n=2}^{N_2}\alpha_n\cot v_n\\
&=\cO(1)-\int_{v_1}^{v_{N_2}}
\cot v\, dv+\sum_{n=2}^{N_2}\cO\left(\frac{\alpha_n^2}{(\sin v_n)^2}\right)\\
&=\cO(1)+\log \sin v_{1}
=\cO(1)-\frac{\beta}{2\beta+1}\log N
\, ,
\end{align}
since $\alpha_n/\sin v_n=\cO(v_n/(n\sin v_n))=\cO(n^{-1})$
and using also Item~(3) at the last line. 
So
\beq
\mathfrak c^{-1}J_{N_2}=-\frac{\beta+1}{2\beta+1} H_{N_2}^{\frac{\beta+1}\beta}\left(\log N+\cO(1)\right)\, .\label{{controlIN2} }
\eeq
It follows from~\eqref{CNANCheck1} combined with the previous estimate and with~\eqref{sumAn_final} that
\beq\label{CNI1}
 H_{N_2}^{\frac{\beta+1}{\beta}}\left(N 
 +\cO(1) \right)= \mathfrak c^{-1}I_{1/\alpha}+\cO(N^{-1})\, .
\eeq
and so
\beq\label{CNprecise}
\boxed{
 H_{N_2}^{\frac{\beta+1}{\beta}}= 
\frac{I_{1/\alpha}}{\mathfrak c N}
\left(1+\cO(N^{-1})\right)\, .}
\eeq

\item \underline{Proof of the second part of Item (4).} 
Using \eqref{deftaun}, we obtain that 
\beq\nonumber
\tau_n\approx \frac{t_n^{-\beta}}{\sin v_n}\approx \frac{t_n^{-\beta}}{n\alpha_n}\approx \frac{t_n^{-\beta}}{n t_n^{-\beta-1}}
\approx t_n/n
\approx   (n^2N^{-\frac 1{\beta+1}})^{-\frac \beta{2\beta+1}}
\, ,
\eeq
ensuring the second part of (4), 
and so
$$\sum_{n=1}^{N_1}\tau_n\approx   N^{\frac 1{\beta+1}}\approx t_{N_1},\quad \sum_{n=N_1}^{N_2}\tau_n\approx N^{\frac{1}{\beta+1}}\approx N\, f_{\mathfrak c}(t_{N_1})\, .$$
\item \underline{Proof of Items (5) and (6).}
For $n\in[2,N_1]$, the curvature $\cK_n=\cK(r_n)$ of the boundary at the position of $x_n$ satisfies
\beq\label{cKn}\cK_n=\frac{f''_{\mathfrak c}(t_n)}{(1+f'_{\mathfrak c}(t_n)^2)^{\frac 32}}=\frac{\mathfrak c(\beta+1)s_n^{\beta+2}}{(1+ \mathfrak c^2s_n^{2\beta+2})^{\frac 32}}=\mathfrak c(\beta+1)s_n^{\beta+2}(1+\cO(s_n^{2\beta+2}))\, .\eeq
It follows from~\eqref{cKn}, \eqref{deftaun}, \eqref{tnn1} and 
\eqref{alphansn}
that
\begin{align}\label{bartaun}
\frac{2\cK_n\tau_n}{\sin\gamma_n}&=\frac{2\mathfrak c^2(\beta+1) (t_n^{-2\beta-2}+t_n^{-\beta-2}t_{n+1}^{-\beta})}{\beta\sin\gamma_n\sin v_n}(1+\cO(s_n^{2\beta+2}))\nonumber\\
&= \frac{2(\beta+1)2\alpha_n^2}{\beta \sin^2 v_n
}(1+\cO(n^{-1}))\, ,
\end{align}
since $|D_n|+s_n^{2\beta+2}=\cO(n^{-1})$. 
It follows from
the definition of $w_n$ in~\eqref{uwn} and from \eqref{wnEn1} 
that
\begin{align}\label{bartaun_oldbis}
\frac{2\cK_n\tau_n}{\sin\gamma_n}
&= \frac{4(\beta+1)}{ \beta w_n^2} \left(1+\cO\left(n^{-1}
+v_n^2
\right)\right)\nonumber\\
&=\frac{4(\beta+1)}{4(2\beta+1)^2 n^2/\beta}\left(1+\cO\left(
n^{-1}+v_n^2
\right)\right)
\, .
\end{align}

This proves Item (5). For Item (6), using~\eqref{deftaun}, we write
$$\frac{\tau_{n+1}}{\tau_{n}}=\frac{s_{n+1}^{\beta}+s_{n+2}^{\beta}}{s_n^{\beta}+s_{n+1}^{\beta}}\cdot\frac{\sin v_n}{\sin v_{n+1}}.$$
Using  \eqref{alphansn},~\eqref{uwn},~\eqref{boundDnan} and \eqref{un-1c}, we obtain that, for $n\in[1,N_2]$,
\begin{align*}
E_1&:=\frac{s_{n+1}^{\beta}+s_{n+2}^{\beta}}{s_n^{\beta}+s_{n+1}^{\beta}}
=\frac{u_{n+1}^{-\frac{\beta}{\beta+1}}+1}{1+u_n^{\frac{\beta}{\beta+1}}}(1+\cO(\alpha_n^2))
=\frac{2-\frac{\beta}{\beta+1}(u_{n+1}-1)}{2+\frac{\beta}{\beta+1}(u_{n}-1)}(1+\cO(
n^{-2}))\\
&=1-\frac{\beta}{(2\beta+1)n}\left(1+\cO\left(\frac{
1}{n}+v_n^2\right)\right)\, .
\end{align*}
Moreover, for $n\in[1,N_2]$,
\begin{align*}
E_2&
:=\frac{\sin v_n}{\sin v_{n+1}}=\frac{\sin v_n}{\sin (v_n+2\alpha_{n+1})}\\
&=\frac{\sin v_n}{\sin v_n+2\alpha_{n+1}\cos v_n+\cO(\alpha_n^2\sin v_n)}\\
&=\frac{1}{1+2\alpha_{n+1}\cot v_n+\cO(\alpha_n^2)}
\end{align*}
and so
\begin{align*}
E_2
&=\frac{1}{1+\frac{2\alpha_{n+1}}{v_n}(1+\cO(v_n^2))+\cO(\alpha_n^2)}\\
&=1-\frac{2\alpha_{n+1}}{v_n}(1+\cO(v_n^2))+\cO(\alpha_n^2+n^{-2})\\
&=1-\frac{2\alpha_{n+1}}{v_n}\left(1 +\cO\left(v_n^2+n^{-1}\right)\right)\, ,
\end{align*}
using $\alpha_{n+1}\approx\alpha_n\approx v_n/n$.
This combined with \eqref{wnEn1} and \eqref{un-1c} (since $\frac{\alpha_{n+1}}{v_n}=(u_nw_n)^{-1}$ due to~\eqref{uwn}) leads to
\begin{align}
E_2
&=1-\frac {2\beta} {2(2\beta+1) n} \left(1+\cO\left(
\frac{
1}n+v_n^2\right)\right) 
,\,\,\,\,n=1,\cdots, [N/2]
\, .
\end{align}
Thus we have for $n\in[1,N_1]$,
\beq\label{ratiotau}\frac{\tau_{n+1}}{\tau_{n}}=E_1\cdot E_2=1-\frac{2\beta}{(2\beta+1)n}\left(1+\cO\left(\frac{
1}n+v_n^2
\right)\right) \, ,\eeq
which ends the proof of Item (6).
\end{itemize}
This ends the proof of Proposition~\ref{prop3}.

\section{Proof of Proposition~\ref{KBbd}}\label{appendKBbd}
As usual in the study of chaotic billiards, we use the $p$-length $|\cdot|_p$ corresponding to $|dx|_p=d(\sin\varphi)=\cos\varphi\, dr$ (in coordinates $x=(r,\varphi)$).
Proposition~\ref{KBbd} will follow from the next lemma. Indeed~\eqref{P1} follows from~\eqref{controllambda} while~\eqref{P3} corresponds to~\eqref{P2}.
\begin{lemma}\label{propdetails}
Assume $V\subset M_N$ is an unstable curve given by $\varphi=\phi_0(r)$. Then $F(V)=\cF^N(V)$ and 
\begin{equation}\label{AA'}
\forall 0\le k\le m\le N,\quad |\cF^m(V)|_p=\int_{\cF^k(V)} \left(\prod_{n=k}^{m-1}
\widetilde \lambda_n\right)\, |dx|_p
\, ,
\end{equation}
with
\beq\label{controllambda}
\widetilde\lambda_0 \approx N^{\frac \beta{(\beta+1)(2\beta+1)}},\quad
\widetilde\lambda_1\cos\varphi_1\approx N^{-\frac\beta{2\beta+1}}\approx 
\widetilde \lambda_{N-1}\cos\varphi_{N-1},\quad \prod_{j=2}^{N-2}\widetilde\lambda_j\approx N\, .
\eeq
In particular, $\widetilde\lambda_{N-1}^{-1}+\widetilde\lambda_1^{-1}=\mathcal O(1)$, and
\beq\label{MMM}
N^{1+\frac \beta{(\beta+1)(2\beta+1)}}|V|_p=\mathcal O\left(  |F(V)|_p\right)\, ,
\eeq
\beq\label{P2}
\frac{|\mathcal FV|_p}{|V|_p}\approx N^{\frac\beta{(\beta+1)(2\beta+1)}}
\quad\mbox{and}\quad
|FV|=|\mathcal F^NV|
\approx N^{-\frac{\beta}{2\beta+1}}|\mathcal F^{N-1}V|\, .
\eeq
\end{lemma}
\begin{proof}
We set $V_k:=\cF^kV$ for all $k=1,..., N$. 
Then $F(V)=V_N$ and each $V_k$ is given by  $d\varphi=\phi_k(r)$. It follows from e.g.~\eqref{DTdiff} (see e.g.~\cite{CM}) that
\begin{equation}\label{AA}
\forall 0\le k\le m\le N,\quad |V_m|_p=\int_{V_k} \left(\prod_{n=k}^{m-1} 
\widetilde \lambda_n\right)\, 
\, \cos\phi_k\, dr 
\, ,
\end{equation}
\beq\label{Eqbn}
\widetilde\lambda_n=1+\lambda_n\quad\mbox{and}\quad\lambda_n:=\frac{\tau_n}{\cos\phi_n}\left(\cK_n+\phi_n'\right)=\frac{2\tau_n\cK_n}{\cos\phi_n}+\frac {\frac{\tau_{n}}{\tau_{n-1}}}{
1+\frac 1{ \lambda_{n-1}   }}\, ,
\eeq
with $\cK_n:=\cK\circ\cF^n$ and $\tau_n:=\tau\circ\cF^n$. 
Note that~\eqref{AA} is a rewritting of~\eqref{AA'}.
Due to Item (4) of Proposition~\ref{prop3}, 
\beq\label{estibL}
\widetilde\lambda_0\approx 1+ \tau_0 \approx N^{\frac \beta{(\beta+1)(2\beta+1)}}\, .
\eeq
ensuring the first estimate in~\eqref{controllambda}. Furthermore it follows from Items (3), (4) and (5) of Proposition~\ref{prop3} that $\lambda_2\approx 1$, and that
\beq\label{bL+1}
\widetilde \lambda_1\approx 1+\frac{t_1 t_1^{-\beta-2}}{\cos\phi_1}+\frac{\mathcal O(1)}{1+\frac 1{\lambda_0}}\approx \frac{N^{-\frac\beta{2\beta+1}}}{\cos\phi_1}
\, .
\eeq
Analogously
\beq\label{estibL_last}
\widetilde\lambda_{N-1}\approx 1+\frac{t_{N-1}t_{N-1}^{-\beta-2}}{\cos\phi_{N-1}}+\mathcal O\left(\frac{t_{N-1}}{\tau_{N-2}}\right)\approx \frac {t_{N-1}^{-\beta-1}}{\cos\phi_{N-1}}\approx \frac{N^{-\frac\beta{2\beta+1}}}{\cos\phi_{N-1}}
\, ,
\eeq
since $\tau_{N-2}\approx t_{N-1}$ (where we used Proposition~\ref{prop3}).
Finally, to prove the last estimate of~\eqref{controllambda}, 
we will prove the existence of $z_+>0$ and $z_-<0$ such that $z_+-z_-=1$ and such that $\lambda_n=\frac{z_+}n+e_n$ for $n\le N_2$ and  $\lambda_n= \frac{z_-}{n-N}+e_n$ for $n>N_2$, with $\sum_{n=2}^{N-1}|e_n|<\infty$. This will imply that
\beq\label{bninterm}
\prod_{n=2}^{N_2}\widetilde\lambda_n\approx  \prod_{n=1}^{N_2}\left(1+ \frac{z_+}{n}\right)\approx N^{z_+}\quad\mbox{and}\quad \prod_{n=N_2+1}^{N-2}\widetilde\lambda_n\approx  \prod_{n=N_2+1}^{N-2}\left(1- \frac{z_-}{N-n}\right)\approx N^{-z_-}\, .
\eeq
Due to Items (5) and (6) of Proposition~\ref{prop3} for $n=2,...,N-2$,
\beq\label{controlinterm}
\frac{2\tau_n\cK_n}
{\cos\phi_n}=\frac{c (1+\cO(\varepsilon_{n,N}))}{A_n^2}\quad
\mbox{and}\quad \frac{\tau_{n}}{\tau_{n-1}}=1-\frac{a}{A_n} (1+\cO(\varepsilon_{n,N}))\, .
\eeq
with $c:=\frac{\beta(\beta+1)}{(2\beta+1)^2}$, $a:=\frac{2\beta}{2\beta+1}$,
$A_n=n$ if $n\le N/2$ and $A_n=n-N$ if $n> N/2$ and with $
\varepsilon_{n,N}=\frac{1}n+v_n^2
$ if $n\le N/2$ and $
\varepsilon_{n,N}=\frac{1}{N-n}+v_n^2
$ if $n\ge N/2$.

\begin{itemize}
\item \underline{Study of $A_n\lambda_n$.}
It follows from~\eqref{Eqbn} combined with~\eqref{controlinterm} that, 
for any $n=2,...,N-1$, 
\begin{align}
\label{Anbn}
A_n\lambda_n&=A_n\frac{2\tau_n\cK_n}{\cos\phi_n}+A_n\frac {\frac{\tau_{n}}{\tau_{n-1}}}{
1+\frac 1{ \lambda_{n-1}   }}\\
&=\frac c {A_n} \left(1+\cO(\varepsilon_{n,N})\right)+A_{n-1}\lambda_{n-1}\frac{1-\frac a{A_n}\left(1+\cO(\varepsilon_{n,N})\right)}{1+\frac{A_{n-1}\lambda_{n-1}}{A_{n-1}}}\frac 1{1-\frac 1{A_{n}}}
\, ,\label{bnrec}
\end{align}
since $A_{n-1}=A_{n}-1$.
Let us prove that
\beq\label{anbnbounded}
\sup_{n=2,...,N-1}|A_n|\lambda_{n}\le M_2<\infty\, .
\eeq
Observe that $\frac{|A_n\lambda_n|}{|A_{n-1}\lambda_{n-1}|}\le   \frac {c'}{|A_n\, A_{n-1}\lambda_{n-1}|} + \frac{1+\frac{c'}{|A_{n}|}}{1+\frac{|A_{n-1}|\lambda_{n-1}}{|A_{n-1}|}}$, for some uniform $c'>0$.\\
Fix some $C_0>1$ such that $C_0>c'+c'/C_0$, then if  $|A_{n-1}\lambda_{n-1}|>C_0$, 
\beq
\frac{|A_n\lambda_n|}{|A_{n-1}\lambda_{n-1}|}\le  \frac {c'}{|A_n|C_0} + \frac{1+\frac{c'}{|A_n|}}{1+\frac{C_0}{|A_{n-1}|}}=1-\frac{C_0-\frac{c'}{C_0}-c'}{|A_n|}+\mathcal O(A_n^{-2})=\mathcal O\left(  1\right)\, ,\eeq
this quantity is even smaller than 1 provided $|A_n|\ge n_0$ for some uniform  $n_0$.
Moreover, if $|A_{n-1}|\lambda_{n-1}\le C_0$, then $|A_n|\lambda_n\le c'+(c'+1)C_0$. Since $\lambda_{2}=\cO(1)$ and $A_{\lfloor N/2\rfloor}\sim A_{\lfloor N/2\rfloor+1}$, we have proved~\eqref{anbnbounded}.  
It follows hence from \eqref{bnrec} that, for all $n\ge 3$, the following holds true
\begin{align*}
A_n\lambda_n&=A_{n-1}\lambda_{n-1}+\frac {c-(A_{n-1}\lambda_{n-1})^2+A_{n-1}\lambda_{n-1}(1-a)}{A_n} +\cO\left(\frac{\varepsilon_{n,N}}{A_n}\right)\\
&=A_{n-1}\lambda_{n-1}-\frac {(A_{n-1}\lambda_{n-1}-z_+)(A_{n-1}\lambda_{n-1}-z_-)}{A_n} +\cO\left(\frac{\varepsilon_{n,N}}{A_n}\right)\, ,
\end{align*}
with $z_{\pm}:=\frac{1-a
\pm\sqrt{(1-a)^2+4c}}2=\frac 12\left(\frac 1{2\beta+1}\pm 1\right)$. In particular $z_+-z_-=1$. 
Let us write $(z_{n,+},z_{n,-})=(z_+,z_-)$ if $n\le N/2$ and 
$(z_{n,+},z_{n,-})=(z_-,z_+)$ if $n>N/2$.
\begin{align*}
A_n\lambda_n-z_{n,+}
&=(A_{n-1}\lambda_{n-1}-z_{n,+})\left(1-\frac {A_{n-1}\lambda_{n-1}
-z_{n,-}}{A_n}\right) +\cO\left(\frac{\varepsilon_{n,N}}{A_n}\right)\, .
\end{align*}
Let $n_1$ be such that $|z_{n,-}/A_n|<1/4$ and $M_2/A_n<1/4$ for all $n=n_1,...,N-n_1$.
Since $\frac{A_{n-1}\lambda_{n-1}}{A_n}\ge 0$ and $\frac{z_{n,-}}{A_n}\le 0$, it follows that, for all $n=n_1,...,N-n_1$,
\begin{align}
\nonumber|A_n\lambda_n-z_{n,+}|
&\le |A_{n-1}\lambda_{n-1}-z_{n,+}|\, \left(1-\frac {|z_{n,-}|}{|A_n|} \right)
 +\cO\left(\frac{\varepsilon_{n,N}}{A_n}\right)\\
&\le |A_{n-1}\lambda_{n-1}-z_{n,+}|\, e^{-\frac {|z_{n,-}|}{|A_n|}} 
 +\cO\left(\frac{\varepsilon_{n,N}}{A_n}\right)\, .\label{majoiter}
\end{align}
\item \underline{Contribution of the entering part.}
\eqref{majoiter} ensures that, for all $n=n_1,...,N_2-1$
\begin{align}
\nonumber&|n\lambda_n-z_{+}|
\le |n_1\lambda_{n_1}-z_+|\prod_{k=n_1+1}^n e^{-\frac {|z_-|}{k}} +\sum_{k=n_1+1}^{n}\cO\left(\frac{\varepsilon_{k,N}}{k}\right)\prod_{m=k}^n e^{-\frac {|z_-|}{m+1}}\\
&\quad\le C |n_1\lambda_{n_1}-z_+| n^{-|z_-|} +\sum_{k=2}^{n}\cO\left(\frac{\varepsilon_{k,N}}{k}\right) (k/n)^{|z_-|}= \cO\left( n^{-|z_-|}   +
v_n^2
\right)\, ,\label{MMMM1}
\end{align}
where we used the fact that
\begin{align*}
\sum_{k=2}^n \frac{\varepsilon_{k,N}}{k} (k/n)^{|z_-|}
   &= \sum_{k=2}^n \left(\frac{1}{k^2}+\frac{v_k^2}k\right) (k/n)^{|z_-|}
    =\mathcal O\left(  n^{-|z_-|}+v_n^2\right)\, ,
\end{align*}
since $\frac{v_k^2}k\le v_n\alpha_k$ and $\sum_{k=2}^{n}\alpha_k\le v_n$. 
Note that $e_{n}:=\lambda_n-\frac{z_{+}}n$ is uniformly normally summable over $n=2,...,N_2-1$ since $ \sum_{n=2}^{N_2-1}\frac{v_n^2}n=\mathcal O\left(\sum_{n=2}^{N_2-1}v_n\alpha_n\right)=\cO(v_n^2)=\cO(1)$.
It follows that
\beq
\prod_{n=2}^{N_2} 
\left(1+\lambda_n\right)=\prod_{n=2}^{N_2} 
\left(1+\frac{z_+
}{n}+e_n\right)
\approx e^{\sum_{n=2}^{N_2} 
\frac{z_+}
{n}+e_n}\approx N^{z_+} \, .
\label{prodentry}
\eeq
\item \underline{Contribution of the exiting part.}
It follows from \eqref{majoiter} that, for $n=N_2,...,N-n_1$
\begin{align}
\nonumber|(N-n)\lambda_n+z_{-}|
&\le |(N-N_2)\lambda_{N_2}+z_{-}|\prod_{k=N_2}^n e^{-\frac {z_+}{N-k}} +\sum_{k=N_2+1}^{n}\cO\left(\frac{\varepsilon_{k,N}}{N-k}\right)\prod_{m=k}^n e^{-\frac {z_+}{N-m}}\\
\nonumber&\le C \left(\frac{N-n}{N}\right)^{|z_+|} +\sum_{k=N_2+1}^{n}\cO\left(\frac{\varepsilon_{k,N}}{N-k}\right) \left(\frac{N-n}{N-k}\right)^{|z_+|}\\
&\le \cO\left( \left(\frac{N-n}{N}\right)^{|z_+|}    + \varepsilon_{n-N,N}\right)
\, ,\label{MMMM2}
\end{align}
where we used the fact that $\sum_{k=N_2}^{n}(N-k)^{-2-z_+}=\mathcal O\left(  (N-n)^{-1-z_+}\right)$ and 
\begin{align*}
\sum_{k=N_2}^n& \frac{v_k^2}{N-k}          \left(\frac{N-n}{N-k}\right)^{|z_+|} =\mathcal O\left(  (N-n)^{z_+}N^{-\frac{2\beta}{2\beta+1}}\sum_{k=N_2}^n(N-k)^{\frac{2\beta}{2\beta+1}-1-z_+}\right) \\
 &=\mathcal O\left( \left(\frac{N-n}{N}\right)^{\frac{2\beta}{2\beta+1} }\right)\approx v_n^2 \, .
\end{align*}
Thus $e_{n}:=\lambda_n+\frac{z_{-}}{N-n}$ is uniformly normally summable over $n=N_2,...,N-2$ since $N^{-z_+}\sum_{n=N_2+1}^N(N-n)^{z_+-1}\approx 1$ and
$\sum_{n=N_2+1}^{N-2} \frac{v_n^2}{N-n}\le v_{N_2+1}\sum_{n=N_2+1}^{N-1}\alpha_{n}\approx 1$, which implies that
\beq
\prod_{n=N_2}^{N-2}
\left(1+\lambda_n\right)=\prod_{n=N_2+1}^{N-2} 
\left(1-\frac{z_-}{N-n}+e_n\right)
\approx \exp\left(\mathcal O(1)-\sum_{n=N_2+1}^{N-2} 
\frac{z_-}{N-n}\right)\approx N^{-z_-}
\, .
\label{prodexit}
\eeq
\item \underline{Conclusion.} 
Due to \eqref{AA}, since our
estimates \eqref{estibL}, \eqref{bL+1}, \eqref{prodentry} and \eqref{prodexit} are uniform on $V_N$ and also in $W$ and $N$,
we obtain the last estimate of~\eqref{controllambda}.
\end{itemize}
It follows from Item~(3) of Proposition~\ref{prop3} that 
$\cos\varphi_1+\cos\varphi_{N-1}=\mathcal O\left(  N^{-\frac\beta{2\beta+1}}\right)$, thus $1=\mathcal O\left( \widetilde\lambda_1\right)$ 
and $1=\mathcal O\left( \widetilde\lambda_{N-1}\right)$. Thus~\eqref{MMM} follows from these two last estimates combined with~\eqref{AA'} and~\eqref{controllambda}. Analogously~\eqref{P2} follows from the two  following estimates:
$$
|\mathcal FV|_p=\int_{V}\widetilde\lambda_0\, |dx_0|_p\approx N^{\frac\beta{(\beta+1)(2\beta+1)}}|V|_p\, ,
$$
and
\begin{align*}
|\mathcal F^NV|&\approx |\mathcal F^NV|_p=\int_{\mathcal F^{N-1}(V)}
\widetilde \lambda_{N-1}|dx_{N-1}|_p\\
&\approx\int_{\mathcal F^{N-1}(V)}
\widetilde \lambda_{N-1}\cos\varphi_{N-1} dx_{N-1}
\approx N^{-\frac{\beta}{2\beta+1}}|\mathcal F^{N-1}V|\, .
\end{align*}
\end{proof}

\section{Distortion bounds}\label{distortion}
Recall that $\vartheta$ has been chosen such that $\vartheta>\beta+1$. 
 Next, let us discuss the regularity of homogeneous unstable curves for the induced map. We fix a constant $C_b$ and denote $\mathcal{W}^u_H$ as the collection of all homogeneous unstable curves 
 for the induced system $(F, M, \mu_M)$, a curve $W$ in $M_N$ being homogeneous if $\mathcal F(W)$ and $\mathcal F^{N-1}(W)$ are both contained in homogeneity strips. This appendix is devoted to the proof of the following result which states that curves in $\mathcal{W}^u_H$ have uniform distortion bounds.
 \begin{proposition}\label{distortionbound}
 There exists $C 
 > 0$ such that, for every homogeneous unstable curve $W \in \mathcal{W}^u_H$,
 assume $W$ is $\mathcal C^{1+\frac 1{\vartheta+1}}$, and
 for every $x,y\in F(W)$,
 $$
 |\ln J_W F^{-1}(x) - \ln J_W F^{-1}(y)| \leq C d_W(x, y)^{\frac 1{\vartheta+1}}\, .
 $$
 \end{proposition}

\begin{proof}
Let $W$ be an unstable curve of $M_N$ such that $\cF (W)$ and $\cF^{N-1}(W)$ are not cut by homogeneity lines.  
It will be useful to notice that this ensures that $\widetilde\lambda_n$ is uniformly bounded on $W$ (by a constant depending on $N$ and on the homogeneity strips containing $\cF W$ and $\cF^{N-1}W$). 
Let $x,y\in W$. The following estimates hold true with uniform constants independent on $W,x,y$. 
Keeping the notations of the two previous Appendices, we write
\begin{align*}
\left|\ln J_WF(x)-\ln J_WF(y)\right|
&=\left|\ln \left(\prod_{n=0}^{N-1}\widetilde\lambda_n(x)\right)-\ln \left(\prod_{n=0}^{N-1}\widetilde\lambda_n(y)\right)\right|\\
&\le \sum_{n=0}^{N-1}\left| \ln \left(\widetilde\lambda_n(x)\right)
-\ln \left(\widetilde\lambda_n(y)\right) \right|\, ,
\end{align*}
with
\beq\label{FFF0}
\widetilde\lambda_n=1+\lambda_n=1+\frac{\tau_n(\kappa_n+\phi'_n)}{\cos\phi_n}\, .
\eeq
We have to control
$$ 
\frac{\widetilde\lambda_n(x)-\widetilde\lambda_n(y)}{\widetilde\lambda_n(x)}=\frac{\lambda_n(x)-\lambda_n(y)}{\widetilde\lambda_n(x)}\, .
$$
It will be worthwhile to recall that
\beq\label{FFF00}
\kappa_n\le \phi'_n=\kappa_n+\frac{\cos\phi_n}{\tau_{n-1}+\frac {1}{\frac {1}{\cos\phi_{n-1}}(\kappa_{n-1}+\phi'_{n-1})}}\le \kappa_n+\frac 1{\tau_{n-1}}\, .
\eeq
We will set $d_e$ for the euclidean distance in $(r,\varphi)$-coordinates and 
$d_p(x_k,y_k)$ for the  $p$-length of the portion of the curve  $\mathcal F^k(W)$ between $x_k$ anf $y_k$. 
\begin{itemize}
\item \underline{A first estimate on the length of $FW$.}
Since $\mathcal F^{N-1}W$ lies in a homogeneity strip, let us say $\mathbf H_k$ with $k\ge N^{\frac\beta{\vartheta(2\beta+1)}}$, 
$$
|\phi_{N-1}(x)-\phi_{N-1}(y)|\le k^{-1-\vartheta}\le (k^{-\vartheta})^{\frac {1+\vartheta}\vartheta}\approx|\cos\phi_{N-1}|^{\frac {1+\vartheta}\vartheta}\le N^{-\frac{\beta(1+\vartheta)}{(2\beta+1)\vartheta}}
\, ,
$$
so
\begin{align*}
d_e(x_{N-1},y_{N-1})&=\mathcal O\left(\frac{|\phi_{N-1}(x)-\phi_{N-1}(y)|}{\phi'_{N-1}}\right)\\
&=\mathcal O\left( N^{-\frac{\beta(1+\vartheta)}{(2\beta+1)\vartheta}}N^{\frac {\beta(\beta+2)}{(\beta+1)(2\beta+1)}}\right)
=\mathcal O\left( N^{\frac{\beta(\vartheta-\beta-1)}{\vartheta(2\beta+1)(\beta+1)}}\right)\, ,
\end{align*}
where we used the fact that $\phi'_{N-1}\approx\kappa_{N-1}$, since $\kappa_{N-1}\le\phi'_{N-1}\le\kappa_{N-1}+\frac {\cos \phi_{N-1}}{\tau_{N-2}}$.
It follows then from~\eqref{P2} in Lemma~\ref{propcvPP} that
\begin{align}\label{estimdN}
d_e(x_{N},y_{N})\approx N^{-\frac\beta{2\beta+1}}d_e(x_{N-1},y_{N-1})
=\mathcal O\left( N^{-\frac{{\beta(\vartheta\beta+\beta+1)}}{\vartheta(2\beta+1)(\beta+1)}}\right)\, .
\end{align}
So 
\beq\label{estimdN2}
d_e(x_N,y_N)=\mathcal O\left( N^{-\frac{{\beta(\vartheta\beta+\beta+1)}}{(\vartheta+1)(2\beta+1)(\beta+1)}}(d_e(x_N,y_N))^{\frac 1{\vartheta+1}}\right)\, .
\eeq
We will use several times the fact that the exponent of $N$ in the above formula is strictly smaller than $-\frac{\beta^2}{(2\beta+1)(\beta+1)}$. 
\item \underline{Contribution of the initial point}. Observe that
$\tau_0\approx N^{\frac \beta{(2\beta+1)(\beta+1)}}$ whereas the other quantities appearing in the expression of $\widetilde\lambda_0$ are uniformly
$\frac 1{\vartheta+1}$-H\"older continuous.   
So the estimate will follow from the following computation
\begin{align*}
\frac{|\tau_{0}(x)-\tau_{0}(y)|}{\tau_{0}(x)}&=\mathcal O\left( t_{1}^{-1}(d_e(x,y)+d_e(x_{1},y_{1}))\right)\\
&= \mathcal O\left(t_1^{-1}\frac{
d_p(x_1,y_1)
}{\cos\varphi_1}\right)
= \mathcal O\left(t_1^{-1}
\frac{
d_p(x_N,y_N))
}{\cos\varphi_1}\prod_{j=1}^{N-1}\widetilde\lambda_j^{-1}\right)
\\
&= \mathcal O\left(t_1^{-1}N^{-1+\frac\beta{2\beta+1}}
\left(d_p(x_{N},y_{N})\right)\right)\,
\end{align*}
where we used Lemma~\ref{propdetails}. Thus
\begin{align*}
\frac{|\tau_{0}(x)-\tau_{0}(y)|}{\tau_{0}(x)}
&= \mathcal O\left(t_1^{-1}N^{-1}N^{\frac\beta{2\beta+1}}
d_e(x_{N},y_{N})\right)=\mathcal O\left( N^{-\frac{\beta^2+3\beta+1}{(\beta+1)(2\beta+1)}}d_e(x_{N},y_{N})\right)\, .
\end{align*}
Therefore, we have proved that
$$\frac{\widetilde\lambda_0(x)-\widetilde\lambda_0(y)}{\widetilde\lambda_0(x)}
=\mathcal O\left(N^{-\frac{\beta^2+3\beta+1}{(\beta+1)(2\beta+1)}}d_e(x_{N},y_{N}) \right).$$
\item \underline{Contribution of the first collision.}
Observe that $1=\mathcal O(\lambda_{1})$, so $\widetilde\lambda_1\approx\lambda_1$. Thus we have to estimate 
$\frac{\lambda_{1}(x)-\lambda_{1}(y)}{\lambda_{1}(x)}$.
Since $\mathcal F(W)$ is contained in a single homogeneity strip, it follows from~\eqref{controlcos} that
\begin{align}
&\frac{|\cos(\phi_{1}(x))-\cos(\phi_{1}(y))|}{\cos(\phi_{1}(x))}=\mathcal O\left(|\phi_{1}(x)-\phi_{1}(y)|^{\frac 1{\vartheta +1}}\right)\nonumber\\
&=\mathcal O\left(\phi'_{1}(x)d_e(x_{1},y_{1})\right)^{\frac 1{\vartheta +1}}\nonumber\\
&=\mathcal O\left(\frac{\cos\phi_{1} \widetilde\lambda_{1}}{\tau_{1}}d_e(x_{1},y_{1})\right)^{\frac 1{\vartheta +1}}\label{FFF1}\\
&=\mathcal O\left(\frac{ \widetilde\lambda_{1}}{\tau_{1}}d_p(x_{1},y_{1})\right)^{\frac 1{\vartheta +1}}\, .\nonumber
\end{align}
where we used also the fact that the slope of $\mathcal F(W)$ is $\phi'_{1}\approx\frac{\cos\phi_{1} \widetilde\lambda_{1}}{\tau_{1}}$
(this follows from~\eqref{FFF0} combined with~\eqref{FFF00} and $\lambda_1\approx\widetilde\lambda_1$). 
Now, it follows from Lemma~\ref{propdetails} that
\begin{align}
&\frac{|\cos(\phi_{1}(x))-\cos(\phi_{1}(y))|}{\cos(\phi_{1}(x))}=\mathcal O\left(\frac{ \widetilde\lambda_{1}}{\tau_{1}}\left(\prod_{j=1}^{N-1}\widetilde \lambda_{j}^{-1}\right)d_e(x_{N},y_{N})\right)^{\frac 1{\vartheta +1}}\nonumber\\
&=\mathcal O\left(\tau_1^{-1}\left(\prod_{j=2}^{N-1}\widetilde \lambda_{j}^{-1}\right)d_e(x_{N},y_{N})\right)^{\frac 1{\vartheta +1}}\nonumber\\
&=\mathcal O\left(N^{-1-\frac\beta{(\beta+1)(2\beta+1)}}d_e(x_{N},y_{N})\right)^{\frac 1{\vartheta +1}}\, .\label{FFF11}
\end{align}
Analogously
\begin{align}
\frac{|\kappa_{1}(x)-\kappa_{1}(y)|}{\kappa_{1}(x)}&=\mathcal O\left(\frac{t_{1}^{-\beta-3}}{t_{1}^{-\beta-2}}d_e(x_{1},y_{1})\right)\nonumber\\
&=\mathcal O\left(N^{-1+\frac {\beta^2}{(\beta+1)(2\beta+1)}}d_e(x_{N},y_{N})\right)\, .\label{FFF2}
\end{align}
using the fact that~\eqref{FFF1} is bounded by~\eqref{FFF11} combined with
$\cos\phi_{1} \widetilde\lambda_{1}\approx N^{-\frac{\beta}{2\beta+1}}$ 
(see~\eqref{controllambda}).
Furthermore
\begin{align*}
\frac{|\tau_{1}(x)-\tau_{1}(y)|}{\tau_{1}(x)}&=\mathcal O\left( t_{1}^{-1}(d_e(x_1,y_1)+d_e(x_{2},y_{2})\right)\, \\
&=\mathcal O\left( d_e(x_N,y_N)N^{-1+\frac{\beta}{2\beta+1}} \right)\, \\
\end{align*}
using the previous estimate of $d_e(x_1,y_1)/t_1$ together with
 \begin{align*}
 d_e(x_2,y_2)&=\mathcal O\left(d_e(x_{N},y_{N})N^{-1}/\cos\phi_2\right)\\
 &=\mathcal O\left(d_e(x_N,y_N)N^{-1+\frac{\beta}{2\beta+1}}  \right)\, .
 \end{align*}
Since $\kappa_{1}$ and $\phi'_{1}$ have same order, it remains to estimate $\frac{|\phi'_{1}(x)-\phi'_{1}(y)|}{t_1^{-\beta-2}}$, which follows from the following formula
$$
\phi'_1
=
\kappa_1
+\frac{\cos\phi_1}{\tau_0+\frac {1}{\frac {1}{\cos\varphi_0}(\kappa_0+\phi'_0)}}\, ,
$$
combined with the previous estimates and the uniform smoothness of $\kappa_0$ and $\phi_0$.
\item \underline{Contribution of the intermediate collisions}. Assume $n=2,...,N-2$.  
Recall that $\lambda_n\approx n^{-1}$. Thus
$$\frac{\widetilde\lambda_n(x)-\widetilde\lambda_n(y)}{\widetilde\lambda_n(x)}\sim\lambda_n(x)-\lambda_n(y)\, .$$
\begin{itemize}
\item \underline{First estimates}. 
\begin{align}\nonumber
|\tau_n(x)-\tau_n(y)|\frac{\kappa_n}{\cos\phi_n}
&=\mathcal O\left(\frac{t_n^{-\beta-2}}{\cos\phi_n}(d_e(x_{n},y_{n})+d_e(x_{n+1},y_{n+1})\right)\\
&=\mathcal O\left(\frac{t_n^{-\beta-2}}{(\cos\phi_n)^2}d_e(x_{N},y_{N})\prod_{j=n+1}^{N-1}\widetilde\lambda_j^{-1}\right)\, .\label{EEE1}
\end{align}
Let us prove the uniform summability of this quantity. We will treat separately the sum over $n<N/2$ (corresponding to the trajectory entering the cusp) and
the sum over $n\ge N/2$ (corresponding to the trajectory exiting the cusp).
If $n<N/2$, then it follows from Appendix~\ref{appendKBbd} (with the notations therein) and from Proposition~\ref{prop3} that 
$$\prod_{j=n}^{N-2}\widetilde\lambda_j^{-1}\approx  N^{z_-}(N/n)^{-z_+}= n^{z_+}N^{-1}$$
and that 
$$ \frac{t_n^{-\beta-2}}{(\cos\phi_n)^2}\approx n^{-\frac{3\beta+2}{2\beta+1}}N^{\frac{\beta^2}{(\beta+1)(2\beta+1)}}\, .$$
Hence, when $n<N/2$, 
\begin{align*}
\sum_{n=2}^{N/2}\frac{t_n^{-\beta-2}}{(\cos\phi_n)^2}\prod_{j=n+1}^{N-1}\widetilde\lambda_j^{-1}&=\mathcal O\left( N^{-1+
\frac{\beta^2}{(\beta+1)(2\beta+1)}} \sum_{n=2}^{N/2}n^{z_+-\frac{3\beta+2}{2\beta+1}}\right)\\
&=\mathcal O\left(N^{-1+\frac{\beta^2}{(\beta+1)(2\beta+1)}}\ln N\right)=\mathcal O(1)\, ,
\end{align*}
since $z_+=\frac{\beta+1}{2\beta+1}$. 
Analogously, when $n>N/2$, it follows from Appendix~\ref{appendKBbd} (with the notations therein) and from Proposition~\ref{prop3} that
$$\forall n>N/2,\quad \prod_{j=n+1}^{N-1}\widetilde\lambda_j^{-1}\approx (N-n)^{z_-}$$
and that
$$ \frac{t_n^{-\beta-2}}{(\cos\phi_n)^2}\approx (N-n)^{-\frac{3\beta+2}{2\beta+1}}N^{\frac{\beta^2}{(\beta+1)(2\beta+1)}}\, .$$
Here, contrary to the entering part, the exponent satisfies $-\frac{3\beta+2}{2\beta+1}+z_-=-2$, so the bare sum is only $\mathcal O\bigl(N^{\frac{\beta^2}{(\beta+1)(2\beta+1)}}\bigr)$ and does \emph{not} reduce to $\mathcal O(1)$. We therefore keep the factor $d_e(x_N,y_N)$ coming from~\eqref{EEE1} and apply~\eqref{estimdN2} (exactly as is done for the term $b_{j,2}$ below), which yields
$$
\sum_{n=N/2}^{N-2}\frac{t_n^{-\beta-2}}{(\cos\phi_n)^2}\,d_e(x_N,y_N)\prod_{j=n+1}^{N-1}\widetilde\lambda_j^{-1}
=\mathcal O\left( N^{\frac{\beta^2}{(\beta+1)(2\beta+1)}}d_e(x_N,y_N) \sum_{m=2}^{N/2}m^{-2}\right)
$$
$$
=\mathcal O\left( N^{\frac{\beta^2}{(\beta+1)(2\beta+1)}}d_e(x_N,y_N)\right)
=\mathcal O\left( N^{-\frac{\beta}{(\vartheta+1)(2\beta+1)(\beta+1)}}d_e(x_N,y_N)^{\frac1{\vartheta+1}}\right) ,
$$
which is uniformly bounded in $N$ (and summable over $N$), where $m=N-n$ and $z_-=-\frac\beta{2\beta+1}$.
This ends the control of the term with $|\tau_n(x)-\tau_n(y)|$. Analogously we have
\begin{align}\nonumber
|\kappa_n(x)-\kappa_n(y)|\frac{\tau_n}{\cos\phi_n}
&=\mathcal O\left(\frac{t_n/n}{\cos\phi_n}t_n^{-\beta-3}d_e(x_{n},y_{n})\right)\\
&=\mathcal O\left(\frac{t_n^{-\beta-2}/n}{(\cos\phi_n)^2}d_e(x_{N},y_{N})\prod_{j=n}^{N-1}\widetilde\lambda_j^{-1}\right)\, .\label{EEE2}
\end{align}
and
\begin{align}
|(\cos(\phi_{n}(x)))^{-1}&-(\cos(\phi_{n}(y)))^{-1}|\kappa_n\tau_n
=\mathcal O\left(\kappa_n\tau_n\frac{\phi'_{n}}{(\cos\phi_n)^2}d_e(x_{n},y_{n})\right)\nonumber\\
&=\mathcal O\left(\kappa_n\tau_n\frac{\phi'_{n}}{(\cos\phi_n)^3}\left(\prod_{j=n}^{N-1}\widetilde \lambda_{j}\right)d_e(x_{N},y_{N})\right)\nonumber\\
&=\mathcal O\left(\kappa_n\frac{1}{n(\cos\phi_n)^2}\left(\prod_{j=n}^{N-1}\widetilde \lambda_{j}\right)d_e(x_{N},y_{N})\right)\, ,\label{EEE3}
\end{align}
where we used the fact that $\frac{\tau_n}{\cos\phi_n}\phi'_n\approx\lambda_n\approx \frac 1n$.
The summability of~\eqref{EEE2} and~\eqref{EEE3} follows from the summability of~\eqref{EEE1}.
\item \underline{Control of $\lambda_n(x)-\lambda_n(y)$}. 
Let us write $\Delta_{x,y}(h):=h(x)-h(y)$. 
Recall that
$$
\lambda_n=\frac{2\tau_n\kappa_n}{\cos\varphi_n}+\frac{\tau_n/\tau_{n-1}}{1+\frac 1{\lambda_{n-1}}}\approx \frac 1n\, .
$$
Thus
$$
\Delta_{x,y}(\lambda_n)= a_n \Delta_{x,y}(\lambda_{n-1})+b_n
$$
with
$$a_n:=\frac{\tau_n(y)}{\tau_{n-1}(y)}\frac 1{(1+\lambda_{n-1}(x))(1+\lambda_{n-1}(y))}=\frac{\tau_n(y)}{\tau_{n-1}(y)}\frac 1{\widetilde\lambda_{n-1}(x)\widetilde\lambda_{n-1}(y)}\, ,$$
which is uniformly bounded 
and
$$
b_n:=\Delta_{x,y}\left(\frac{2\tau_n\kappa_n}{\cos\varphi_n}\right)+
\frac{\Delta_{x,y}(\tau_n/\tau_{n-1})}{1+\frac 1{\lambda_{n-1}}}
=b_{n,1}+b_{n,2}\, ,
$$
writing $b_{n,j}$ for the $j$-th term of the middle hand side of the above formula. 
Therefore, 
\beq\label{decompDeltalambda}
\Delta_{x,y} (\lambda_n)=(a_n...a_{2})\Delta_{x,y}(\lambda_1)+\sum_{j=2}^n(a_n...a_{j+1})b_j\, ,
\eeq
with
$$
a_n....a_{j+1}\approx\frac{\tau_n}{\tau_j}\prod_{i=j}^{n-1}\widetilde\lambda_i^{-2}\, .
$$
In particular, the sum over $n$ of the first term of the right hand side of~\eqref{decompDeltalambda} is controlled via the following estimate
using Proposition~\ref{prop3} and~\eqref{bninterm}:
\begin{align*}
\sum_{n=2}^{N}(a_n...a_{2})&=\mathcal O\left( \sum_{n=2}^{N/2}\frac{t_n/n}{t_1}n^{-2z_+}
+ \sum_{n=N/2}^{N}\frac{t_n/(N-n)}{t_1}N^{-2z_+}(N/(N-n))^{2z_-} \right)\\
&=\mathcal O\left( \sum_{n=2}^{N/2}n^{\frac 1{2\beta+1}-1-2z_+}+N^{-2}\sum_{n=N/2}^{N}(N-n)^{\frac 1{2\beta+1}-1-2z_-}(N-n)^{-2z_+}\right)\\
&=\mathcal O(1)\, ,
\end{align*}
where we used also the fact that $z_+-z_-=1$ and that $2z_+-1=\frac 1{2\beta+1}$.\\
It remains to control the sum over $n$ of the second term of the right hand side of~\eqref{decompDeltalambda}.\\
Assume that $2\le j\le n\le N/2$. Then
\begin{align*}
b_{j,2}&=\mathcal O(1/j
)\Delta_{x,y}(\tau_j/\tau_{j-1})
=\mathcal O(1/j
) \Delta_{x,y}\left(\frac{t_{j+1}^{-\beta}+t_{j}^{-\beta}}{t_{j}^{-\beta}+t_{j-1}^{-\beta}}\frac{\cos\phi_{j-1}}{\cos\phi_j}\right)\\
&=\mathcal O\left(
\frac 1j
\left(\frac{|\Delta_{x,y}(t_j)|}{t_j}+\frac{|\Delta_{x,y}(\phi_{j-1})|}{\cos\phi_j}\right)\right)\\
&=\mathcal O\left(  \frac 1j
\left(\frac 1{t_j}+\frac{\phi'_j}{\cos\phi_j}\right)\right)(d_e(x_j,y_j)+d_e(x_{j-1},y_{j-1}))\\
&=\mathcal O\left(  \frac 1j
\left(\frac 1{t_j}+\frac{t_j^{-\beta-2}}{\cos\phi_j}\right)\right)\frac{d_p(x_j,y_j)}{\cos\phi_j}=\mathcal O\left(  \frac 1{j
t_j\cos\phi_j}\right)d_p(x_j,y_j)\, ,
\end{align*}
where we used the fact that $\cos\phi_j\approx \gamma_j\approx jt_j^{-\beta-1}$  (see Items~(3) and ~(4) of Proposition~\ref{prop3}). 
Thus
\begin{align*}
& (a_{n}...a_{j+1})b_{j,2}=\mathcal O\left(\frac 1{jt_j\cos\phi_j}
d_e(x_N,y_N)\frac{\tau_n}{\tau_j}\prod_{i'=j}^{n-1}\widetilde\lambda_{i'}^{-2}\prod_{i=j}^{N-1}\widetilde\lambda_i^{-1}\right)\\
&=\mathcal O\left(\frac 1{jt_j(j/N)^{\frac\beta{2\beta+1}}}
d_e(x_N,y_N)\frac{t_n/n}{t_j/j}(j/n)^{2z_+}(j^{z_+}/N)\right)\\
&=\mathcal O\left(\frac 1{j(jN^{\frac\beta{\beta+1}})^{\frac 1{2\beta+1}}(j/N)^{\frac\beta{2\beta+1}}}
d_e(x_N,y_N)(n/j)^{\frac 1{2\beta+1}-1-2z_+}j^{z_+}/N\right)\, .
\end{align*}
Recall that $z_+=\frac{\beta+1}{2\beta+1}$. 
Summing over $n$, we obtain
\begin{align*}
\sum_{n=j}^{N/2}& (a_{n}...a_{j+1})b_{j,2}=\mathcal O\left(N^{-1+\frac{\beta^2}{(\beta+1)(2\beta+1)}}\frac 1{j^{1+\frac{\beta+1}{2\beta+1}}}
d_e(x_N,y_N)j^{3z_+-\frac 1{2\beta+1}+1}j^{\frac 1{2\beta+1}-2z_+}\right)\\
&=\mathcal O\left(N^{-1+\frac{\beta^2}{(\beta+1)(2\beta+1)}}
d_e(x_N,y_N)\right)
\, .
\end{align*}
Summing over $j$ and using~\eqref{estimdN2}, we obtain
\begin{align*}
\sum_{j=2}^{N/2}\sum_{n=j}^{N/2}& (a_{n}...a_{j+1})b_{j,2}=\mathcal O\left(N^{\frac{\beta^2}{(2\beta+1)(\beta+1)}}\, d_e(x_N,y_N)  \right)\\
&=\mathcal O\left(N^{-\frac{\beta}{(\vartheta+1)(2\beta+1)(\beta+1)}}(d_e(x_{N},y_{N})^{\frac 1{\vartheta+1}}\right)\, .
\end{align*}
The control of the sums $\sum_{n=N/2}^N\sum_{j=2}^{N/2}$ and 
$\sum_{n=N/2}^N\sum_{j=N/2}^{n}$ go the same way up to replace $j$ and $n$ by respectively $\min(j,N-j)$ and $\min(n,N-n)$, and if, $j>N/2$, $N^{-1}j^{z_+}$ by $N^{z_-}(N-j)^{-z_-}$ and $(j/n)^{z_+}$ by $((N-j)/(N-n))^{-z_-}$, recalling that $z_+-z_-=1$.
We deduce the analogous estimate for $b_{n,1}$ using~\eqref{EEE2},~\eqref{EEE3} with~\eqref{EEE1}.
\end{itemize}

\item \underline{Contribution of the last collision}. Observe that $1=\mathcal O(\widetilde\lambda_{N-1})$. Thus we have to estimate 
$\frac{\lambda_{N-1}(x)-\lambda_{N-1}(y)}{\lambda_{N-1}(x)}$. 
Recall from~\eqref{Eqbn} that
$$ 
\lambda_{N-1}=2\frac{\tau_{N-1}\kappa_{N-1}}{\cos\phi_{N-1}
}+\frac{\frac{\tau_{N-1}}{\tau_{N-2}}}{1+\frac 1{\lambda_{N-2}}}
\, .
$$
Thus the estimates of $\lambda_{N-1}(x)-\lambda_{N-1}(y)$ will follow from the previous estimates (at order $N-2$) combined with the following estimates at order $N-1$.\\
Since $\mathcal F^{N-1}W$ is contained in a homogeneity strip, it follows from~\eqref{controlcos} that
\begin{align*}
&\frac{|\cos(\phi_{N-1}(x))-\cos(\phi_{N-1}(y))|}{\cos(\phi_{N-1}(x))}=\mathcal O\left(|\phi_{N-1}(x)-\phi_{N-1}(y)|^{\frac 1{\vartheta +1}}\right)\\
&=\mathcal O\left(\phi'_{N-1}d_e(x_{N-1},y_{N-1})\right)^{\frac 1{\vartheta +1}}\\
&=\mathcal O\left(\phi'_{N-1}\widetilde \lambda_{N-1}^{-1}d_e(x_{N},y_{N})/\cos\phi_{N-1}\right)^{\frac 1{\vartheta +1}}\\
&=\mathcal O\left(N^{-\frac{\beta}{(2\beta+1)(\beta+1)}}d_e(x_{N},y_{N})\right)^{\frac 1{\vartheta +1}}
\end{align*}
where we used the fact that the slope of $\mathcal F^{N-1}(W)$ is $\phi'_{N-1}\approx\kappa_{N-1}\approx N^{-\frac{\beta(\beta+2)}{(2\beta+1)(\beta+1)}}$
and that $(\widetilde\lambda_{N-1}\cos\phi_{N-1})^{-1}\approx N^{\frac{\beta}{2\beta+1}}$,
and Item~(4) of Proposition~\ref{prop3}. 
Furthermore
\begin{align}
&\frac{|\kappa_{N-1}(x)-\kappa_{N-1}(y)|}{\kappa_{N-1}(x)}=\mathcal O\left(\frac{t_{N-1}^{-\beta-3}}{t_{N-1}^{-\beta-2}}d_e(x_{N-1},y_{N-1})\right)\nonumber\\
&=\mathcal O\left(t_{N-1}^{-1}N^{\frac\beta{2\beta+1}}d_e(x_{N},y_{N})\right)=\mathcal O\left(N^{\frac\beta{2\beta+1}\left(1-\frac 1{\beta+1}\right)}d_e(x_{N},y_{N})\right)\nonumber\\
&=\mathcal O\left(N^{\frac{\beta^2}{(2\beta+1)(\beta+1)}}d_e(x_{N},y_{N})\right)\, .\label{estimkappafinal}
\end{align}
This combined with~\eqref{estimdN2} leads to
\begin{align*}
\frac{|\kappa_{N-1}(x)-\kappa_{N-1}(y)|}{\kappa_{N-1}(x)}
&=\mathcal O\left(N^{\frac{(\vartheta+1)\beta^2}{(\vartheta+1)(2\beta+1)(\beta+1)}}d_e(x_{N},y_{N})^{\frac 1{\vartheta+1}}N^{-\frac{{\beta(\vartheta\beta+\beta+1)}}{(\vartheta+1)(2\beta+1)(\beta+1)}}\right)\\
&=\mathcal O\left(N^{-\frac{\beta}{(\vartheta+1)(2\beta+1)(\beta+1)}}(d_e(x_{N},y_{N})^{\frac 1{\vartheta+1}}\right)\, ,
\end{align*}
and, analogously,
\begin{align*}
\frac{|\tau_{N-1}(x)-\tau_{N-1}(y)|}{\tau_{N-1}(x)}&=\mathcal O\left( t_{N-1}^{-1}(d_e(x_N,y_N)+d_e(x_{N-1},y_{N-1})\right)\\
&=\mathcal O\left(N^{-\frac{\beta}{(\vartheta+1)(2\beta+1)(\beta+1)}}(d_e(x_{N},y_{N})^{\frac 1{\vartheta+1}}\right)\, .
\end{align*}
\end{itemize}
\end{proof}

\medskip\noindent\textbf{Acknowledgements}.
M.~Lenci is partially supported by the PRIN Grant 2022NTKXCX of the Ministry of University and Research (MUR), Italy. This research is also part of his activity within the Gruppo Nazionale di Fisica Matematica, INdAM, Italy.
F.~P\`ene is supported by the ANR projects GALS (ANR-23-CE40-0001) and RAWABRANCH (ANR-23-CE40-0008), and conduced this work within the framework of the Henri Lebesgue Center ANR-11-LABX-0020-01.  H.-K.~Zhang is partially supported by the National Science Foundation (DMS-2220211), as well as Simons
Foundation Collaboration Grants for Mathematicians (706383).
We thank the CIRM in Marseille-Luminy, the Scuola Normale Superiore di Pisa, the ESI in Wien, the Erd\"os Center in Budapest, the universities of Amherst (UMass), Bologna and Brest for their hospitality.

\bigskip\noindent\textbf{AI Declaration}.
No AI was used in the preparation of this manuscript.

\end{document}